\documentclass[11pt,onecolumn]{IEEEtran}
\usepackage{graphicx}
\usepackage{amsmath}
\usepackage{mathrsfs}
\usepackage{mathtools}
\usepackage{amssymb}
\usepackage{float}
\usepackage{url}
\usepackage{verbatim}
\usepackage{enumerate}
\usepackage{layout}
\usepackage{bm}

\usepackage[bottom=0.8in]{geometry}

\usepackage[colorlinks,bookmarks]{hyperref}
\hypersetup{
  colorlinks   = true,
  urlcolor     = blue,
  linkcolor    = blue,
  citecolor    = magenta
}

\newtheorem{proposition}{Proposition}
\newtheorem{corollary}{Corollary}

\newtheorem{definition}{Definition}
\newtheorem{example}{Example}
\newtheorem{remark}{Remark}

\newcommand{\naturals}{\ensuremath{\mathbb{N}}}

\newcommand{\Reals}{\ensuremath{\mathbb{R}}}

\newcommand{\tr}{\mathrm{tr}}

\newcommand{\set}{\ensuremath{\mathcal}}

\newcommand{\OneTo}[1]{[#1]}

\newcommand{\card}[1]{|#1|}
\newcommand{\bigcard}[1]{\bigl|#1\bigr|}

\RequirePackage{bbm}

\DeclareMathOperator{\Vertex}{\mathsf{V}}
\DeclareMathOperator{\Edge}{\mathsf{E}}
\DeclareMathOperator{\Adjacency}{\mathbf{A}}
\DeclareMathOperator{\AllOne}{\mathbf{J}}
\DeclareMathOperator{\Identity}{\mathbf{I}}

\DeclareMathOperator{\Kneser}{\mathsf{K}}
\DeclareMathOperator{\Complete}{\mathsf{K}}

\DeclareMathOperator{\Path}{\mathsf{P}}
\DeclareMathOperator{\Cycle}{\mathsf{C}}
\DeclareMathOperator{\SRG}{\mathsf{srg}}
\DeclareMathOperator{\Clique}{\omega}
\DeclareMathOperator{\Chromatic}{\chi}
\newcommand{\Gr}[1]{\mathsf{#1}}
\newcommand{\CGr}[1]{\overline{\mathsf{#1}}}
\newcommand{\V}[1]{\Vertex(#1)}
\newcommand{\E}[1]{\Edge(#1)}
\newcommand{\A}{\Adjacency}
\newcommand{\J}[1]{\AllOne_{#1}}
\newcommand{\I}[1]{\Identity_{#1}}
\newcommand{\AG}[1]{\Adjacency(#1)}

\newcommand{\indnum}[1]{\alpha(#1)}
\newcommand{\clnum}[1]{\Clique(#1)}
\newcommand{\chrnum}[1]{\Chromatic(#1)}

\newcommand{\Eigval}[2]{\lambda_{#1}(#2)}
\newcommand{\CoG}[1]{\Complete_{#1}}
\newcommand{\CoBG}[2]{\Complete_{#1,#2}}
\newcommand{\KG}[2]{\Kneser(#1,#2)}

\newcommand{\PathG}[1]{\Path_{#1}}
\newcommand{\CG}[1]{\Cycle_{#1}}

\begin{document}
\thispagestyle{empty}
\setcounter{page}{1}
\setlength{\baselineskip}{1.15\baselineskip}

\title{\huge{Observations on the Lov\'{a}sz $\theta$-Function,
Graph Capacity, Eigenvalues, and Strong Products}\\[0.2cm]}
\author{Igal Sason\\[0.2cm]
{\em \normalsize{Dedicated to my friend and former teacher, Professor Emeritus
Abraham (Avi) Berman, in the occasion of his eightieth birthday\\[0.2cm]
{\bf{Citation}}:
{\em I. Sason, ``Observations on the Lov\'{a}sz $\theta$-function, graph capacity,
eigenvalues, and strong products,'' Entropy, vol.~25, no.~1, paper~104, pp.~1--41, January 2023.
\newline DOI: \url{https://doi.org/10.3390/e25010104}.}}}
\thanks{
I. Sason is with the Viterbi Faculty of Electrical and Computer Engineering,
and the Department of Mathematics at the Technion - Israel
Institute of Technology, Haifa 3200003, Israel (e-mail: eeigal@technion.ac.il).}}

\maketitle
\thispagestyle{empty}

\vspace*{-0.6cm}
\begin{small}
\begin{abstract}
This paper provides new observations on the Lov\'{a}sz $\theta$-function
of graphs. These include a simple closed-form expression of that function
for all strongly regular graphs, together with upper and lower bounds on that
function for all regular graphs. These bounds are expressed in terms of the
second-largest and smallest eigenvalues of the adjacency matrix of
the regular graph, together with sufficient conditions for equalities (the
upper bound is due to Lov\'{a}sz, followed by a new sufficient condition
for its tightness). These results are shown to be useful in many ways,
leading to the determination of the exact value of the Shannon capacity
of various graphs, eigenvalue inequalities, and bounds on the clique
and chromatic numbers of graphs. Since the Lov\'{a}sz $\theta$-function
factorizes for the strong product of graphs, the results are also
particularly useful for parameters of strong products or strong powers of
graphs. Bounds on the smallest and second-largest eigenvalues of strong
products of regular graphs are consequently derived, expressed as functions
of the Lov\'{a}sz $\theta$-function (or the smallest eigenvalue) of each
factor. The resulting lower bound on the second-largest eigenvalue of a
$k$-fold strong power of a regular graph is compared to the Alon--Boppana
bound; under a certain condition, the new bound is superior in its exponential
growth rate (in $k$). Lower bounds on the chromatic number of strong
products of graphs are expressed in terms of the order and the Lov\'{a}sz
$\theta$-function of each factor. The utility of these bounds is exemplified,
leading in some cases to an exact determination of the chromatic numbers
of strong products or strong powers of graphs. The present research paper
is aimed to have tutorial value as well.
\end{abstract}
\vspace*{-0.2cm}
\end{small}

\section{Introduction}
\label{section: introduction}

The notion of the graph capacity in Shannon's problem of zero-error communication \cite{Shannon56}
had a significant impact on the development of information theory and graph theory, including the
introduction of perfect graphs by C. Berge \cite{Berge73},
strong graph products (or powers) (\cite{Alon02,Sabidussi60}),
the introduction of the Lov\'{a}sz $\theta$-function of a graph as a computable upper bound on its
Shannon capacity (\cite{Knuth94,Lovasz79_IT}), the rank-bound by Haemers \cite{Haemers79}, and
other important follow-up works that are surveyed, e.g., in \cite{Alon02,Alon19,Jurkiewicz14,KornerO98}.

In graph theory, there are four central sorts of graph products, each with its own applications
and theoretical interpretations. The reader is referred to the excellent handbook \cite{HammackIK11},
which presents the rich and fertile field of graph products. Strong product of graphs is one of the
most extensively studied sorts of graph products, and there exists a polynomial-time algorithm that
finds the unique prime factorization of any connected graph with that type of multiplication
\cite{FeigenbaumS92}. Strong powers of graphs are also fundamental in information theory.
Their information-theoretic significance stems from the notion of the Shannon capacity of graphs
for error-free communication \cite{Shannon56}, and the Witsenhausen rate \cite{Witsenhausen76}
in the zero-error source coding problem with perfect side information at the receiver. Properties
of strong products and strong powers of graphs, and bounds on their independence numbers and
chromatic numbers have been extensively studied, e.g., in \cite{Acin17,Alon02,Alon19,AlonL06,
BohmanH03,BohmanHN09,BohmanHN13,Borowiecki72,BrimkovCCL04,Brimkov07,ErdosMT71,EsperetW22,
FeigenbaumS92,GeethaS18,Godsil16,GuoW90,Haemers79,Hales73,HammackIK11,HammackIK16,HuTS18,
Jurkiewicz14,Klavzar93,Klavzar96,KlavzarM94,LixinFH15,Lovasz79_IT,McElieceP71,Sabidussi60,
Shannon56,SonnemannK74,Vesztergombi79,Vesztergombi81,Witsenhausen76,XuR13}.

The present work continues the above paths of research. It provides some new observations on
the Lov\'{a}sz $\theta$-function of regular graphs, calculation of the Shannon capacity of some
strongly regular graphs, bounds on eigenvalues of graphs (in particular, the
second-largest and smallest eigenvalues of the adjacency matrix, which play a key role in spectral
graph theory), bounds related to Ramanujan graphs, and strong products of graphs.
The analysis in the present work mainly relies on the notion of the Lov\'{a}sz $\theta$-function
of graphs \cite{Lovasz79_IT}.
The paper includes a thorough review of the backgrounds relevant to this work with suitable
references or explanations, which also serve to motivate the presentation of the results in
this work and to put them into perspective.
The presentation in this research paper is consequently aimed to have tutorial value as well.

\vspace*{0.1cm}
The results obtained in this work are outlined as follows:
\begin{enumerate}[(1)]
\item
A known upper bound on the Lov\'{a}sz $\theta$-function of a regular graph is expressed
in terms of the smallest eigenvalue of its adjacency matrix \cite{Lovasz79_IT}.
A key result in this work provides a lower bound on the Lov\'{a}sz $\theta$-function of
a regular graph, which is expressed in terms of the second-largest eigenvalue of its
adjacency matrix. New sufficient conditions for attaining these bounds are specified
(Proposition~\ref{prop1: bounds on theta}).
\item
A simple and closed-form expression of the Lov\'{a}sz $\theta$-function is derived for all strongly
regular graphs (Corollary~\ref{cor4: Lovasz number for srg}).
\item
Eigenvalue inequalities are derived, which relate the smallest and second-largest eigenvalues of a regular
graph. They hold with equality if and only if the graph is strongly regular
(Corollaries~\ref{cor0: ineq. for graph eigvals} and~\ref{cor1: ineq. for graph eigvals}).
\item
The Shannon capacity of several strongly regular graphs is determined
(Section~\ref{subsection: graph capacities - s.r.g}).
\item
Bounds on parameters of regular graphs, and in particular of Ramanujan graphs, are derived
(Corollaries~\ref{cor1: clique and chromatic nums}--\ref{cor3: clique number - Ramanujan}).
\item
Bounds on the smallest and the second-largest eigenvalues of strong products of regular graphs
are derived, which are expressed in terms of calculable parameters of its factors
(Proposition~\ref{prop: bounds on eigvals - strong products}).
\item
A new lower bound on the second-largest eigenvalue of a $k$-fold strong power
of a regular graph is compared to the Alon--Boppana bound. Under a certain condition,
the former bound shows an improvement in its exponential growth rate
as a function of $k$ (Section~\ref{subsection: eigenvalues of strong products}).
\item
Every non-complete and non-empty connected regular graph, whose
Lov\'{a}sz $\theta$-function is below a certain value, is proved to have the property that
almost all its strong powers are highly non-Ramanujan (Proposition~\ref{prop.: not Ramanujan}).
\item
Lower bounds on the chromatic number of strong products of graphs are
expressed in terms of the order and Lov\'{a}sz $\theta$-function of each
factor (Proposition~\ref{prop. chromatic numbers}). Their utility
is exemplified, while also leading to exact chromatic numbers in some cases.
\end{enumerate}

The paper is structured as follows: Section~\ref{section: preliminaries}
provides notation and a thorough review of the backgrounds relevant to this work.
Section~\ref{section: results} provides the results of this work, followed by
examples and discussions. It is composed of five subsections
that address issues related to the Lov\'{a}sz $\theta$-function,
Shannon capacity of graphs, Ramanujan graphs, the second-largest
and smallest eigenvalues of strong products or strong powers of graphs,
and the chromatic numbers of such graphs. Section~\ref{section: proofs}
proves the results in Section~\ref{section: results}.

\section{Preliminaries}
\label{section: preliminaries}

This section provides essential notation and preliminaries
for this paper. The following standard notation in set theory is used:
$\naturals = \{1, 2, \ldots\}$ is the set of natural numbers,
$\Reals$ is the set of real numbers, and
$\OneTo{n} \triangleq \{1, \ldots, n\}$ with $n \in \naturals$.
The cardinality of a set $\set{A}$ is a measure of its number
of elements; it is denoted by $\card{\set{A}}$, and (by definition)
it is equal to its number of elements if $\set{A}$ is a finite set.

Let $\Gr{G}$ be a graph, and let $\V{\Gr{G}}$ and $\E{\Gr{G}}$
denote, respectively, the sets of vertices and edges in $\Gr{G}$.
The {\em order} and {\em size} of a graph $\Gr{G}$ are
defined to be $\card{\V{\Gr{G}}}$ and $\card{\E{\Gr{G}}}$, respectively.
A graph $\Gr{G}$ is said to be {\em finite} if $\V{\Gr{G}}$ is a finite set.
A pair of vertices are {\em adjacent} in a graph $\Gr{G}$ if these
two vertices are the endpoints of an edge $e \in \E{\Gr{G}}$.

A graph is called {\em simple} if it has no loops (i.e., it has no edge with
identical endpoints), and if there are no multiple edges between
any pair of adjacent vertices. A graph $\Gr{G}$ is said to be
{\em undirected} if its edges have no directions; otherwise, it is
a directed graph, a.k.a. a {\em digraph}. Throughout this paper, it is
assumed that the graphs under consideration are finite, undirected, and simple.

An {\em empty graph} is an edgeless graph.
A graph $\Gr{G}$ is said to be an {\em $r$-partite graph} if its vertex set $\V{\Gr{G}}$
is a disjoint union of $r$ subsets such that every pair of vertices that are
elements of an identical subset are non-adjacent.
If $r=2$, then it is a {\em bipartite graph}.

A {\em walk} in a graph $\Gr{G}$ is a sequence of its vertices such that
(by definition) every pair of consecutive vertices are adjacent.
A {\em path} in a graph is a walk with no repeated vertices (in other words,
a path is a walk along the vertices of $\Gr{G}$ such that no vertex can be
visited twice). The {\em length of a path} is defined as
its number of edges. Hence, $\Path = [v_1, \ldots, v_\ell]$ is a path
in a graph $\Gr{G}$ if $\{v_i, v_{i+1} \} \in \E{\Gr{G}}$ for all
$i \in \OneTo{\ell-1}$, and all the vertices in the sequence
$\{v_i\}_{i=1}^{\ell}$ are distinct; the endpoints of the path $\Path$
are $v_1$ and $v_\ell$, and its length is equal to $\ell-1$. A {\em cycle}
$\Cycle$ in a graph $\Gr{G}$ is obtained by adding an edge to a path
$\Path$ such that it gives a closed walk (i.e., a walk whose endpoints
are identical). The cycle $\Cycle = [v_1, \ldots, v_\ell, v_1]$ is of length
$\ell$, which is obtained by adding the edge $e = \{v_\ell, v_1\} \in \E{\Gr{G}}$
to the above $(\ell-1)$--length path $\Path$; the two identical endpoints of the cycle
$\Cycle$ are the vertex $v_1$.

A graph is said to be {\em connected} if every two distinct vertices in that
graph are connected by a path (otherwise, it is a {\em disconnected graph}).
The {\em distance} between a pair of distinct and connected vertices in a graph is defined
to be the length of the shortest path whose two endpoints are the given pair of vertices.
The distance between two disconnected vertices in a graph is defined to be infinite,
and the distance between a vertex and itself is set to be zero. The {\em diameter} of
a graph is the maximum distance between any pair of vertices in that graph.
A connected and finite graph has a finite diameter, and the diameter of a
disconnected graph is (by definition) infinite.

The following standard notation in graph theory is used:
\begin{enumerate}[(a)]
\item $\CoG{n}$ denotes the {\em complete graph} on $n \in \naturals$ vertices,
where every pair of distinct vertices are adjacent; hence,
$\CoG{1}$ is an empty graph with a single vertex.
\item $\CoBG{n}{m}$ denotes the {\em complete bipartite graph}, which is a bipartite graph
consisting of a vertex set that is a disjoint union of two finite sets $\set{V}_1$
and $\set{V}_2$ of cardinalities $n$ and $m$, respectively, and a set of edges that are all
the possible connections of a vertex in $\set{V}_1$ and a vertex in $\set{V}_2$.
\item $\PathG{n}$ denotes an {\em $(n-1)$--length path} with $n \in \naturals$, which is a graph
on $n$ vertices that forms a path of length $n-1$; in particular, $\PathG{1}=\CoG{1}$.
\item $\CG{n}$ denotes an {\em $n$-length cycle}, which is a graph on $n \geq 3$ vertices
that forms a cycle of length $n$.
\item $\KG{m}{r}$ denotes the {\em Kneser graph} with integers $1 \leq r \leq m$.
It has $n = \binom{m}{r}$ vertices, represented by all $r$-subsets of $\OneTo{m}$.
Two vertices are adjacent in that graph if they are represented by disjoint $r$-subsets.
The graph $\KG{m}{r}$, provided that it has more than one vertex, is a connected graph
if and only if either $m > 2r$ or $(m,r)=(2,1)$.
\end{enumerate}

A {\em subgraph} is a graph that exists within another graph. More formally,
$\Gr{F}$ is a subgraph of a graph $\Gr{G}$ if $\V{\Gr{F}} \subseteq \V{\Gr{G}}$
and $\E{\Gr{F}} \subseteq \E{\Gr{G}}$. If $\Gr{F}$ is a subgraph of $\Gr{G}$,
we write $\Gr{F} \subseteq \Gr{G}$. A {\em spanning subgraph} is obtained by
edge deletions from the original graph, while its vertex set is left unchanged.
An {\em induced subgraph} is obtained by removing vertices from the original graph,
followed by the deletion of the edges that are adjacent to these removed vertices.

A {\em clique} in a graph $\Gr{G}$ is a subset of pairwise adjacent vertices
in $\Gr{G}$ (in other words, the induced subgraph of $\Gr{G}$ on that subset
is a complete subgraph). The maximum size of a clique in $\Gr{G}$ is
called the {\em clique number} of $\Gr{G}$, and it is denoted by $\clnum{\Gr{G}}$.
Similarly, an {\em independent set} (a.k.a. {\em coclique}) in $\Gr{G}$ is a
subset of pairwise non-adjacent vertices. The maximum size of an independent
set is called the {\em independence number} of $\Gr{G}$, and it is denoted by
$\indnum{\Gr{G}}$.
A {\em proper vertex coloring} of $\Gr{G}$ is a coloring of its vertices such
that no pair of adjacent vertices are assigned the same color. The smallest
required number of colors for a proper coloring of the vertices in $\Gr{G}$
is called the {\em chromatic number} of $\Gr{G}$, and it is denoted by
$\chrnum{\Gr{G}}$.

Let $\Gr{G} = (\V{\Gr{G}}, \E{\Gr{G}})$ be a finite, undirected, and simple graph
of order $\card{\V{\Gr{G}}} = n$. Define the {\em adjacency matrix}
$\A = \AG{\Gr{G}}$ of the graph $\Gr{G}$ to be an $n$-times-$n$
symmetric matrix such that $\A = (a_{i,j})_{1 \leq i,j \leq n}$ with
$a_{i,j} = 1$ if $\{i,j\} \in \E{\Gr{G}}$, and $a_{i,j} = 0$ otherwise
(hence, the diagonal elements of $\A$ are in particular zeros).
Let the eigenvalues of $\A$ (a.k.a. the eigenvalues of $\Gr{G}$)
be given in decreasing order by
\begin{IEEEeqnarray}{rCl}
\label{eq:08.11.22a1}
\Eigval{\max}{\Gr{G}} = \Eigval{1}{\Gr{G}} \geq \Eigval{2}{\Gr{G}}
\geq \ldots \geq \Eigval{n}{\Gr{G}} = \Eigval{\min}{\Gr{G}}.
\end{IEEEeqnarray}
The {\em spectrum} of $\Gr{G}$ consists of the eigenvalues
of $\A$, including their multiplicities. The terms
$\Eigval{1}{\Gr{G}}, \Eigval{2}{\Gr{G}}$, and $\Eigval{n}{\Gr{G}}$
are referred to as the {\em largest, second-largest, and
smallest eigenvalues} of the graph $\Gr{G}$, respectively.
The second-largest and smallest eigenvalues play a key role in {\em spectral
graph theory}, and the interested reader is referred to, e.g.,
\cite{BrouwerH11,CioabaM22,CvetkovicRS09,Spielman12,Stanic15}.

The number of edges that are incident to a vertex in a graph is called the
{\em degree} of the vertex, and a graph is said to be {\em regular} if
its vertices have an identical degree. A regular graph whose all vertices
have a fixed degree $d$ is called a {\em $d$-regular} graph,
and it is said to have {\em valency} $d$.
A $d$-regular graph has a largest eigenvalue $\Eigval{1}{\Gr{G}} = d$,
with the all-ones column vector (of length $n$) as an eigenvector.

Let $\Gr{G}$ be a finite, simple and undirected graph. The {\em complement}
of $\Gr{G}$, denoted by $\CGr{G}$, is defined to have the same vertex
set as $\Gr{G}$, and (by definition) any pair of distinct vertices in
$\V{\Gr{G}}$ are adjacent in $\CGr{G}$ if and only if they are non-adjacent
in the graph $\Gr{G}$.
Hence, we have the equalities
$\card{\E{\Gr{G}}} + \card{\E{\CGr{G}}} = \tfrac12 \, n(n-1)$,
and $\V{\Gr{G}} = \V{\CGr{G}}$. Furthermore, let $\J{n}$ and $\I{n}$
denote the $n$-times-$n$ all-ones and identity matrices, respectively. Then,
by definition, it follows that the adjacency matrices of the graph $\Gr{G}$
and its complement $\CGr{G}$ are related by the equality
\begin{align}
\label{eq:adjacency matrices}
\AG{\CGr{G}} = \J{n} - \I{n} - \AG{\Gr{G}}.
\end{align}

Let $\Gr{G}$ be a $d$-regular graph on $n$ vertices, and let $\CGr{G}$ be the
complement graph. Then, by \eqref{eq:adjacency matrices}, the spectra of $\Gr{G}$
and $\CGr{G}$ are related as follows
\cite[Section~1.3.2]{BrouwerH11}:
\begin{IEEEeqnarray}{rCl}
\label{eq: 21.11.2022a1}
&& \Eigval{1}{\CGr{G}} = n-d-1 = n-1-\Eigval{1}{\Gr{G}}, \\[0.1cm]
\label{eq: 21.11.2022a2}
&& \Eigval{\ell}{\CGr{G}} = -1 - \Eigval{n+2-\ell}{\Gr{G}}, \quad \ell = 2, \ldots, n.
\end{IEEEeqnarray}
Specifically, setting $\ell=n$ in \eqref{eq: 21.11.2022a2} gives
\begin{IEEEeqnarray}{rCl}
\label{eq: 21.11.2022a3}
\Eigval{\min}{\CGr{G}} = -1 - \Eigval{2}{\Gr{G}}.
\end{IEEEeqnarray}

A graph is called {\em acyclic} if it has no cycles, and a connected acyclic
graph is called a {\em tree}. A tree on $n$ vertices has $n-1$ edges, and
for every pair of distinct vertices in a tree, there is a unique path joining
them (see, e.g., \cite[Theorem~5.1.2]{CioabaM22}). A {\em leaf} in a tree is
a vertex of degree~1, and every tree contains at least two leaves (a tree
is therefore not a regular graph, unless it is the complete graph on two vertices
$\CoG{2}$). A disjoint union of trees is a {\em forest} (i.e., a
forest is an acyclic graph whose connected components are trees). In particular,
a deletion of an edge from a tree gives a forest of two trees. The interested
reader is referred to, e.g., \cite[Chapter~5]{CioabaM22} for further properties and
analysis related to trees and forests.

Two simple graphs $\Gr{G} = (\V{\Gr{G}}, \E{\Gr{G}})$ and
$\Gr{H} = (\V{\Gr{H}}, \E{\Gr{H}})$ are said to be {\em isomorphic}
if there is a bijection (i.e., a one-to-one and onto mapping)
$f \colon \V{\Gr{G}} \to \V{\Gr{H}}$ that preserves adjacency
and non-adjacency, i.e., $\{u,v\} \in \E{\Gr{G}}$ if and only if
$\{ f(u), f(v) \} \in \E{\Gr{H}}$. The notation $\Gr{G} \cong \Gr{H}$
denotes that $\Gr{G}$ and $\Gr{H}$ are isomorphic graphs. An
isomorphism from a graph to itself is called an {\em automorphism}
of the graph.

A graph is called {\em self-complementary} if $\Gr{G}$ and $\CGr{G}$
are isomorphic graphs. These include, e.g., the trivial complete graph on one vertex
$\CoG{1}$, the length-3 path $\PathG{4}$, and the 5-cycle graph $\CG{5}$.
If $\Gr{G}$ is a self-complementary
graph of order $n$, then the size of $\Gr{G}$ is
$m = \card{\E{\Gr{G}}} = \tfrac12 \binom{n}{2} = \tfrac14 \, n(n-1)$.
Since only $n$ or $n-1$ can be even, either
$n \equiv 0 \, ( \hspace*{-0.25cm} \mod 4)$ or $n \equiv 1 \, ( \hspace*{-0.25cm} \mod 4)$.
For every such $n$, there exists a recursive algorithm for constructing a
self-complementary graph of order $n$ (see \cite[Exercise~29]{ChartrandLZ15}).
More explicitly, if $\Gr{G}$ is a self-complementary graph of order $n$, then its
disjoint union with the length-3 path $\PathG{4} = [v_1, v_2, v_3, v_4]$, where each of
the vertices $v_2$ and $v_3$ in $\PathG{4}$ is connected to all the vertices in
$\Gr{G}$, gives a self-complementary graph of order $n+4$. Starting with a graph
$\Gr{G}$ that is equal to $\PathG{4}$ or $\CG{5}$ (for graph orders of $n=4$ or $n=5$,
respectively) gives, by the above suggested recursive construction, a self-complementary
graph of order $n$ for all integers $n>1$ such that
$n \equiv 0 \, ( \hspace*{-0.25cm} \mod 4)$ or $n \equiv 1 \, ( \hspace*{-0.25cm} \mod 4)$,
respectively.

A graph $\Gr{G}$ is {\em vertex-transitive} if for every two vertices of
$\Gr{G}$, there is an automorphism of $\Gr{G}$ that maps one vertex to the other.
Similarly, $\Gr{G}$ is said to be an {\em edge-transitive} graph if for every two
edges $\{u_1,v_1\}$ and $\{u_2, v_2\}$ of $\Gr{G}$, there is an automorphism
$f \colon \V{\Gr{G}} \to \V{\Gr{G}}$ of the graph $\Gr{G}$, such that
$\{u_2, v_2\}= \{f(u_1), f(v_1) \}$. A vertex-transitive graph is
necessarily regular, but the opposite does not hold in general. Unlike a
vertex-transitive graph, an edge-transitive graph is not necessarily regular.
Types of graph transitivity are studied, e.g., in \cite[Chapters~3--4]{GodsilR}.

All the eigenvalues of a $d$-regular graph $\Gr{G}$ are bounded in absolute
value by $d$, and the largest eigenvalue is equal to $d$ with an
eigenvector that is equal to the all-ones column vector (see, e.g.,
\cite[Proposition~12.1.1]{CioabaM22}). The multiplicity
of the largest eigenvalue of a $d$-regular graph is~1 if and only if $\Gr{G}$
is connected (see, e.g., \cite[Theorem~4.5.2 and Proposition~12.1.1]{CioabaM22}).
A graph $\Gr{G}$ has its eigenvalues symmetric around zero (including their
multiplicities) if and only if $\Gr{G}$ is bipartite (see, e.g.,
\cite[Theorem~4.3.2]{CioabaM22}). Hence, the smallest eigenvalue of a
$d$-regular bipartite graph is equal to $-d$. Moreover, if $\Gr{G}$ is $d$-regular
and connected, then $\Gr{G}$ is bipartite if and only if $-d$ is an eigenvalue
of its adjacency matrix (see, e.g., \cite[Proposition~12.1.1]{CioabaM22}).
For a $d$-regular graph $\Gr{G}$, let
\begin{IEEEeqnarray}{rCl}
\label{eq:29.11.22}
\lambda(\Gr{G}) \triangleq \max_{\ell: \, \Eigval{\ell}{\Gr{G}} \neq \pm d} \; \bigl| \Eigval{\ell}{\Gr{G}} \bigr|.
\end{IEEEeqnarray}
A connected $d$-regular graph $\Gr{G}$ is called {\em Ramanujan} if
\begin{IEEEeqnarray}{rCl}
\label{eq:Ramanujan}
\lambda(\Gr{G}) \leq 2 \sqrt{d-1}.
\end{IEEEeqnarray}
The reason for the expression in the right-hand side of \eqref{eq:Ramanujan}
is related to the Alon--Boppana bound, which addresses the question
of how small can the second-largest eigenvalue be for a connected
$d$-regular graph or for a sequence of connected $d$-regular graphs
whose number of vertices tends to infinity (the value of the parameter $d$ is kept fixed).
It states that for every connected $d$-regular graph $\Gr{G}$ on $n$ vertices, with $d \geq 3$,
\begin{IEEEeqnarray}{rCl}
\label{eq: AB LB}
\Eigval{2}{\Gr{G}} \geq 2 \sqrt{d-1} - O\bigl((\log_d n)^{-1}\bigr).
\end{IEEEeqnarray}
A non-asymptotic version of \eqref{eq: AB LB} appears in \cite{Nilli91}
(see also \cite[Theorem~12.2.1]{CioabaM22}).
The Alon--Boppana bound was first mentioned in \cite[p.~95]{Alon86},
and it was analyzed, e.g., in \cite{AbbeARS20,AbbeR21,AbbeR22,Cioaba06,
Cioaba08,Friedman91,Friedman93,Friedman08,Li01,Nilli91,Nilli04}.
Moreover, all eigenvalues of a tree, with maximum degree $d \geq 2$,
are (in absolute value) at most $2 \sqrt{d-1}$ (see \cite[Theorem~1]{Murty20}).
Examples of Ramanujan graphs include:
\begin{enumerate}[(a)]
\item The complete $d$-regular graph $\CoG{d+1}$, with $d \geq 2$, whose
eigenvalues are equal to $d$ with multiplicity~1, and $-1$ with multiplicity $d$;
\item The complete bipartite graph $\CoBG{d}{d}$, with $d \geq 2$, is a $d$-regular graph
whose two nonzero eigenvalues are $\pm d$ (each of multiplicity 1), and its
other $2d-2$ eigenvalues are zeros.
\item The Petersen graph, which is isomorphic to the Kneser graph $\KG{5}{2}$,
is a Ramanujan graph since it is $3$-regular with the distinct eigenvalues
3, $-1$, and $-2$.
\end{enumerate}
The interested reader is referred to \cite{Murty20} for a recent survey
paper on Ramanujan graphs, and to references therein (see, e.g., \cite{MarcusSS15}
that proves the existence of infinite families of bipartite Ramanujan
graphs of every degree greater than~2, followed by the extension of that result to
bi-regular bipartite graphs; the proof in \cite{MarcusSS15} uses
an original technique for controlling the eigenvalues of some random matrices).

Let $\Gr{G}$ be a $d$-regular graph of order $n$. The graph $\Gr{G}$
is said to be a {\em strongly regular} graph if there exist nonnegative
integers $\lambda$ and $\mu$ such that the following two conditions hold:
\begin{itemize}
\item Every pair of adjacent vertices have exactly $\lambda$ common neighbors;
\item Every pair of distinct and non-adjacent vertices have exactly $\mu$ common neighbors.
\end{itemize}
Such a strongly regular graph is denoted by $\SRG(n,d,\lambda,\mu)$.
Some basic properties of strongly regular graphs are next introduced, which also serve in our analysis.
\begin{enumerate}[(a)]
\item The complement of a strongly regular graph is also strongly regular.
More explicitly, the complement of $\SRG(n,d,\lambda,\mu)$ is given by
$\SRG(n,n-d-1,n-2d+\mu-2,n-2d+\lambda)$.

\item The four parameters of a strongly regular graph $\SRG(n,d,\lambda,\mu)$ satisfy the relation
\begin{IEEEeqnarray}{rCl}
\label{eq: relation pars. SRG}
(n-d-1) \, \mu = d \, (d-\lambda-1).
\end{IEEEeqnarray}

\item A strongly regular graph $\SRG(n,d,\lambda,\mu)$ has at most three distinct eigenvalues.
If it is connected, then $\Eigval{1}{\Gr{G}} = d$ (multiplicity~1), and the other two distinct
eigenvalues are
\begin{IEEEeqnarray}{rCl}
\label{eq: eigs s.r.g.}
p_{1,2} = \tfrac12 \, \biggl[ \lambda - \mu \pm \sqrt{ (\lambda-\mu)^2 + 4(d-\mu) } \, \biggr],
\end{IEEEeqnarray}
whose respective multiplicities are given by
\begin{IEEEeqnarray}{rCl}
\label{eq: multiplicities s.r.g}
m_{1,2} = \frac12 \Biggl[ n-1 \mp \frac{2d+(n-1)(\lambda-\mu)}{\sqrt{(\lambda-\mu)^2+4(d-\mu)}} \Biggr].
\end{IEEEeqnarray}
Since multiplicities of eigenvalues must be nonnegative integers, their expressions
in \eqref{eq: multiplicities s.r.g} provide further constraints on the values of
$n,d,\lambda$ and $\mu$ (in addition to equality \eqref{eq: relation pars. SRG}).

\item A connected regular graph with exactly three distinct eigenvalues is strongly regular.

\item A strongly regular graph $\SRG(n,d,\lambda,\mu)$, with $\mu>0$, is a connected graph
whose diameter is equal to~2. This holds since two non-adjacent vertices have $\mu>0$ common
neighbors, so the distance between any pair of non-adjacent
vertices is equal to~2. This can be also explained by spectral graph theory since
the diameter of a connected graph is strictly smaller than the number
of its distinct eigenvalues (see \cite[Theorem~4.4.1]{CioabaM22}).
In light of that, the above claim about the diameter holds for all graphs that are
connected and strongly regular since these graphs only have three distinct eigenvalues.

\item If $\mu=0$, the strongly regular graph is disconnected, and it is a disjoint union of
equal-sized complete graphs (i.e., a disjoint union of cliques of the same size). A disjoint
union of an arbitrary number $\ell \geq 2$ of equal-sized complete graphs, $\CoG{d+1}$, has
the parameters $\SRG((d+1)\ell, d, d-1, 0)$. In that case, $d = p_1$ (see \eqref{eq: eigs s.r.g.}),
so the largest and second-largest eigenvalues coincide (by \eqref{eq: multiplicities s.r.g},
that common eigenvalue has multiplicity $m_1+1 = \ell$ in the graph spectrum).
A strongly regular graph $\Gr{G}$ is called {\em primitive} if both $\Gr{G}$ and its
complement $\CGr{G}$ are connected graphs. Otherwise, $\Gr{G}$ is called {\em imprimitive}.
An imprimitive graph is either a disjoint union of equal-sized complete graphs
or its complement, which is a complete multipartite graph. A strongly regular graph
$\Gr{G}$ is imprimitive if and only if 0 or~$-1$ is an eigenvalue of $\Gr{G}$.

\item
Let $\Gr{G}$ be a primitive strongly regular graph $\SRG(n,d,\lambda,\mu)$ with the
largest eigenvalue $d$ (multiplicity~1), second-largest eigenvalue
$r = p_1$ (multiplicity~$m_1$), and smallest eigenvalue
$s = p_2$ (multiplicity~$m_2$). By \eqref{eq: 21.11.2022a1} and \eqref{eq: 21.11.2022a2},
the complement $\CGr{G}$ is a primitive strongly regular graph, having the largest eigenvalue
$n-d-1$ (multiplicity~1), second-largest eigenvalue $-1-s$ (multiplicity~$m_2$), and
smallest eigenvalue $-1-r$ (multiplicity~$m_1$). Each of these primitive strongly regular
graphs has three distinct eigenvalues.
\end{enumerate}
The reader is referred to \cite{BrouwerM22}, which is focused on properties and constructions of
strongly regular graphs.

Part of this work is focused on the strong product of graphs, which is defined as follows.
\begin{definition}
\label{def: strong product}
Let $\Gr{G}_1$ and $\Gr{G}_2$ be two graphs. The {\em strong product}
$\Gr{G} = \Gr{G}_1 \boxtimes \Gr{G}_2$ is a graph whose vertex set is
$\V{\Gr{G}} = \V{\Gr{G}_1} \times \V{\Gr{G}_2}$ (a Cartesian
product), and distinct vertices $\{u_1, u_2\}$ and $\{v_1, v_2\}$ in
$\Gr{G}$ are adjacent if one of the following three conditions is satisfied:
\begin{enumerate}[(a)]
\item $u_1 = v_1$ and $\{u_2, v_2\} \in \E{\Gr{G}_2}$,
\item $\{u_1, v_1\} \in \E{\Gr{G}_1}$ and $u_2 = v_2$,
\item $\{u_1, v_1\} \in \E{\Gr{G}_1}$ and $\{u_2, v_2\} \in \E{\Gr{G}_2}$.
\end{enumerate}
\end{definition}
A strong product of graphs is commutative in the sense that
\begin{eqnarray}
\label{eq:commutative}
\Gr{G}_1 \boxtimes \Gr{G}_2 \cong \Gr{G}_2 \boxtimes \Gr{G}_1,
\end{eqnarray}
and it is also associative in the sense that
\begin{eqnarray}
\label{eq:associative}
(\Gr{G}_1 \boxtimes \Gr{G}_2) \boxtimes \Gr{G}_3 \cong \Gr{G}_1 \boxtimes (\Gr{G}_2 \boxtimes \Gr{G}_3),
\end{eqnarray}
for every three graphs $\Gr{G}_1, \Gr{G}_2$ and $\Gr{G}_3$ (see \cite[Proposition~4.1]{HammackIK11}).

This paper relies on the {\em Lov\'{a}sz $\theta$-function}
of a graph \cite{Lovasz79_IT}, which is next introduced.

\begin{definition}
\label{def: orthogonal representation}
Let $\Gr{G}$ be a simple graph.
An {\em orthogonal representation} of $\Gr{G}$ in the $d$-dimensional Euclidean space $(\Reals^d)$ assigns to each
vertex $i \in \V{\Gr{G}}$ a vector ${\bf{u}}_i \in \Reals^d$ such that ${\bf{u}}_i^{\mathrm{T}} {\bf{u}}_j = 0$
if $\{i, j\} \notin \E{\Gr{G}}$. In other words, the vertices of a simple graph are assigned vectors in $\Reals^d$
such that the vectors that are assigned to any pair of distinct and non-adjacent vertices of that graph are orthogonal.
An {\em orthonormal representation} of $\Gr{G}$ is an orthogonal representation of that graph such that
all representing vectors have unit length.
\end{definition}
\begin{remark}
\label{remark: orthogonal representation}
In an orthogonal representation of a graph $\Gr{G}$, non-adjacent vertices are mapped to orthogonal
vectors, although adjacent vertices are not necessarily mapped to non-orthogonal vectors. If the latter
condition is satisfied, then the orthogonal representation is said to be {\em faithful}.
\end{remark}

\begin{definition}
\label{def: Lovasz theta function}
Let $\Gr{G}$ be a finite, undirected and simple graph. Its {\em Lov\'{a}sz $\theta$-function} is given~by
\begin{IEEEeqnarray}{rCl}
\label{eq: Lovasz theta function}
\theta(\Gr{G}) \triangleq \min_{\bf{u}, \bf{c}} \, \max_{i \in \V{\Gr{G}}} \,
\frac1{\bigl( {\bf{c}}^{\mathrm{T}} {\bf{u}}_i \bigr)^2} \, ,
\end{IEEEeqnarray}
where the minimum is taken over all orthonormal representations $\{{\bf{u}}_i: i \in \V{\Gr{G}} \}$ of $\Gr{G}$,
and all unit vectors ${\bf{c}}$.
The unit vector $\bf{c}$ is called the {\em handle} of the orthonormal representation.
\end{definition}

By the Cauchy-Schwarz inequality,
$\bigl|{\bf{c}}^{\mathrm{T}} {\bf{u}}_i \bigr| \leq \|{\bf{c}} \| \, \|{\bf{u}}_i \| = 1$, so
$\theta(\Gr{G}) \geq 1$ with equality if and only if $\Gr{G}$ is a complete graph.

The Lov\'{a}sz $\theta$-function of a graph can be written as a semidefinite program, which
satisfies strong duality (\cite{GrotschelLS81,Knuth94,Lovasz86} and \cite[Section~11.2]{Lovasz19}).
This enables to compute the value of $\theta(\Gr{G})$ in polynomial
time \cite{GrotschelLS81}. More precisely, there is an algorithm that computes, for every graph
$\Gr{G}$ and every $\varepsilon > 0$, a real number $t$ such that
$\bigl|\theta(\Gr{G}) - t \bigr| < \varepsilon$, where the running time of the algorithm
is polynomial in $n \triangleq \card{\V{\Gr{G}}}$ and $\log \bigl( \frac1\varepsilon \bigr)$ \cite[Theorem~11.11]{Lovasz19}.

The following properties of the Lov\'{a}sz $\theta$-function are used throughout this paper:
\begin{enumerate}[(a)]
\item The sandwich theorem (\cite{GrotschelLS81}, \cite{Knuth94}, \cite[Lemma~3.2.4]{Lovasz86}, \cite[Theorem~11.1]{Lovasz19})
is stated in the two equivalent forms
\begin{IEEEeqnarray}{rCl}
\label{eq1a: sandwich}
&& \indnum{\Gr{G}} \leq \theta(\Gr{G}) \leq \chrnum{\CGr{G}}, \\[0.1cm]
\label{eq1b: sandwich}
&& \clnum{\Gr{G}} \leq \theta(\CGr{G}) \leq \chrnum{\Gr{G}}.
\end{IEEEeqnarray}
\item \cite[Theorem~7]{Lovasz79_IT}: The Lov\'{a}sz $\theta$-function factorizes
for the strong product of graphs, i.e.,
\begin{IEEEeqnarray}{rCl}
\label{eq: Lovasz79 - Theorem 7}
\theta(\Gr{G}_1 \boxtimes \Gr{G}_2) = \theta(\Gr{G}_1) \, \theta(\Gr{G}_2).
\end{IEEEeqnarray}
\item \cite[Corollary~2]{Lovasz79_IT} and \cite[Theorem 8]{Lovasz79_IT}:
\begin{IEEEeqnarray}{rCl}
\label{eq: Lovasz79 - Corollary 2}
\theta(\Gr{G}) \, \theta(\CGr{G}) \geq \card{\V{\Gr{G}}},
\end{IEEEeqnarray}
with equality in \eqref{eq: Lovasz79 - Corollary 2} if the graph $\Gr{G}$
is vertex-transitive.
\item \cite[Theorem~9]{Lovasz79_IT}: Let $\Gr{G}$ be a $d$-regular graph
of order $n$. Then,
\begin{IEEEeqnarray}{rCl}
\label{eq: Lovasz79 - Theorem 9}
\theta(\Gr{G}) \leq -\frac{n \, \Eigval{n}{\Gr{G}}}{d - \Eigval{n}{\Gr{G}}},
\end{IEEEeqnarray}
where $\Eigval{n}{\Gr{G}} < 0$, unless $\Gr{G}$ is an empty graph.
Equality holds in \eqref{eq: Lovasz79 - Theorem 9} if $\Gr{G}$ is edge-transitive.
\item Two simple observations relating the Lov\'{a}sz $\theta$-functions of a graph and its subgraphs:
\begin{itemize}
\item If $\Gr{F}$ is a spanning subgraph of a graph $\Gr{G}$, then $\theta(\Gr{F}) \geq \theta(\Gr{G})$.
\item If $\Gr{F}$ is an induced subgraph of a graph $\Gr{G}$, then $\theta(\Gr{F}) \leq \theta(\Gr{G})$.
\end{itemize}
\item \cite[Theorem~2]{Acin17}: Although unrelated to the analysis in this paper,
another interesting property of the Lov\'{a}sz $\theta$-function is given by the identity
\begin{eqnarray}
\label{eq2: Acin17 - Theorem 2}
\sup_{\Gr{H}} \frac{\indnum{\Gr{G} \boxtimes \Gr{H}}}{\theta(\Gr{G} \boxtimes \Gr{H})} = 1,
\end{eqnarray}
which holds for every simple, finite, and undirected graph $\Gr{G}$, where the supremum
is taken over all such graphs $\Gr{H}$.
This shows that the leftmost inequality in \eqref{eq1a: sandwich} can be made arbitrarily tight
by looking at the strong product of the given graph $\Gr{G}$ with a suitable graph $\Gr{H}$.
\end{enumerate}

The Shannon capacity of a simple, finite and undirected graph $\Gr{G}$ was
introduced in \cite{Shannon56} to determine the maximum
information rate that enables error-free communication. To that end, a discrete
memoryless communication channel is represented by a {\em confusion graph} $\Gr{G}$
that is constructed as follows. The vertices in the graph are represented by the
symbols of the input alphabet to that channel, and any two distinct vertices in that
graph are adjacent if the corresponding two input symbols are not distinguishable
by the channel (in the sense that there exists an output symbol such that the
transition probabilities from each of these two input symbols to that output symbol
are strictly positive). This means that the exact values of the positive
transition probabilities of the channel, as well as the output alphabet of the
channel, are irrelevant to the construction of the confusion graph $\Gr{G}$.
The rationality in doing so is the interest to pictorially represent (by a graph)
all those pairs of input symbols that are not distinguishable by the channel.
Consider a transmission of $k$-length strings.
The {\em $k$-th confusion graph} of the channel is defined as
\begin{eqnarray}
\label{eq: graph power}
\Gr{G}^{\boxtimes \, k} \triangleq \underset{k-1 \,
\text{strong products}}{\underbrace{\Gr{G} \boxtimes \ldots \boxtimes \Gr{G}}},
\end{eqnarray}
which is the $k$-fold strong power of $\Gr{G}$. This is because the independence number of
$\Gr{G}^{\boxtimes \, k}$ is equal to the maximum number of $k$-length strings at the channel
input that can be transmitted with error-free communication (indeed, a pair
of non-adjacent vertices in $\Gr{G}^{\boxtimes \, k}$ represent $k$-length strings that can
be distinguished by the channel, as a result of having a common position
in these two input strings where the corresponding two symbols at that
position are distinguishable by the channel). Consequently, the maximum
{\em information rate per symbol} that is achievable by using input strings of length $k$ is equal to
$\frac1k \, \log \indnum{\Gr{G}^{\boxtimes \, k}} = \log \sqrt[k]{\indnum{\Gr{G}^{\boxtimes \, k}}}$,
for all $k \in \naturals$ (i.e., it is the logarithm of the maximum number of $k$-length input
strings that are distinguishable by the channel, normalized by the length $k$).
The {\em Shannon capacity of a graph} $\Gr{G}$ is defined to be the (exponent of the)
maximum information rate per symbol that is achievable with error-free communication,
where the transmission takes place over a discrete memoryless channel whose confusion
graph is equal to $\Gr{G}$, and the length of the input strings to the channel is unlimited.
It is denoted by $\Theta(\Gr{G})$ (recall that the Lov\'{a}sz $\theta$-function is denoted by $\theta(\Gr{G})$).
Taking the supremum over $k$, the Shannon capacity of $\Gr{G}$ is given by (see \cite{Shannon56},
and \cite[Chapter~42]{AignerZ18})
\begin{eqnarray}
\label{eq:Shannon capacity}
\Theta(\Gr{G}) = \sup_{k \in \naturals} \sqrt[k]{\indnum{G^{\boxtimes \, k}}}.
\end{eqnarray}
The Shannon capacity can be rarely computed exactly
(see, e.g., \cite{Alon02, Alon19, Haemers79, HuTS18, Jurkiewicz14, KornerO98, Lovasz79_IT, Shannon56}).
Analytical observations that also explain why it is, in general, even difficult to approximate it
are addressed in \cite{AlonL06} and \cite{GuoW90}. Calculable upper bounds on $\Theta(\Gr{G})$
were derived by Shannon \cite{Shannon56}, Lov\'{a}sz \cite{Lovasz79_IT}, Haemers \cite{Haemers79},
and more recently by Hu {\em et al.} \cite{HuTS18}. The Lov\'{a}sz $\theta$-function $\theta(\Gr{G})$
is a calculable upper bound on the graph capacity $\Theta(\Gr{G})$, i.e.,
\begin{eqnarray}
\label{eq1: capacity bounds}
\indnum{\Gr{G}} \leq \Theta(\Gr{G}) \leq \theta(\Gr{G}),
\end{eqnarray}
where the leftmost inequality in \eqref{eq1: capacity bounds} follows from \eqref{eq:Shannon capacity}
(by setting $k=1$), and the rightmost inequality in \eqref{eq1: capacity bounds} is \cite[Theorem 1]{Lovasz79_IT}.
In regard to the rightmost inequality in \eqref{eq1: capacity bounds}, it is also
shown in \cite{Alon19} that the Lov\'{a}sz $\theta$-function of a graph, $\theta(\Gr{G})$,
cannot be upper bounded by any function of its Shannon capacity $\Theta(\Gr{G})$.
As mentioned above, the computational task of the Lov\'{a}sz $\theta$-function, $\theta(\Gr{G})$, is
in general feasible by semidefinite programming. Fundamental graph parameters such as its
Shannon capacity $\Theta(\Gr{G})$, independence number $\indnum{\Gr{G}}$, clique number
$\clnum{\Gr{G}}$, and chromatic number $\chrnum{\Gr{G}}$ are all NP-hard problems.
The polynomial-time computability of the Lov\'{a}sz $\theta$-function
of a graph makes inequalities \eqref{eq1a: sandwich}, \eqref{eq1b: sandwich},
and \eqref{eq1: capacity bounds} very useful in obtaining polynomial-time computable upper bounds
on the independence number, clique number, and the Shannon capacity of a graph,
as well as having such a computable lower bound on the chromatic number.

\section{Theorems, Discussions and Examples}
\label{section: results}
The present section provides the results of this work, followed by
examples and discussions. It is composed of five subsections
that address issues related to the Lov\'{a}sz $\theta$-function,
Shannon capacity of graphs, Ramanujan graphs, eigenvalues,
and chromatic numbers of strong graph products or strong graph powers.

\subsection{Bounds on Lov\'{a}sz $\theta$-function, and an exact result
for strongly regular graphs}
\label{subsection: bounds on theta}

Let $\Gr{G}$ be a $d$-regular graph on $n$ vertices, and let $\CGr{G}$
be the complement graph of $\Gr{G}$ that is an $(n-d-1)$-regular graph of order $n$.
An upper bound on $\theta(\Gr{G})$ and a lower bound on $\theta(\CGr{G})$
were obtained by Lov\'{a}sz, expressed in terms of the smallest eigenvalue
of the adjacency matrix of $\Gr{G}$ (see \cite[Theorem~9]{Lovasz79_IT}
and \cite[Corollary~3]{Lovasz79_IT}).
The novelties in the next result (Proposition~\ref{prop1: bounds on theta})
are as follows:
\begin{enumerate}[(a)]
\item It forms a counterpart of a bound by Lov\'{a}sz \cite[Theorem~9]{Lovasz79_IT},
providing a lower bound on $\theta(\Gr{G})$ and an upper bound on
$\theta(\CGr{G})$ that are both expressed in terms of the second-largest
eigenvalue of the adjacency matrix of $\Gr{G}$.
\item It asserts that these two pairs of upper and lower bounds on
$\theta(\Gr{G})$ and $\theta(\CGr{G})$ are tight for the family of
strongly regular graphs. This gives a simple
closed-form expression of the Lov\'{a}sz $\theta$-function of a
strongly regular graph $\SRG(n,d,\lambda,\mu)$ (and the
complement graph) as a function of its four parameters.
\item Further sufficient conditions for the tightness of these bounds
are provided.
\end{enumerate}

\vspace*{0.1cm}
\begin{proposition}
\label{prop1: bounds on theta}
Let $\Gr{G}$ be a $d$-regular graph of order $n$, which is a non-complete
and non-empty graph. Then, the following
bounds hold for the Lov\'{a}sz $\theta$-function of $\Gr{G}$ and its
complement $\CGr{G}$:
\begin{enumerate}[(a)]
\item
\begin{IEEEeqnarray}{rCl}
\label{eq:21.10.22a1}
\frac{n-d+\Eigval{2}{\Gr{G}}}{1+\Eigval{2}{\Gr{G}}} \leq \theta(\Gr{G})
\leq -\frac{n \Eigval{n}{\Gr{G}}}{d - \Eigval{n}{\Gr{G}}}.
\end{IEEEeqnarray}
\begin{itemize}
\item Equality holds in the leftmost inequality of \eqref{eq:21.10.22a1} if $\CGr{G}$
is both vertex-transitive and edge-transitive, or if $\Gr{G}$ is a strongly regular graph;
\item Equality holds in the rightmost inequality of \eqref{eq:21.10.22a1} if $\Gr{G}$
is edge-transitive, or if $\Gr{G}$ is a strongly regular graph.
\end{itemize}

\item
\begin{IEEEeqnarray}{rCl}
\label{eq:21.10.22a2}
1 - \frac{d}{\Eigval{n}{\Gr{G}}} \leq \theta(\CGr{G})
\leq \frac{n \bigl(1+\Eigval{2}{\Gr{G}}\bigr)}{n-d+\Eigval{2}{\Gr{G}}}.
\end{IEEEeqnarray}
\begin{itemize}
\item Equality holds in the leftmost inequality of \eqref{eq:21.10.22a2}
if $\Gr{G}$ is both vertex-transitive and edge-transitive, or if $\Gr{G}$ is
a strongly regular graph;
\item Equality holds in the rightmost inequality of \eqref{eq:21.10.22a2}
if $\CGr{G}$ is edge-transitive, or if $\Gr{G}$ is a strongly regular graph.
\end{itemize}
\end{enumerate}
\end{proposition}
\begin{IEEEproof}
See Section~\ref{subsubsection: proof of prop1}.
\end{IEEEproof}

\vspace*{0.1cm}
\begin{remark}
In light of the sufficient conditions for each of the four inequalities in
Proposition~\ref{prop1: bounds on theta} to hold with equality, define the
following subfamilies of regular graphs:
\begin{itemize}
\item Let $\mathcal{G}_1$ be the family of graphs $\Gr{G}$ such that $\CGr{G}$ is
both vertex-transitive and edge-transitive;
\item Let $\mathcal{G}_2$ be the family of regular and edge-transitive graphs;
\item Let $\mathcal{G}_3$ be the family of graphs $\Gr{G}$ such that $\CGr{G}$
is regular and edge-transitive;
\item Let $\mathcal{G}_4$ be the family of graphs that are
both vertex-transitive and edge-transitive;
\item Let $\mathcal{G}_5$ be the family of the strongly regular graphs.
\end{itemize}
We next show by explicit examples, obtained by using the SageMath software
\cite{SageMath}, that none of the families $\mathcal{G}_1$, $\mathcal{G}_2$
$\mathcal{G}_3$ and $\mathcal{G}_4$ is included in the family $\mathcal{G}_5$,
and vice versa.
\begin{enumerate}[(a)]
\item The Cameron graph is a strongly regular graph
$\SRG(231, 30, 9, 3)$ (see \cite[Section~10.54]{BrouwerM22}). Its complement
is vertex-transitive (hence, regular), but not edge-transitive. This shows that
$\mathcal{G}_5 \not\subseteq \mathcal{G}_3$, so also
$\mathcal{G}_5 \not\subseteq \mathcal{G}_1$.
\item The complement of the Cameron graph is a strongly regular graph
$\SRG(231, 200, 172, 180)$; it is vertex-transitive (hence, regular),
but not edge-transitive. This shows that $\mathcal{G}_5 \not\subseteq \mathcal{G}_2$,
so also $\mathcal{G}_5 \not\subseteq \mathcal{G}_4$.
\item The Foster graph is 3-regular on 90~vertices (see \cite[p.~305]{BrouwerM22}),
which is vertex-transitive and edge-transitive, but it is not strongly regular. This shows that
$\mathcal{G}_4 \not\subseteq \mathcal{G}_5$, so also $\mathcal{G}_2 \not\subseteq \mathcal{G}_5$.
\item The complement of the Foster graph is an 86-regular graph on 90~vertices,
whose complement (i.e., the Foster graph) is vertex-transitive and edge-transitive,
but it is not strongly regular. This shows that
$\mathcal{G}_1 \not\subseteq \mathcal{G}_5$, so also
$\mathcal{G}_3 \not\subseteq \mathcal{G}_5$.
\end{enumerate}
\end{remark}

The next result provides a closed-form expression of the Lov\'{a}sz $\theta$-function
for strongly regular graphs (and their strongly regular complements).
This result relies on Proposition~\ref{prop1: bounds on theta} and the closed-form
expressions of the distinct eigenvalues of a strongly regular graph.

\vspace*{0.1cm}
\begin{corollary}
\label{cor4: Lovasz number for srg}
Let $\Gr{G}$ be a strongly regular graph with parameters $\SRG(n, d, \lambda, \mu)$.
Then,
\begin{IEEEeqnarray}{rCl}
\label{eq:30.10.22a1}
\theta(\Gr{G}) = \frac{n \, (t+\mu-\lambda)}{2d+t+\mu-\lambda},
\end{IEEEeqnarray}
and
\begin{IEEEeqnarray}{rCl}
\label{eq: Lovasz equality for srg}
\theta(\CGr{G}) &=& \frac{n}{\theta(\Gr{G})}  \\
\label{eq:30.10.22a2}
&=& 1 + \frac{2d}{t+\mu-\lambda},
\end{IEEEeqnarray}
where
\begin{IEEEeqnarray}{rCl}
\label{eq: t - srg}
t \triangleq \sqrt{(\mu-\lambda)^2 + 4(d-\mu)}.
\end{IEEEeqnarray}
Furthermore, if $2d + (n-1) \, (\lambda-\mu) \neq 0$, then $\theta(\Gr{G})$
and $\theta(\CGr{G})$ are rational numbers.
\end{corollary}
\begin{IEEEproof}
See Section~\ref{subsubsection: proof of cor4}.
\end{IEEEproof}

\vspace*{0.1cm}
\begin{remark}
By \eqref{eq: Lovasz equality for srg}, if $\Gr{G}$ is
a strongly regular graph on $n$ vertices, then
$\theta(\Gr{G}) \, \theta(\CGr{G}) = n$. This relation is also
known to hold if the graph $\Gr{G}$ is vertex-transitive
\cite[Theorem~8]{Lovasz79_IT}. It should be noted that
not all the strongly regular graphs are necessarily
vertex-transitive, so the observation here is not implied by
\cite[Theorem~8]{Lovasz79_IT}.
As a counter example for strongly regular graphs that are not
vertex-transitive, consider the Chang graphs.
These are three of the existing four non-isomorphic
strongly regular graphs with
parameters $\SRG(28, 12, 6, 4)$ \cite[Section~10.11]{BrouwerM22}
(the fourth such graph, denoted by $\Gr{T}_8$, is the line graph
of the complete graph on 8~vertices $\CoG{8}$). The three Chang
graphs are not vertex-transitive and also not edge-transitive
(in contrast to $\Gr{T}_8$ that is vertex-transitive and edge-transitive),
as it can be verified by the SageMath software \cite{SageMath}.
\end{remark}

\vspace*{0.1cm}
\begin{remark}
The 5-cycle $\CG{5}$ is a strongly regular graph
$\SRG(5,2,0,1)$. Its Lov\'{a}sz $\theta$-function
coincides with its Shannon capacity, being equal to $\sqrt{5}$
(see \cite[Theorem~2]{Lovasz79_IT}). Although it is an irrational number,
it is consistent with Corollary~\ref{cor4: Lovasz number for srg}
since $2d + (n-1) \, (\lambda-\mu) = 2 \cdot 2 + 4 (0-1) = 0$.
\end{remark}

\subsection{Eigenvalue inequalities, strongly regular graphs, and Ramanujan graphs}
\label{subsection: EIGs and Ramanujan graphs}

The present subsection relies on Proposition~\ref{prop1: bounds on theta}, with
the following contributions:
\begin{enumerate}[(a)]
\item Derivation of inequalities that relate the second-largest and smallest
eigenvalues of a regular graph. These inequalities hold with equality if
and only if the graph is strongly regular.
\item Derivation of bounds on parameters of Ramanujan graphs.
\item A more general result is presented for a sequence of regular graphs whose
degrees scale sub-linearly with the orders of these graphs, and their orders
tend to infinity.
\end{enumerate}

This subsection is composed of
Corollaries~\ref{cor0: ineq. for graph eigvals}--\ref{cor3: clique number - Ramanujan},
which all rely on Proposition~\ref{prop1: bounds on theta}. It starts by providing
eigenvalue inequalities.

\vspace*{0.1cm}
\begin{corollary}
\label{cor0: ineq. for graph eigvals}
Let $\Gr{G}$ be a $d$-regular graph of order~$n$, which is non-complete and non-empty. Then,
\begin{IEEEeqnarray}{rCl}
\label{eq1: cor0}
\Eigval{n}{\Gr{G}} \leq -\frac{d \, (n-d+\Eigval{2}{\Gr{G}})}{d + (n-1) \, \Eigval{2}{\Gr{G}}},
\end{IEEEeqnarray}
or equivalently,
\begin{IEEEeqnarray}{rCl}
\label{eq2: cor0}
\Eigval{2}{\Gr{G}} \geq -\frac{d \, (n-d+\Eigval{n}{\Gr{G}})}{d + (n-1) \, \Eigval{n}{\Gr{G}}}.
\end{IEEEeqnarray}
Furthermore, \eqref{eq1: cor0} and \eqref{eq2: cor0} hold with equality if
and only if $\Gr{G}$ is a strongly regular graph.
\end{corollary}
\begin{IEEEproof}
See Section~\ref{subsubsection: proof of cor0}.
\end{IEEEproof}

The next result introduces, in part, Nordhaus–Gaddum type inequalities for the second-largest and smallest
eigenvalues of regular graphs, which are tight for all strongly regular graphs. Regarding
Nordhaus–Gaddum type inequalities, the interested reader is referred to \cite{Nikiforov07, NikiforovY14},
and \cite{NordhausG56}.

\vspace*{0.1cm}
\begin{corollary}
\label{cor1: ineq. for graph eigvals}
Let $\Gr{G}$ be a $d$-regular graph of order~$n$, which is non-complete and non-empty, and let
\begin{IEEEeqnarray}{rCl}
\label{eq: 17.11.22a1}
g_\ell(\Gr{G}) \triangleq \Eigval{\ell}{\CGr{G}}
- \frac{d (n-d+\Eigval{\ell}{\Gr{G}})}{d+(n-1) \Eigval{\ell}{\Gr{G}}},
\qquad \forall \, \ell \in \OneTo{n}.
\end{IEEEeqnarray}
The following holds:
\begin{enumerate}[(a)]
\item
\begin{IEEEeqnarray}{rCl}
\label{eq: 17.11.22a2}
g_n(\Gr{G}) \leq -1 \leq g_2(\Gr{G}),
\end{IEEEeqnarray}
and the two inequalities in \eqref{eq: 17.11.22a2} hold with equality if
and only if $\Gr{G}$ is strongly regular.
\item
If $\Gr{G}$ is a strongly regular graph, then the number of distinct values
in the sequence $\{g_\ell(\Gr{G})\}_{\ell=1}^n$ is either~2 or 3, and
\begin{itemize}
\item it is equal to~2 if the multiplicities of the second-largest and smallest
eigenvalues of $\Gr{G}$ are identical in the subsequence $(\Eigval{2}{\Gr{G}}, \ldots, \Eigval{n}{\Gr{G}})$;
\item it is otherwise equal to~3.
\end{itemize}
\item
If $\Gr{G}$ is self-complementary, then
\begin{IEEEeqnarray}{rCl}
\label{eq: 17.11.22a3}
&& \Eigval{2}{\Gr{G}} \geq \tfrac12 \bigl( \sqrt{n}-1 \bigr), \\[0.1cm]
\label{eq: 17.11.22a4}
&& \Eigval{n}{\Gr{G}} \leq -\tfrac12 \bigl( \sqrt{n}+1 \bigr).
\end{IEEEeqnarray}
\item
If $\Gr{G}$ is self-complementary and strongly regular, then
\eqref{eq: 17.11.22a3} and \eqref{eq: 17.11.22a4} hold with equality.
\end{enumerate}
\end{corollary}
\begin{IEEEproof}
See Section~\ref{subsubsection: proof of cor1 - eigenvals ineq.}.
\end{IEEEproof}

\vspace*{0.1cm}
\begin{example}
Item~(b) of Corollary~\ref{cor1: ineq. for graph eigvals} refers to the
dichotomy in the number of distinct values in the sequence $\{g_\ell(\Gr{G})\}_{\ell=1}^n$.
This statement applies to all strongly regular graphs (either
connected or disconnected). We believe that the following example contributes
to its clarity, in addition to its formal proof in
Section~\ref{subsubsection: proof of cor1 - eigenvals ineq.}.
Let $\Gr{G}$ be a disjoint union of the three complete graphs $\CoG{2}$,
which gives a disconnected strongly regular graph.
Its complement is the complete 3-partite graph $\CGr{G} = \CoG{2,2,2}$, so $n=6$, and
\begin{eqnarray}
\label{eq:17.12.22a1}
&& \{\Eigval{\ell}{\Gr{G}}\}_{\ell=1}^6 = (1, \; 1, \; 1, -1, -1, -1), \\
\label{eq:17.12.22a2}
&& \{\Eigval{\ell}{\CGr{G}}\}_{\ell=1}^6 = (4, \; 0, \; 0, \; 0, -2, -2),  \\
\label{eq:17.12.22a3}
&& \{g_\ell(\Gr{G})\}_{\ell=1}^6 = (3, -1, -1, \; 1, -1, -1).
\end{eqnarray}
By \eqref{eq:17.12.22a1}, the multiplicities of
$\Eigval{2}{\Gr{G}}$ and $\Eigval{n}{\Gr{G}}$ in the spectrum of $\Gr{G}$ are identical, but these
multiplicities are distinct in the subsequence $\{\Eigval{\ell}{\Gr{G}}\}_{\ell=2}^6 = (1, \; 1, -1, -1, -1)$.
Hence, the fact that the sequence $\{g_\ell(\Gr{G})\}_{\ell=1}^6$ gets three distinct values is indeed
consistent with the claim in Item~(b) of Corollary~\ref{cor1: ineq. for graph eigvals}.
\end{example}

\vspace*{0.1cm}
\begin{example}
Let $\Gr{G}$ be the Hall-Janko graph, which is a strongly regular graph with parameters
$\SRG(100, 36, 14, 12)$ (see Section~10.32 of \cite{BrouwerM22}). As a numerical verification
of Item~(b) in Corollary~\ref{cor1: ineq. for graph eigvals} (and its proof
in Section~\ref{subsubsection: proof of cor1 - eigenvals ineq.}),
the sequence $\{g_{\ell}(\Gr{G})\}_{\ell=1}^{100}$ in \eqref{eq: 17.11.22a1}
gets the three distinct values: $n-d-2=62$ (at $\ell=1$), $-1$ (for $2 \leq \ell \leq 37$
or $65 \leq \ell \leq 100$), and~9 (for $38 \leq \ell \leq 64$). The third value (9)
is attained by the sequence $\{g_{\ell}(\Gr{G})\}$ twenty-seven times. By
the proof of Item~(b) in Corollary~\ref{cor1: ineq. for graph eigvals}, the
multiplicity of the third value (27) is equal to the absolute value of the difference between the
multiplicities of the second-largest and smallest eigenvalues of the graph
$\Gr{G}$. The spectrum of the graph $\Gr{G}$ is given by
$36^{1} 6^{36} (-4)^{63}$ (this can be verified by \eqref{eq: eigs s.r.g.}
and \eqref{eq: multiplicities s.r.g}), and the above difference (in absolute value)
is indeed equal to $|63-36| = 27$.
Next, let $\Gr{G}$ be the 5-cycle graph $\CG{5}$, which is a strongly regular graph with
parameters $\SRG(5,2,0,1)$.
The second-largest and smallest eigenvalues of $\Gr{G}$
are equal to $\tfrac12 (\sqrt{5}-1)$ and $-\tfrac12 (\sqrt{5}+1)$, respectively,
and their multiplicities coincide, being both equal to~2. In light of Item~(b) in
Corollary~\ref{cor1: ineq. for graph eigvals},
the sequence $\{g_{\ell}(\Gr{G})\}_{\ell=1}^{5}$ in \eqref{eq: 17.11.22a1}
gets only two distinct values: $n-d-2=1$ at $\ell=1$, and $-1$ for $2 \leq \ell \leq 5$.
\end{example}

\vspace*{0.1cm}
\begin{remark}
\label{remark: cond. -> srg}
We discuss here an implication of the conditions for equalities
in Proposition~\ref{prop1: bounds on theta}.
Let $\Gr{G}$ be a $d$-regular graph.
Inequality \eqref{eq1: cor0} holds with equality if and only if both
inequalities in \eqref{eq:21.10.22a1} hold with equality. By Item~(a)
of Proposition~\ref{prop1: bounds on theta}, the leftmost inequality
in \eqref{eq:21.10.22a1} holds with equality if $\CGr{G}$ is both
vertex-transitive and edge-transitive, and the rightmost inequality
in \eqref{eq:21.10.22a1} holds with equality if $\Gr{G}$ is edge-transitive.
Combining both sufficient conditions for the two inequalities in
\eqref{eq:21.10.22a1} to hold with equality, it follows that a sufficient
condition for equality in \eqref{eq1: cor0} is given by the requirement
that $\Gr{G}$ and $\CGr{G}$ are both vertex-transitive and edge-transitive
(recall that $\CGr{G}$ is vertex-transitive if and only if $\Gr{G}$ is so).

By Corollary~\ref{cor0: ineq. for graph eigvals},
inequality \eqref{eq1: cor0} holds with equality if and only if $\Gr{G}$ is
strongly regular. By a comparison of the former (sufficient) condition with
the latter (necessary and sufficient) condition for equality to hold in
\eqref{eq1: cor0}, it follows that if $\Gr{G}$ and $\CGr{G}$ are both
vertex-transitive and edge-transitive, then $\Gr{G}$ is strongly regular.

It should be noted, with gratitude to an anonymous reviewer, that
a stronger result can be obtained by replacing the requirement of the
vertex-transitivity with the weaker condition of regularity.
Namely, if $\Gr{G}$ is regular, and $\Gr{G}$ and $\CGr{G}$ are both edge-transitive,
then $\Gr{G}$ is strongly regular. This can be shown as follows.
\begin{itemize}
\item
($\Gr{G}$ is edge transitive) $\Rightarrow$
(every edge in $\Gr{G}$ is contained in the same number of
triangles) $\Leftrightarrow$ (every pair of adjacent vertices
in $\Gr{G}$ has the same number of common neighbors);
\item
($\CGr{G}$ is edge transitive) $\Rightarrow$
(for every edge $\{u,v\} \in \E{\CGr{G}}$, the same number of vertices
are not adjacent in $\CGr{G}$ to either $u$ or $v$) $\Leftrightarrow$
(every pair of non-adjacent vertices in $\Gr{G}$ has the
same number of common neighbors);
\item
$\Gr{G}$ is regular (by assumption);
\end{itemize}
and these observations correspond to the conditions in the definition of a strongly
regular graph.
\end{remark}

\vspace*{0.1cm}
\begin{example}
\label{example: ve-srg}
The Schl\"{a}fli graph $\Gr{G}_1$ is 16-regular on~27 vertices. Both $\Gr{G}_1$
and its complement are edge-transitive.
By Remark~\ref{remark: cond. -> srg}, the graph $\Gr{G}_1$ is indeed strongly
regular; its parameters are given by $\SRG(27,16,10,8)$ \cite[Section~10.10]{BrouwerM22}.

Consider, on the other hand, the Shrikhande graph $\Gr{G}_2$ which is 6-regular on
16~vertices. It is a strongly regular graph with parameters $\SRG(16,6,2,2)$
\cite[Section~10.6]{BrouwerM22}. The
graph $\Gr{G}_2$ is edge-transitive, but its complement
is {\em not} edge-transitive (this was verified by the
SageMath software \cite{SageMath}). In addition, the Cameron graph $\Gr{G}_3$ is a
strongly regular graph $\SRG(231, 30, 9, 3)$ \cite[Section~10.54]{BrouwerM22}.
It can be verified that it is edge-transitive, and that its
complement is not edge-transitive.
This shows that the family of regular graphs $\Gr{G}$ with the property that $\Gr{G}$
and its complement $\CGr{G}$ are both edge-transitive
is a {\em strict} subset of the family of strongly regular graphs.
\end{example}

\vspace*{0.1cm}
\begin{remark}
\label{remark: 06.11.22a}
In continuation to Example~\ref{example: ve-srg}, strongly-regular graphs $\Gr{G}$
such that $\Gr{G}$ and $\CGr{G}$ are both edge-transitive include, e.g., the
Hall-Janko, Hoffman-Singleton, Mesner, Petersen, Schl\"{a}fli, Sims-Gewirtz, and
Suzuki graphs (this has been verified by the SageMath software \cite{SageMath};
for the introduction of these graphs, the reader is referred to \cite{BrouwerM22}).
It also includes the infinite families of Paley and Peisert graphs (\cite{Mullin09,Peisert01}),
which are self-complementary and arc-transitive graphs. (All self-complementary and
arc-transitive graphs are strongly regular).
\end{remark}

\vspace*{0.1cm}
\begin{remark}
In connection to Remark~\ref{remark: cond. -> srg}, two related
statements have been proved by Neumaier:
\begin{enumerate}[(a)]
\item A connected, edge-transitive and strongly regular graph is
vertex-transitive \cite[Lemma~1.3]{Neumaier80}.
\item A vertex-transitive and edge-transitive graph containing a regular clique
is strongly regular (see \cite[Corollary~2.4]{Neumaier81}). (A clique $\mathcal{C}$
is called regular if every vertex not in $\mathcal{C}$ is adjacent to the same
positive number of vertices in $\mathcal{C}$).
\end{enumerate}
Strongly regular graphs that are both vertex- and edge-transitive are studied in
\cite{MorrisPS09}.
\end{remark}

\vspace*{0.1cm}
\begin{corollary}
\label{cor1: clique and chromatic nums}
Let $\{\Gr{G}_{\ell}\}_{\ell \in \naturals}$ be a sequence of
graphs where $\Gr{G}_{\ell}$ is $d_\ell$-regular of order
$n_\ell$, and
\begin{IEEEeqnarray}{rCl}
\label{eq: 2 limits}
\underset{\ell \to \infty}{\lim} n_\ell = \infty, \qquad
\underset{\ell \to \infty}{\lim} \, \frac{d_\ell}{n_\ell} = 0.
\end{IEEEeqnarray}
Then,
\begin{IEEEeqnarray}{rCl}
\label{eq:20.10.22b1}
&& \limsup_{\ell \to \infty} \, \clnum{\Gr{G}_{\ell}} \leq a, \\
\label{eq:20.10.22b2}
&& \liminf_{\ell \to \infty} \,  \frac{\chrnum{\CGr{G}_\ell}}{n_\ell} \geq \frac1a,
\end{IEEEeqnarray}
with
\begin{IEEEeqnarray}{rCl}
\label{eq:25.10.22b1}
a \triangleq 1
+ \limsup_{\ell \to \infty} \, \lfloor \Eigval{2}{\Gr{G}_{\ell}} \rfloor.
\end{IEEEeqnarray}
\end{corollary}
\begin{IEEEproof}
See Section~\ref{subsubsection: proof of cor1 - clique and chromatic nums}.
\end{IEEEproof}

\vspace*{0.1cm}
\begin{corollary}
\label{cor2: clique and chromatic nums}
Let $\{\Gr{G}_{\ell}\}_{\ell \in \naturals}$ be a sequence of
Ramanujan $d$-regular graphs where $d \in \naturals$ is fixed,
$\Gr{G}_{\ell}$ is a graph on $n_\ell$ vertices, and
$\underset{\ell \to \infty}{\lim} n_{\ell} = \infty$. Then,
\begin{IEEEeqnarray}{rCl}
\label{eq:21.10.22c1}
&& \limsup_{\ell \to \infty} \, \clnum{\Gr{G}_{\ell}} \leq 1+ \lfloor 2 \sqrt{d-1} \rfloor,  \\[0.1cm]
\label{eq:21.10.22d5}
&& \liminf_{\ell \to \infty} \, \frac{\theta(\Gr{G}_\ell)}{n_\ell} \geq \frac1{1+2 \sqrt{d-1}}, \\[0.1cm]
\label{eq:21.10.22c2}
&& \liminf_{\ell \to \infty} \,  \frac{\chrnum{\CGr{G}_\ell}}{n_\ell} \geq \frac1{1+ \lfloor 2 \sqrt{d-1} \rfloor}.
\end{IEEEeqnarray}
\end{corollary}
\begin{IEEEproof}
See Section~\ref{subsubsection: proof of cor2}.
\end{IEEEproof}

In continuation to Corollary~\ref{cor2: clique and chromatic nums}, the following result provides
non-asymptotic bounds on some graph parameters.

\vspace*{0.1cm}
\begin{corollary}
\label{cor3: clique number - Ramanujan}
Let $\Gr{G}$ be a Ramanujan $d$-regular graph on $n$ vertices. Then,
\begin{IEEEeqnarray}{rCl}
\label{eq:21.10.22b1}
&& \clnum{\Gr{G}} \leq \bigg\lfloor \frac{n \bigl(1+2 \sqrt{d-1}\,\bigr)}{n-d+2 \sqrt{d-1}} \bigg\rfloor \, ,  \\[0.15cm]
\label{eq:21.10.22d3}
&& \theta(\Gr{G}) \geq \frac{n-d+2 \sqrt{d-1}}{1+2 \sqrt{d-1}} \, , \\[0.15cm]
\label{eq:21.10.22b3}
&& \chrnum{\CGr{G}} \geq \bigg\lceil \frac{n-d+2 \sqrt{d-1}}{1+2 \sqrt{d-1}} \, \bigg\rceil.
\end{IEEEeqnarray}
\end{corollary}
\begin{IEEEproof}
See Section~\ref{subsubsection: proof of cor3}.
\end{IEEEproof}

\vspace*{0.1cm}
\begin{remark}
Inequalities~\eqref{eq:20.10.22b1}, \eqref{eq:21.10.22c1}, and \eqref{eq:21.10.22b1} can
be also obtained from \cite[Theorem~2.1.3]{Haemers79_thesis}, which states that
for an arbitrary simple, finite, and undirected graph $\Gr{G}$ on $n$ vertices, whose maximal degree
is given by $\Delta(\Gr{G})$, the following bound on the clique number of $\Gr{G}$ holds:
\begin{IEEEeqnarray}{rCl}
\label{eq: Haemers - Theorem 2.1.3}
\clnum{\Gr{G}} \leq \frac{n \bigl(d+\Eigval{1}{\Gr{G}} \, \Eigval{2}{\Gr{G}} \bigr)}{dn-\Delta(\Gr{G})^2+\Eigval{1}{\Gr{G}} \, \Eigval{2}{\Gr{G}}}.
\end{IEEEeqnarray}
For a $d$-regular graph $\Gr{G}$ on $n$ vertices, we have $\Eigval{1}{\Gr{G}} = d = \Delta(\Gr{G})$,
which then specializes \eqref{eq: Haemers - Theorem 2.1.3} to (see \cite[Theorem~2.1.4]{Haemers79_thesis})
\begin{IEEEeqnarray}{rCl}
\label{eq: Haemers - Theorem 2.1.4}
\clnum{\Gr{G}} \leq \frac{n \bigl(1+\Eigval{2}{\Gr{G}} \bigr)}{n-d+\Eigval{2}{\Gr{G}}}.
\end{IEEEeqnarray}
Inequality~\eqref{eq:21.10.22b1} can be also obtained from \eqref{eq: Haemers - Theorem 2.1.4},
combined with the satisfiability of the inequality
$\bigl| \Eigval{2}{\Gr{G}} \bigr| \leq 2 \sqrt{d-1}$ if $\Gr{G}$ is a
Ramanujan $d$-regular graph, together with the fact that the right-hand side
of \eqref{eq: Haemers - Theorem 2.1.4} is monotonically increasing in the parameter $\Eigval{2}{\Gr{G}}$.
\end{remark}

\subsection{Bounds on eigenvalues of strong products of regular graphs}
\label{subsection: eigenvalues of strong products}

A small second-largest eigenvalue of the adjacency matrix of a regular graph implies
that the graph is a good expander (see, e.g., \cite[Theorem~12.1.2]{CioabaM22},
and \cite{Alon86,Dowling16,HooryLW06,Kahale95,Krivelevich19,Lubotzky12}).

The Alon--Boppana bound in \eqref{eq: AB LB} is a lower bound on the second-largest
eigenvalue of a connected regular graph (see Section~\ref{section: preliminaries}).
By the Alon--Boppana bound, for every sequence $\{\Gr{G}_k\}_{k=1}^{\infty}$ of connected
$d$-regular graphs, with a fixed integer $d \geq 3$ and orders tending to infinity (i.e.,
$\underset{k \to \infty}{\lim} \card{\V{\Gr{G}_k}} = \infty$), the second-largest
eigenvalues of their adjacency matrices satisfy
\begin{IEEEeqnarray}{rCl}
\label{eq: asympt. AB-LB}
\liminf_{k \to \infty} \Eigval{2}{\Gr{G}_k} \geq 2 \sqrt{d-1},
\end{IEEEeqnarray}
(see, e.g., \cite[Theorem~12.1.2]{CioabaM22}).
This lower bound is asymptotically tight, and its tightness can be strengthened beyond
the second-largest eigenvalue. More explicitly, by Serre's theorem \cite{Li01}, for every
fixed integer $d \geq 3$ and (an arbitrarily small) $\varepsilon > 0$, there exists a positive
constant $c = c(\varepsilon, d)$ such that every $d$-regular graph on $n$ vertices
has at least $cn$ eigenvalues that are larger than or equal to $2 \sqrt{d-1} - \varepsilon$.
In other words, Serre's theorem states that a non-vanishing fraction of the $n$ eigenvalues of every $d$-regular graph
has the property of satisfying the Alon--Boppana lower bound within the desired
accuracy (see a simplified proof in \cite[Theorem~12.2.3]{CioabaM22} or
\cite[Theorem~1]{Cioaba06}). Analogous theorems, concerning the least eigenvalues
of $d$-regular graphs, also hold under an additional hypothesis that the graphs do not have odd cycles
below a certain length (see \cite[Section~4]{Cioaba06}). It overall justifies the definition
of Ramanujan graphs (see \eqref{eq:Ramanujan}) as connected $d$-regular graphs
whose all non-trivial eigenvalues are (in absolute value) at most $2 \sqrt{d-1}$.

For a $k$-fold strong power of a $d$-regular graph, the degree is increased exponentially in $k$,
being equal to $d_k = (1+d)^k-1$. It is therefore of interest
to obtain an alternative lower bound on the second-largest eigenvalue of strong products
of regular graphs. Its derivation is motivated by the significance
of strong products, and in particular strong powers of a given graph:
\begin{enumerate}[(1)]
\item The graph capacity in Shannon's problem of zero-error communication \cite{Shannon56}
is given in \eqref{eq:Shannon capacity}, which is expressed in terms of the independence numbers
of all $k$-fold strong powers of the graph (with $k \in \naturals$);
\item The Witsenhausen rate \cite{Witsenhausen76} in the zero-error source coding problem,
with perfect side information at the receiver, is expressed in a dual form to \eqref{eq:Shannon capacity},
where the independence numbers of $k$-fold strong powers of a graph (with $k \in \naturals$)
are replaced by their chromatic numbers, and the supremum over $k$ is replaced by an infimum
(see \cite[Section~3]{Alon02});
\item There exists a polynomial-time algorithm that finds the unique prime factorization of
any connected graph under the operation of strong graph multiplication \cite{FeigenbaumS92}.
\end{enumerate}

It is demonstrated in this section that, under a certain condition, the suggested lower bound
on the second-largest eigenvalue of the $k$-fold strong power of a regular graph offers a
larger exponential growth rate in $k$, as compared to the Alon--Boppana bound.

The following proposition provides a lower bound on the second-largest eigenvalue, and
an upper bound on the smallest eigenvalue of the adjacency matrices of strong products of regular graphs.
Their derivation relies on Proposition~\ref{prop1: bounds on theta},
jointly with the factorization property in \eqref{eq: Lovasz79 - Theorem 7}.
Both bounds are expressed in terms
of the Lov\'{a}sz $\theta$-function of each factor. This enables to obtain analytical
bounds on the second-largest and smallest eigenvalues of a $k$-fold strong power of a
regular graph, with a low computational complexity that is not affected by $k$.
This stays in contrast to the computational complexity of these eigenvalues,
which significantly increases with $k$.

\vspace*{0.1cm}
\begin{proposition}
\label{prop: bounds on eigvals - strong products}
Let $\Gr{G}_1, \ldots, \Gr{G}_k$ be regular graphs such that, for all $\ell \in \OneTo{k}$,
the graph $\Gr{G}_\ell$ is $d_\ell$-regular of order $n_\ell$.
The following bounds hold for their strong product:
\begin{enumerate}[(a)]
\item Unless all $\Gr{G}_\ell$ (with $\ell \in \OneTo{k}$) are complete graphs, then
\begin{IEEEeqnarray}{rCl}
\label{eq: eig2 - LB1}
\Eigval{2}{\Gr{G}_1 \boxtimes \ldots \boxtimes \Gr{G}_k}
& \geq & \frac{\overset{k}{\underset{\ell=1}{\prod}} n_\ell
- \overset{k}{\underset{\ell=1}{\prod}} (1+d_\ell)}{ \overset{k}{\underset{\ell=1}{\prod}}
\theta(\Gr{G}_\ell)-1} - 1 \\[0.2cm]
\label{eq: eig2 - LB2}
& \geq & \frac{\overset{k}{\underset{\ell=1}{\prod}} n_\ell
- \overset{k}{\underset{\ell=1}{\prod}} (1+d_\ell)}{ \overset{k}{\underset{\ell=1}{\prod}}
\biggl(-\dfrac{n_\ell \, \Eigval{\min}{\Gr{G}_\ell}}{d_\ell
- \Eigval{\min}{\Gr{G}_\ell}} \biggr) -1} - 1,
\end{IEEEeqnarray}
and inequality \eqref{eq: eig2 - LB2} holds with equality if, for
all $\ell \in \OneTo{k}$, the regular graph $\Gr{G}_\ell$ is either edge-transitive
or strongly regular.

\item Unless all $\Gr{G}_\ell$ (with $\ell \in \OneTo{k}$) are empty graphs, then
\begin{IEEEeqnarray}{rCl}
\label{eq: eig_min - UB1}
\Eigval{\min}{\Gr{G}_1 \boxtimes \ldots \boxtimes \Gr{G}_k}
& \leq & -\frac{\overset{k}{\underset{\ell=1}{\prod}} (1+d_\ell) - 1}{\overset{k}{\underset{\ell=1}{\prod}}
\biggl( \dfrac{n_\ell}{\theta(\Gr{G}_\ell)} \biggr) - 1} \, .
\end{IEEEeqnarray}
If each regular graph $\Gr{G}_\ell$ (with $\ell \in \OneTo{k}$) is either edge-transitive
or strongly regular, then \eqref{eq: eig_min - UB1} can be expressed in an equivalent form as
\begin{IEEEeqnarray}{rCl}
\label{eq: eig_min - UB2}
\Eigval{\min}{\Gr{G}_1 \boxtimes \ldots \boxtimes \Gr{G}_k}
& \leq & -\frac{\overset{k}{\underset{\ell=1}{\prod}} (1+d_\ell) - 1}{\overset{k}{\underset{\ell=1}{\prod}}
\biggl( 1 - \dfrac{d_\ell}{\Eigval{\min}{\Gr{G}_\ell}} \biggr) - 1}.
\end{IEEEeqnarray}
\end{enumerate}
\end{proposition}
\begin{IEEEproof}
See Section~\ref{subsubsection: proof of prop. 2}.
\end{IEEEproof}

\vspace*{0.1cm}
\begin{remark}
As a sanity check, it would be in place to verify that the lower bound on the second-largest
eigenvalue in the right-hand side of \eqref{eq: eig2 - LB1} is smaller than or equal to the
largest eigenvalue of the strong product:
\begin{IEEEeqnarray}{rCl}
\label{eq: eig_max}
\Eigval{\max}{\Gr{G}_1 \boxtimes \ldots \boxtimes \Gr{G}_k} = \overset{k}{\underset{\ell=1}{\prod}} (1+d_\ell) - 1.
\end{IEEEeqnarray}
First, equality~\eqref{eq: eig_max} holds since $\Gr{G}_1 \boxtimes \ldots \boxtimes \Gr{G}_k$
is $d$-regular with a value of $d$ that is equal to the right-hand side of \eqref{eq: eig_max}.
Straightforward algebra reveals that the required inequality we wish to assert readily follows
from the inequalities
\begin{IEEEeqnarray}{rCl}
\label{eq: two ineq.}
\theta(\Gr{G}_\ell) \geq \alpha(\Gr{G}_\ell) \geq \frac{n_\ell}{1+d_{\ell}}, \quad \ell \in \OneTo{k}.
\end{IEEEeqnarray}
Indeed, the first inequality in \eqref{eq: two ineq.} holds since the Lov\'{a}sz $\theta$-function of a
graph is an upper bound on its graph capacity (see \cite[Theorem~1]{Lovasz79_IT}),
and (by definition) the graph capacity is larger than
or equal to the independence number of the graph. The second inequality in \eqref{eq: two ineq.}
holds by Wei's inequality \cite{Wei81}, which provides a lower bound on the independence number
or the clique number of a finite simple graph as a function of the degrees of its vertices
(see, e.g., \cite[p.~287]{AignerZ18} or \cite[p.~100]{AlonS16} for a nice probabilistic proof
of these inequalities; for some further such bounds, see \cite{Griggs83}).
For a simple graph $\Gr{G}$ of order $n$, where vertex $i \in \OneTo{n}$ is of degree $d_i$,
Wei's bound states that
\begin{IEEEeqnarray}{rCl}
\label{eq: Wei's bound}
\alpha(\Gr{G}) \geq \sum_{i=1}^n \frac1{1+d_i}, \qquad \omega(\Gr{G}) \geq \sum_{i=1}^n \frac1{n-d_i},
\end{IEEEeqnarray}
so the first inequality in \eqref{eq: Wei's bound} is specialized to the second inequality
in \eqref{eq: two ineq.} for a $d_\ell$-regular graph $\Gr{G}_\ell$ of order $n_\ell$.
It should be noted that the pair of inequalities in \eqref{eq: two ineq.} also imply that,
unless not all graphs $\{\Gr{G}_\ell\}$ are empty, the upper bound on
$\Eigval{\min}{\Gr{G}_1 \boxtimes \ldots \boxtimes \Gr{G}_k}$ in the right-hand side
of \eqref{eq: eig_min - UB1} is smaller than or equal to $-1$. It is a desired property of
this upper bound since the smallest eigenvalue of a non-empty and finite regular graph
is smaller than or equal to~$-1$, while attaining this value if the graph is complete.
\end{remark}

\vspace*{0.1cm}
\begin{corollary}
\label{corollary: bounds on eigvals - strong powers}
Let $\Gr{G}$ a $d$-regular graph of order $n$. Then, for all $k \in \naturals$,
the second-largest and smallest eigenvalues of the $k$-fold strong power of
$\Gr{G}$ satisfy the inequalities
\begin{IEEEeqnarray}{rCl}
\label{eq: eig2 LB - graph powers}
\Eigval{2}{\Gr{G}^{\boxtimes \, k}} \geq \frac{n^k - (1+d)^k}{\theta(\Gr{G})^k-1} - 1,
\end{IEEEeqnarray}
and
\begin{IEEEeqnarray}{rCl}
\label{eq: eig_min UB - graph powers}
\Eigval{\min}{\Gr{G}^{\boxtimes \, k}} \leq
-\frac{(1+d)^k-1}{\Bigl(\frac{n}{\theta(\Gr{G})}\Bigr)^k - 1}.
\end{IEEEeqnarray}
\end{corollary}
\begin{IEEEproof}
See Section~\ref{subsubsection: proof of cor. - bounds on eigs., strong product}.
\end{IEEEproof}

\vspace*{0.1cm}
\begin{example}
\label{example: C5}
By Corollary~\ref{corollary: bounds on eigvals - strong powers},
the second-largest eigenvalue of the $k$-fold strong power of the
5-cycle (pentagon) graph $\CG{5}$ satisfies
\begin{IEEEeqnarray}{rCl}
\label{eq: eig2 LB - C5 powers}
\Eigval{2}{\CG{5}^{\boxtimes \, k}} \geq \frac{5^k - 3^k}{5^{\frac{k}{2}}-1} - 1, \qquad k \in \naturals,
\end{IEEEeqnarray}
which holds since $\theta(\CG{5}) = \sqrt{5}$ (it is also the Shannon capacity of the pentagon
\cite[Theorem~2]{Lovasz79_IT}), and $\CG{5}$ is 2-regular ($d=2$).
The lower bound on $\Eigval{2}{\CG{5}^{\boxtimes \, k}}$ in the right-hand side of \eqref{eq: eig2 LB - C5 powers}
scales asymptotically (for large $k$) like $5^{\frac{k}{2}}$. It is next compared with the Alon--Boppana lower bound.
The $k$-fold strong power of $\CG{5}$ is a $d_k$-regular graph with $d_k = (1+d)^k - 1 = 3^k-1$.
The Alon--Boppana lower bound in \eqref{eq: AB LB} is slightly smaller than $2 \sqrt{d_k-1}$, which scales asymptotically like
$2 \cdot 3^{\frac{k}{2}}$. This exemplifies an improvement in the exponential growth rate of the lower bound in
\eqref{eq: eig2 LB - C5 powers}, as compared to the Alon--Boppana lower bound.

We next compare the exact values of $\Eigval{2}{\CG{5}^{\boxtimes \, k}}$, for $1 \leq k \leq 5$,
with their lower bounds in the right-hand side of \eqref{eq: eig2 LB - C5 powers} (these exact values
were calculated numerically by the SageMath software \cite{SageMath}, and their numerical computation
for $k>5$ seems to be a difficult task). The exact values of $\Eigval{2}{\CG{5}^{\boxtimes \, k}}$, for
$k = 1, \ldots, 5$, are equal to $0.6180, \, 3.8541, \, 13.5623, \, 42.6869, \, 130.0608$, respectively
(with 4~digit decimal precision), in comparison to the lower bound in the right-hand side of
\eqref{eq: eig2 LB - C5 powers} that is equal to $0.6180, \, 3.0000, \, 8.6264, \, 21.3333$ and $51.4938$,
respectively.
\end{example}

The following result refers to connected regular graphs that are non-complete and non-empty, and that their
Lov\'{a}sz $\theta$-function is below a certain value (that value depends on the order and valency of the
regular graph). It is asserted that every such a graph has the property that almost all its strong powers
are non-Ramanujan graphs. The derivation of that result relies on inequality~\eqref{eq: eig2 LB - graph powers}.

\vspace*{0.1cm}
\begin{proposition}
\label{prop.: not Ramanujan}
Let $\Gr{G}$ be a connected $d$-regular graph on $n$ vertices, which is non-empty and non-complete. If
\begin{IEEEeqnarray}{rCl}
\label{eq: 09.11.22b1}
\theta(\Gr{G}) < \frac{n}{\sqrt{d+1}},
\end{IEEEeqnarray}
then there exists $k_0 \in \naturals$ such that, for all $k \geq k_0$,
the $k$-fold strong power $\Gr{G}^{\boxtimes \, k}$ is (highly) non-Ramanujan.
An explicit closed-form expression for the value of $k_0$ is given by
\begin{IEEEeqnarray}{rCl}
\label{eq: k_0}
k_0 = \max\left\{ 3, \left\lceil \frac{\log \bigl(2 + (d+1)^{-\frac{3}{2}} \bigr) \,
+ \log \Bigl( \frac{n^3}{n^3-(d+1)^3} \Bigr)}{\log
\Bigl( \frac{n}{\theta(\Gr{G}) \, \sqrt{d+1}} \Bigr)} \right\rceil \right\}.
\end{IEEEeqnarray}
This holds, in particular, for all finite graphs that are self-complementary
and vertex-transitive (they all satisfy the condition
in \eqref{eq: 09.11.22b1} if $n>1$), with a value of $k_0$
which, respectively, is equal to~5, 4 or~3 if $n=5$, $n=9$ or
$n \geq 13$ with $n \equiv 1 \, ( \hspace*{-0.25cm} \mod 4)$.
\end{proposition}
\begin{IEEEproof}
See Section~\ref{subsection: proof of proposition on non-Ramanujan graphs}.
\end{IEEEproof}

\vspace*{0.1cm}
\begin{remark}
\label{remark: n's for sc-ve graphs}
A necessary and sufficient condition for the existence of a self-complementary
graph on $n$ vertices is that $n \equiv 0 \, ( \hspace*{-0.25cm} \mod 4)$ or
$n \equiv 1 \, ( \hspace*{-0.25cm} \mod 4)$ (see \cite[pp.~16--17]{ChartrandLZ15}).
A self-complementary and $d$-regular graph of order $n$ satisfies $d = \frac{n-1}{2}$,
which implies that $n$ needs to be odd. Since a vertex-transitive graph is regular,
the option of $n \equiv 0 \, ( \hspace*{-0.25cm} \mod 4)$ is rejected for
graphs of order $n$ that are self-complementary and vertex-transitive. This implies that the
order $n$ of such graphs must satisfy
$n \equiv 1 \, ( \hspace*{-0.25cm} \mod 4)$. For $n=1$ and $n=5$, there exist
graphs of order $n$ that are self-complementary and vertex-transitive; they are,
respectively, given by $\CoG{1}$ and $\CG{5}$.
Graphs that are self-complementary and vertex-transitive, and approaches for their
construction, received attention in the literature (see, e.g.,
\cite{Li98,LiR14,LiRS21,Lovasz79_IT,Peisert01,Rao14,Rao16}).
\end{remark}

\vspace*{0.1cm}
\begin{remark}
Proposition~\ref{prop.: not Ramanujan} was inspired by Example~\ref{example: C5},
referring to the 5-cycle graph $\CG{5}$ that is self-complementary and vertex-transitive.
It can be numerically verified, with the SageMath software \cite{SageMath}, that
$\CG{5}$ and $\CG{5}^{\boxtimes \, 2}$ are Ramanujan graphs, and then the higher
strong powers $\CG{5}^{\boxtimes \, 3}$, $\CG{5}^{\boxtimes \, 4}$,
$\CG{5}^{\boxtimes \, 5}$ are non-Ramanujan graphs. Proposition~\ref{prop.: not Ramanujan}
shows that all strong powers $\CG{5}^{\boxtimes \, k}$, with $k \geq 5$, are
non-Ramanujan graphs. An expression for $k_0$, as it is given in \eqref{eq: k_0},
was derived in order to reduce the minimal analytical value of $k_0$ for
which the $k$-fold strong power of $\CG{5}$ is asserted to be
non-Ramanujan for all $k \geq k_0$. It started with an initial value of~$k_0=8$
for $n=5$, with a more simple initial expression for $k_0$, and it was reduced
to~$k_0 = 5$ with the closed-form expression in \eqref{eq: k_0}.
Proposition~\ref{prop.: not Ramanujan}, and the numerical experimentation as above,
gives that $\CG{5}^{\boxtimes \, k}$ is a Ramanujan graph if and only if $k=1$ or $k=2$.
\end{remark}

\vspace*{0.1cm}
\begin{example}
Consider the connected Kneser graphs $\Gr{G} = \KG{m}{r}$ with $m > 2r$ and $m,r \in \naturals$. By their
construction and \cite[Theorem~13]{Lovasz79_IT},
\begin{IEEEeqnarray}{rCl}
\label{eq: Kneser par.}
n = \binom{m}{r}, \quad d = \binom{m-r}{r}, \quad \theta(\Gr{G}) = \binom{m-1}{r-1}.
\end{IEEEeqnarray}
The expression of $\theta(\Gr{G})$ in \eqref{eq: Kneser par.} relies on the proof
of \cite[Theorem~13]{Lovasz79_IT},
which shows that $\alpha(\Gr{G}) = \Theta(\Gr{G}) = \theta(\Gr{G})$ (it holds by combining the
Erd\"{o}s-Ko-Rado theorem for finite sets, which serves to determine the independence number of
the graph, together with the upper bound on $\theta(\Gr{G})$ in \cite[Theorem~9]{Lovasz79_IT}).
Straightforward algebra shows that the condition in \eqref{eq: 09.11.22b1} is not
necessarily satisfied for the set of parameters in \eqref{eq: Kneser par.}. For example, by selecting
$m = 2r+1$, the condition in \eqref{eq: 09.11.22b1} is satisfied if and only if
$1 \leq r \leq 3$ (e.g., it is satisfied by the Petersen graph, which corresponds to $r=2$).
The condition in \eqref{eq: 09.11.22b1} is violated for $r \geq 4$ so, for these
values of the parameter $r$, the exponential growth rate of the Alon--Boppana lower bound on
$\Eigval{2}{\Gr{G}^{\boxtimes \, k}}$ is larger than the exponent of the suggested lower bound
in the right-hand side of \eqref{eq: eig2 LB - graph powers}.
\end{example}

\vspace*{0.1cm}
\begin{remark}
Several methods for obtaining upper bounds on the smallest eigenvalue of a graph
$\Gr{G}$ rely on the identity $\Eigval{\min}{\Gr{G}} =
\underset{{\bf{x} \neq 0}}{\min} \, \frac{{\bf{x}}^T \A \bf{x}}{\bf{x}^T \bf{x}}$
(see, e.g., \cite{DesaiR}).
The dimensions of the adjacency matrix of a $k$-fold strong power of a graph
$\Gr{G}$ grow exponentially in $k$, so the computational complexity
of such bounds is typically very high for strong graph powers. However, for strong powers,
the eigenvalue bounds in Corollary~\ref{corollary: bounds on eigvals - strong powers}
are analytical, and their computational complexity is not affected by $k$.
The reader is also referred to recent works which derive lower bounds on the smallest
eigenvalue of a graph based on graph decompositions \cite{CioabaEG20}, and
lower bounds on the smallest eigenvalue of regular graphs containing many copies of
a smaller fixed subgraph \cite{CioabaG22}. Such lower bounds on the smallest eigenvalue
of a regular graph transform to upper bounds on the second-largest eigenvalue of
regular graphs by using equality~\eqref{eq: 21.11.2022a3}. Some upper bounds on the second-largest
eigenvalue of connected graphs, with conditions for their attainability, are provided
in \cite{ZhaiLW12}. The paper \cite{Hong93} surveys (less recent) bounds on the eigenvalues
of simple~graphs.
\end{remark}

\subsection{Lower bounds on the chromatic numbers of strong products}
\label{subsection: LB on the chromatic number of strong products}

The problem of relating the chromatic number of a graph to its eigenvalues dates back to
Haemers \cite[Section~2.2]{Haemers79_thesis}, Hoffman \cite{Hoffman70} and Wilf \cite{Wilf67},
obtaining upper and lower bounds on the chromatic numbers of graphs in terms of
their largest, second-largest and smallest eigenvalues.

This section presents lower bounds on the chromatic number of a strong product of graphs
(and its complement) in terms of the Lov\'{a}sz $\theta$-functions (or smallest
eigenvalues) of its factors.

\vspace*{0.1cm}
\begin{proposition}
\label{prop. chromatic numbers}
The following lower bounds on chromatic numbers of strong products hold:
\begin{enumerate}[(a)]
\item
Let $\Gr{G}_1, \ldots, \Gr{G}_k$ be $k$ simple graphs,
$\bigcard{\V{\Gr{G}_\ell}} = n_\ell$ for $\ell \in \OneTo{k}$,
and $\Gr{G} = \Gr{G}_1 \boxtimes \ldots \boxtimes \Gr{G}_k$. Then,
\begin{IEEEeqnarray}{rCl}
\label{eq:26.10.22a1}
&& \chrnum{\Gr{G}} \geq \Bigg\lceil \prod_{\ell=1}^k \frac{n_\ell}{\theta(\Gr{G}_\ell)} \Bigg\rceil, \\
\label{eq:26.10.22a2}
&& \chrnum{\CGr{G}} \geq \Bigg\lceil \prod_{\ell=1}^k \theta(\Gr{G}_\ell) \Biggr\rceil.
\end{IEEEeqnarray}

\item
Let $\Gr{G}_1, \ldots, \Gr{G}_k$ be regular graphs, where $\Gr{G}_\ell$ is $d_\ell$-regular
of order $n_\ell$ for all $\ell \in \OneTo{k}$. Then,
\begin{IEEEeqnarray}{rCl}
\label{eq:26.10.22a3}
\chrnum{\Gr{G}} &\geq& \Bigg\lceil \prod_{\ell=1}^k \frac{n_\ell}{\theta(\Gr{G}_\ell)} \Bigg\rceil \\[0.1cm]
\label{eq:26.10.22a4}
&\geq& \Bigg\lceil \prod_{\ell=1}^k \biggl( 1 - \frac{d_\ell}{\Eigval{\min}{\Gr{G}_\ell}} \biggr) \bigg\rceil,
\end{IEEEeqnarray}
and inequality \eqref{eq:26.10.22a4} holds with equality if each regular graph $\Gr{G}_\ell$
is either edge-transitive or strongly regular.

\item
If, for all $\ell \in \OneTo{k}$, $\Gr{G}_\ell$ is $d_\ell$-regular, and it is either
edge-transitive or strongly regular, then
\begin{IEEEeqnarray}{rCl}
\label{eq:26.10.22a5}
\prod_{\ell=1}^k \biggl( 1 - \frac{d_\ell}{\Eigval{\min}{\Gr{G}_\ell}} \biggr)
\geq 1 - \frac{d(\Gr{G})}{\Eigval{\min}{\Gr{G}}},
\end{IEEEeqnarray}
where
\begin{IEEEeqnarray}{rCl}
\label{eq:26.10.22a6}
d(\Gr{G}) = \prod_{\ell=1}^k (1+d_\ell) - 1
\end{IEEEeqnarray}
is the valency of the regular graph $\Gr{G} = \Gr{G}_1 \boxtimes \ldots \boxtimes \Gr{G}_k$,
and $\Eigval{\min}{\Gr{G}}$ is its smallest eigenvalue.

\item
Let $\Gr{G}_1, \ldots, \Gr{G}_k$ be regular graphs, where $\Gr{G}_\ell$ is $d_\ell$-regular
of order $n_\ell$ for all $\ell \in \OneTo{k}$.
\begin{enumerate}[(1)]
\item If, for all $\ell \in \OneTo{k}$, the graph $\Gr{G}_\ell$ is either
vertex-transitive or strongly regular, then
the lower bound on $\chrnum{\Gr{G}}$ in the right-hand side
of \eqref{eq:26.10.22a1} is larger than or equal to the lower bound
$\overset{k}{\underset{\ell=1}{\prod}} \clnum{\Gr{G}_\ell}$.  \vspace*{0.1cm}
\item If, for all $\ell \in \OneTo{k}$, the graph $\Gr{G}_\ell$ is
either (i) both vertex-transitive and edge-transitive, or (ii)~strongly
regular, then the lower bound on $\chrnum{\Gr{G}}$ in the
right-hand side of \eqref{eq:26.10.22a4} is larger than or equal to the
lower bound $\overset{k}{\underset{\ell=1}{\prod}} \clnum{\Gr{G}_\ell}$.
\end{enumerate}

\item
Let, for all $\ell \in \OneTo{k}$, the graph $\Gr{G}_\ell$ be $d_\ell$-regular
on $n_\ell$ vertices, and suppose that it is either
edge-transitive or strongly regular. Then,
\begin{IEEEeqnarray}{rCl}
\label{eq:26.10.22a7}
\chrnum{\CGr{G}} \geq \Bigg\lceil \prod_{\ell=1}^k
\biggl( -\frac{n_\ell \, \Eigval{\min}{\Gr{G}_\ell}}{d_\ell
- \Eigval{\min}{\Gr{G}_\ell}} \biggr) \Bigg\rceil.
\end{IEEEeqnarray}

\item
Let, for all $\ell \in \OneTo{k}$, $\Gr{G}_\ell$ be a
self-complementary graph on $n_\ell$ vertices that is
either vertex-transitive or strongly regular.
Let $n \triangleq \overset{k}{\underset{\ell=1}{\prod}} n_\ell$
be the order of $\Gr{G} = \Gr{G}_1 \boxtimes \ldots \boxtimes \Gr{G}_n$.
Then,
\begin{IEEEeqnarray}{rCl}
\label{eq:26.10.22a8.1}
&& \chrnum{\Gr{G}} \geq \big\lceil \sqrt{n} \, \big\rceil, \\[0.1cm]
\label{eq:26.10.22a8.2}
&& \chrnum{\CGr{G}} \geq \big\lceil \sqrt{n} \, \big\rceil.
\end{IEEEeqnarray}
\end{enumerate}
\end{proposition}
\begin{IEEEproof}
See Section~\ref{subsubsection: proof of prop. - chromatic numbers of strong products}.
\end{IEEEproof}

\vspace*{0.1cm}
\begin{remark}
The following inequality is proved in \cite[Theorem~11]{Acin17}:
\begin{IEEEeqnarray}{rCl}
\label{eq: Acin17 - Theorem 11}
\chrnum{\CGr{G}_1 \boxtimes \CGr{G}_2} \geq \theta(\Gr{G}_1) \, \theta(\Gr{G}_2).
\end{IEEEeqnarray}
(Analogous inequalities to \eqref{eq: Acin17 - Theorem 11} were derived by Hales
\cite{Hales73}, and by McEliece and Posner \cite{McElieceP71}; see \cite[Theorem~10]{Acin17}).
Inequality \eqref{eq:26.10.22a1} can be obtained by combining \eqref{eq: Lovasz79 - Theorem 7},
\eqref{eq: Lovasz79 - Corollary 2} and \eqref{eq: Acin17 - Theorem 11}. This can be done by first
replacing $\Gr{G}_1$ and $\Gr{G}_2$ in \eqref{eq: Acin17 - Theorem 11} with
$\Gr{G}'_1 \triangleq \CGr{G_1}$ and $\Gr{G}'_2 \triangleq \CGr{\Gr{G}_2 \boxtimes \ldots \boxtimes \Gr{G}_k}$,
respectively, then relying on \eqref{eq: Lovasz79 - Corollary 2} to get the inequalities
$\theta(\Gr{G}'_1) \geq \frac{n_1}{\theta(\Gr{G}_1)}$ and
$\theta(\Gr{G}'_2) \geq \frac{n_2 \ldots n_k}{\theta(\CGr{G'_2})}$, and finally relying on
\eqref{eq: Lovasz79 - Theorem 7} to get the equality $\theta(\CGr{G'_2}) = \theta(\Gr{G}_2) \ldots \theta(\Gr{G}_k)$.
Our two simple proofs of \eqref{eq:26.10.22a1} are, however, easier than this one, and they do not require to
rely on \eqref{eq: Acin17 - Theorem 11}.
\end{remark}

\vspace*{0.1cm}
\begin{remark}
Inequality~\eqref{eq: Acin17 - Theorem 11} differs from \eqref{eq:26.10.22a2}, although the lower
bounds in the right-hand sides are similar for $k=2$.
The difference between \eqref{eq:26.10.22a2} and \eqref{eq: Acin17 - Theorem 11}
is that the left-hand side of \eqref{eq:26.10.22a2} refers to the chromatic
number of $\CGr{\Gr{G}_1 \boxtimes \Gr{G}_2}$, whereas the left-hand side of
\eqref{eq: Acin17 - Theorem 11} refers to $\CGr{\Gr{G}_1} \boxtimes \CGr{\Gr{G}_2}$.
The result in \eqref{eq:26.10.22a2} appears in Proposition~\ref{prop. chromatic numbers}
for completeness since it is an immediate consequence of the sandwich theorem in
\eqref{eq1a: sandwich}, and the identity in \eqref{eq: Lovasz79 - Theorem 7}.
\end{remark}

\vspace*{0.1cm}
\begin{remark}
Under the assumptions of Item~(c) in Proposition~\ref{prop. chromatic numbers},
the lower bound on $\chrnum{\Gr{G}}$ in the right-hand side of \eqref{eq:26.10.22a4}
is larger than or equal to Hoffman's lower bound (see \eqref{eq:26.10.22a5}).
This holds in addition to the high complexity in computing $\Eigval{\min}{\Gr{G}}$ in the
right-hand side of \eqref{eq:26.10.22a5}, even for relatively small
values of $k$. Hoffman's bound was originally proved for regular graphs
\cite{Hoffman70}, and it was later extended by Haemers to general simple and finite graphs
(see, \cite[Proposition~3.5.3]{BrouwerH11}, \cite[Theorem~2.1.3]{Haemers79_thesis},
\cite[Corollary~8.10 and Theorem~8.11]{Stanic15}), and to hypergraphs~\cite{FilmusGL21}.
A recent perspective on Hoffman's bound appears in \cite{Haemers21}.
\end{remark}

\vspace*{0.1cm}
\begin{remark}
We refer to Item~(f) of Proposition~\ref{prop. chromatic numbers}. If
$\Gr{G}$ is self-complementary of order $n$, then
\begin{IEEEeqnarray}{rCl}
\label{eq:26.10.22a9}
\chrnum{\Gr{G}} \geq \big\lceil \sqrt{n} \, \big\rceil.
\end{IEEEeqnarray}
Indeed, for every graph $\Gr{G}$,
\begin{IEEEeqnarray}{rCl}
\label{eq:26.10.22a10}
\chrnum{\Gr{G}} \, \chrnum{\CGr{G}} && \geq \theta(\Gr{G}) \, \theta(\CGr{G}) \\
\label{eq:26.10.22a11}
&& \geq n,
\end{IEEEeqnarray}
where inequality~\eqref{eq:26.10.22a10} holds since, by \eqref{eq1a: sandwich} and \eqref{eq1b: sandwich},
$\chrnum{\Gr{G}} \geq \theta(\CGr{G})$ and $\chrnum{\CGr{G}} \geq \theta(\Gr{G})$,
and inequality \eqref{eq:26.10.22a11} holds by \cite[Corollary~2]{Lovasz79_IT}. The chromatic
number is invariant under isomorphism, so $\chrnum{\Gr{G}} = \chrnum{\CGr{G}}$
if $\Gr{G}$ is self-complementary. This gives \eqref{eq:26.10.22a9}
from \eqref{eq:26.10.22a11}. It should be noted, however, that the result
in Item~(f) (see \eqref{eq:26.10.22a8.1} and \eqref{eq:26.10.22a8.2}) is
not implied by \eqref{eq:26.10.22a9}. This is because a strong product of
self-complementary graphs, where each factor is also either vertex-transitive
or strongly regular, does not necessarily give a self-complementary graph.
As a counter example, let $\Gr{G}_1 = \Gr{G}_2 = \CG{5}$. The pentagon $\CG{5}$
is a self-complementary, vertex-transitive and strongly regular graph, whereas
$ \Gr{G} = \CG{5} \boxtimes \CG{5} $ is not a self-complementary graph by
\cite[Theorem~1.1]{LiR14}; according to that theorem, $\Gr{G}$ would have been
self-complementary if the strong product of $\Gr{G}_1$ and $\Gr{G}_2$ had
been replaced by their lexicographic product (these last two observations
were also verified by the SageMath software).
\end{remark}

The next result specializes Item~(a) of Proposition~\ref{prop. chromatic numbers}
(see \eqref{eq:26.10.22a1}) to strong products of
strongly regular graphs. This result is obtained by relying on
Corollary~\ref{cor4: Lovasz number for srg} that provides a closed-form expression
for the Lov\'{a}sz $\theta$-function of each factor.

\vspace*{0.1cm}
\begin{corollary}
\label{cor8: chrnum strong product srg}
Let $\Gr{G}_1, \ldots, \Gr{G}_k$ be strongly regular graphs
with parameters $\SRG(n_\ell, d_\ell, \lambda_\ell, \mu_\ell)$
for $\ell \in \OneTo{k}$ (they need not be distinct). Then,
the chromatic number of their strong product satisfies
\begin{IEEEeqnarray}{rCl}
\label{eq:02.11.2022b1}
&& \Bigg\lceil \prod_{\ell = 1}^k \biggl( 1 + \frac{2d_\ell}{t_\ell+\mu_\ell-\lambda_\ell} \biggr) \Bigg\rceil
\leq \chrnum{\Gr{G}_1 \boxtimes \ldots \boxtimes \Gr{G}_k}
\leq \prod_{\ell = 1}^k  \chrnum{\Gr{G}_k},
\end{IEEEeqnarray}
where $\{t_\ell\}_{\ell=1}^k$ in the leftmost term of \eqref{eq:02.11.2022b1} is given by
\begin{IEEEeqnarray}{rCl}
\label{eq:02.11.2022b2}
t_\ell \triangleq \sqrt{(\lambda_\ell - \mu_\ell)^2 + 4(d_\ell - \mu_\ell)}, \quad \ell \in \OneTo{k}.
\end{IEEEeqnarray}
The leftmost term in \eqref{eq:02.11.2022b1} is also larger than or equal to the
product of the clique numbers of the factors $\{\Gr{G}_\ell\}_{\ell=1}^k$.
\end{corollary}
\begin{IEEEproof}
See Section~\ref{subsubsection: proof of cor. chr. strong prod. srg}.
\end{IEEEproof}

The next two examples present strong products of strongly regular graphs, whose chromatic
numbers are exactly determined by Corollary~\ref{cor8: chrnum strong product srg}.

\vspace*{0.1cm}
\begin{example}
\label{example: G_1,2,3}
Let $\Gr{G}_1, \Gr{G}_2, \Gr{G}_3$ be the Schl\"{a}fli, Shrikhande,
and Hall-Janko graphs, respectively. These are strongly regular graphs
whose parameters are $\SRG(27,16,10,8)$, $\SRG(16,6,2,2)$, and
$\SRG(100, 36, 14, 12)$, respectively. Their chromatic numbers are equal to
$\chrnum{\Gr{G}_1} = 9$, $\chrnum{\Gr{G}_2} = 4$, and $\chrnum{\Gr{G}_3} = 10$.
Consider the chromatic number of the strong product of arbitrary
nonnegative powers of $\Gr{G}_1$, $\Gr{G}_2$ and $\Gr{G}_3$.
It can be verified that, for all such strong products, the upper and lower bounds in
Corollary~\ref{cor8: chrnum strong product srg} coincide, so for all
integers $k_1, k_2, k_3 \geq 0$,
\begin{IEEEeqnarray}{rCl}
\label{eq:27.10.22a3}
\chrnum{\Gr{G}_1^{\boxtimes \, k_1} \boxtimes \Gr{G}_2^{\boxtimes \, k_2} \boxtimes \Gr{G}_3^{\boxtimes \, k_3}}
= 9^{k_1} \, 4^{k_2} \, 10^{k_3}.
\end{IEEEeqnarray}
For comparison, the lower bound that is given by the product of the clique numbers
of each factor is equal to $6^{k_1} 3^{k_2} 4^{k_3}$
(since $\clnum{\Gr{G}_1} = 6$, $\clnum{\Gr{G}_2} = 3$, and $\clnum{\Gr{G}_3} = 4$).
This shows that it is significantly looser than the tight lower bound in the
right-hand side of \eqref{eq:27.10.22a3}.
\end{example}

\vspace*{0.1cm}
\begin{example}
\label{example: chang graphs}
Consider the three non-isomorphic Chang graphs. These are strongly regular graphs
with the same set of parameters $\SRG(28,12,6,4)$. The clique number of one of these
graphs is equal to~5, and the clique numbers of the other two graphs are equal to~6.
Let us denote these graphs by $\Gr{G}_1$, $\Gr{G}_2$ and $\Gr{G}_3$, such that
$\clnum{\Gr{G}_1} = 5$, $\clnum{\Gr{G}_2} = 6$, and $\clnum{\Gr{G}_3} = 6$.
The chromatic numbers of all these three graphs are similar, and they are equal to~7,
i.e., $\chrnum{\Gr{G}_1} = \chrnum{\Gr{G}_2} = \chrnum{\Gr{G}_3} = 7$ (the clique and
chromatic numbers of the three Chang graphs are easy to verify with the SageMath
software \cite{SageMath}).
Let $k_1, k_2$ and $k_3$ be arbitrary nonnegative integers.
By Corollary~\ref{cor8: chrnum strong product srg},
\begin{IEEEeqnarray}{rCl}
\label{eq:14.11.2022b1}
\chrnum{\Gr{G}_1^{\boxtimes \, k_1} \boxtimes \Gr{G}_2^{\boxtimes \, k_2} \boxtimes \Gr{G}_3^{\boxtimes \, k_3}}
= 7^{k_1 + k_2 + k_3}
\end{IEEEeqnarray}
due to the coincidence of the upper and lower bounds in \eqref{eq:02.11.2022b1}.
For comparison, the lower bound on the chromatic number in the left-hand side of
\eqref{eq:14.11.2022b1}, which is given by the product of the clique numbers
of each factor, is equal to $5^{k_1} 6^{k_2+k_3}$.
\end{example}

The next example presents a power product of a vertex and edge-transitive regular graph,
whose chromatic number is exactly determined by Item~(b) of Proposition~\ref{prop. chromatic numbers}.

\vspace*{0.1cm}
\begin{example}
\label{example: Perkel}
The present example provides the exact value of $\chrnum{\Gr{G}^{\boxtimes \, k}}$ where $\Gr{G}$
is the Perkel graph, and $k \in \naturals$. The Perkel graph is 6-regular on 57~vertices, and it
is both vertex-transitive and edge-transitive. The clique and chromatic numbers of $\Gr{G}$ are
equal to $\clnum{\Gr{G}}=2$ and $\chrnum{\Gr{G}}=3$, respectively, and the smallest eigenvalue
of (the adjacency matrix of) $\Gr{G}$ is equal to $\Eigval{\min}{\Gr{G}} = -3$. By
Item~(b) of Proposition~\ref{prop. chromatic numbers}, for all $k \in \naturals$,
\begin{IEEEeqnarray}{rCl}
\label{eq: chrnum Perkel}
\chrnum{\Gr{G}^{\boxtimes \, k}} \geq \biggl( 1 - \frac{d}{\Eigval{\min}{\Gr{G}}} \biggr)^k = 3^k,
\end{IEEEeqnarray}
which can be compared here to the simple upper and lower bounds on $\chrnum{\Gr{G}^{\boxtimes \, k}}$
whose values are given by $\chrnum{\Gr{G}}^k=3^k$ and $\clnum{\Gr{G}}^k=2^k$, respectively. The
coincidence of the improved lower bound in \eqref{eq: chrnum Perkel} and the upper bound gives that,
for all $k \in \naturals$,
\begin{IEEEeqnarray}{rCl}
\label{eq2: chrnum Perkel}
\chrnum{\Gr{G}^{\boxtimes \, k}} = 3^k.
\end{IEEEeqnarray}
\end{example}

The next two examples illustrate numerically Part~2 of Item~(d) in Proposition~\ref{prop. chromatic numbers}.

\vspace*{0.1cm}
\begin{example}
\label{example: Suzuki}
Let $\Gr{G}$ be the Suzuki graph, which is a strongly regular graph with parameters
$\SRG(1782, 416, 100, 96)$ (see \cite[Section~10.83]{BrouwerM22}).
The lower bound in Corollary~\ref{cor8: chrnum strong product srg} (see the leftmost
term in \eqref{eq:02.11.2022b1}) gives that $\chrnum{\Gr{G}^{\boxtimes \, k}} \geq 27^k$
for all $k \in \naturals$. For comparison, since $\clnum{\Gr{G}} = 6$, the lower bound
that is based on the clique number of $\Gr{G}$ gives $\chrnum{\Gr{G}^{\boxtimes \, k}} \geq 6^k$
for all $k \in \naturals$. The exact value of $\chrnum{\Gr{G}}$ is not available for the Suzuki
graph, so the upper bound in the rightmost term of \eqref{eq:02.11.2022b1} is unknown.
The improvement in the exponential lower bound on the chromatic number of the strong
powers $\Gr{G}^{\boxtimes \, k}$ is, however, significant since it is increased from $6^k$ to $27^k$.
\end{example}

\vspace*{0.1cm}
\begin{example}
\label{example: Gosset}
Let $\Gr{G}$ be the Gosset graph. We apply here Item~(b) of Proposition~\ref{prop. chromatic numbers}
in order to obtain an improved lower bound on $\chrnum{\Gr{G}^{\boxtimes \, k}}$ for all $k \in \naturals$.
The graph $\Gr{G}$ is 27-regular on 56~vertices (i.e., $d=27$ and $n=56$); it is both
vertex-transitive and edge-transitive, and it is also not strongly regular.
The clique and chromatic numbers of $\Gr{G}$ are equal to $\clnum{\Gr{G}}=7$ and $\chrnum{\Gr{G}}=14$,
respectively, and the smallest eigenvalue of (the adjacency matrix of) $\Gr{G}$ is equal to
$\Eigval{\min}{\Gr{G}} = -3$. By Item~(b) of Proposition~\ref{prop. chromatic numbers}
(note that the edge-transitivity of $\Gr{G}$ implies that \eqref{eq:26.10.22a4} holds with equality),
it follows that for all $k \in \naturals$,
\begin{IEEEeqnarray}{rCl}
\chrnum{\Gr{G}^{\boxtimes \, k}} &\geq& \biggl( 1 - \frac{d}{\Eigval{\min}{\Gr{G}}} \biggr)^k \\
\label{eq: chrnum Gosset}
&=& 10^k,
\end{IEEEeqnarray}
which can be compared here to the simple upper and lower bounds on $\chrnum{\Gr{G}^{\boxtimes \, k}}$
whose values are given by $\chrnum{\Gr{G}}^k=14^k$ and $\clnum{\Gr{G}}^k=7^k$, respectively.
\end{example}

All the factors of the strong products in Examples~\ref{example: G_1,2,3}--\ref{example: Gosset}
are either strongly regular or otherwise, they are both vertex-transitive and
edge-transitive graphs.
Consequently, by Item~(b) of Proposition~\ref{prop. chromatic numbers},
the two lower bounds in the right-hand sides of \eqref{eq:26.10.22a3} and \eqref{eq:26.10.22a4}
coincide. Furthermore, by Part~2 of Item~(d) in Proposition~\ref{prop. chromatic numbers},
they also offer an improvement over the simple lower bound that is equal to the product of
the clique numbers of each factor of the strong product.
The next example shows that such an improvement does not necessarily take place if the factors
of the strong product of the regular graphs are not vertex-transitive and edge-transitive,
while also not being strongly regular.

\vspace*{0.1cm}
\begin{example}
\label{example: Frucht}
Let $\Gr{G}$ be the Frucht graph, which is 3-regular on 12~vertices
(i.e., $d=3$ and $n=12$). This graph is not strongly regular,
and also not vertex-transitive or edge-transitive.
The clique and chromatic numbers of this graph are both equal to~3
(i.e., $\clnum{\Gr{G}} = \chrnum{\Gr{G}} = 3$), so
\begin{IEEEeqnarray}{rCl}
\label{eq1: chrnum Frucht}
\chrnum{\Gr{G}^{\boxtimes \, k}} = 3^k, \quad k \in \naturals.
\end{IEEEeqnarray}
The smallest eigenvalue of the adjacency matrix of $\Gr{G}$ is equal to
$\Eigval{\min}{\Gr{G}} = -2.33866$. By \eqref{eq:26.10.22a4}, for all $k \in \naturals$,
\begin{IEEEeqnarray}{rCl}
\label{eq2: chrnum Frucht}
\chrnum{\Gr{G}^{\boxtimes \, k}} & \geq & \biggl( 1 - \frac{d}{\Eigval{\min}{\Gr{G}}} \biggr)^k \\[0.1cm]
\label{eq3: chrnum Frucht}
& = & 2.28278^k.
\end{IEEEeqnarray}
This shows that the lower bound on the chromatic number in the right-hand side of \eqref{eq:26.10.22a4}
may be looser than the simple lower bound that is equal to the product of the clique numbers of each factor.
This may happen if some of the factors, which are regular graphs, are not either strongly regular,
or both vertex-transitive and edge-transitive.
\end{example}

\subsection{The Shannon capacity of strongly regular graphs}
\label{subsection: graph capacities - s.r.g}

The Lov\'{a}sz $\theta$-function of a strongly regular graph is expressed in
closed-form in Corollary~\ref{cor4: Lovasz number for srg}, being also an upper bound on the
Shannon capacity of such a graph \cite[Theorem~1]{Lovasz79_IT}.
The independence number of a graph is, on the other hand, a lower bound on the Shannon
capacity. The Shannon capacity of such a graph is therefore determined if these
upper and lower bounds coincide. The following examples show such cases of coincidence
for strongly regular graphs. The first one reproduces (in a different way) the
Shannon capacity of the Petersen graph, a result dating back to Lov\'{a}sz \cite{Lovasz79_IT}.
The rest of the examples provide new results that determine the Shannon capacity of some
strongly regular graphs.

\vspace*{0.1cm}
\begin{example}
\label{example: Petersen}
Let $\Gr{G}$ be the Petersen graph, whose Shannon capacity is equal to~4.
It is a special case of \cite[Theorem~13]{Lovasz79_IT},
which gives the Shannon capacity of Kneser graphs
\begin{IEEEeqnarray}{rCl}
\label{eq: cap. Kneser graph}
\Theta\bigl( \KG{m}{r} \bigr) = \binom{m-1}{r-1}, \quad m \geq 2r.
\end{IEEEeqnarray}
Then, the capacity of the Petersen graph is obtained in \cite[Corollary~6]{Lovasz79_IT}
by viewing it as an isomorphic graph to the Kneser graph $\KG{5}{2}$.

As an alternative way to determine its Shannon capacity, the Petersen graph is
a strongly regular graph with parameters $\SRG(10,3,0,1)$ \cite[Section~10.3]{BrouwerM22}.
By Corollary~\ref{cor4: Lovasz number for srg}, it follows that $\theta(\Gr{G}) = 4$. This
value coincides with the independence number $\indnum{\Gr{G}} = 4$, so $\Theta(\Gr{G})=4$.
\end{example}

\vspace*{0.1cm}
\begin{example}
\label{example: Shrikhande}
Let $\Gr{G}$ be the Shrikhande graph, which is strongly regular with the parameters $\SRG(16,6,2,2)$
\cite[Section~10.6]{BrouwerM22}.
By Corollary~\ref{cor4: Lovasz number for srg}, it can be verified that
$\theta(\Gr{G}) = 4$.
Its chromatic number is $\chrnum{\Gr{G}} = 4$, and its
independence number is $\indnum{\Gr{G}} = 4$ (so
$\indnum{\Gr{G}} \chrnum{\Gr{G}} = n$, which implies that the~$n=16$
vertices in $\Gr{G}$ are partitioned into four color classes, where each
color class is a largest independent set in $\Gr{G}$ of size~4).
Hence, the graph capacity is equal to $\Theta(\Gr{G}) = 4$.

The capacity of every graph $\Gr{G}$ of order $n$ that is self-complementary and vertex-transitive
is $\Theta(\Gr{G}) = \sqrt{n}$ (see \cite[Theorem~12]{Lovasz79_IT}). It should be noted,
however, that although the capacity of the Shrikhande graph is a case where
$\Theta(\Gr{G}) = \sqrt{n}$ (with $n=16$), it does not follow from \cite[Theorem~12]{Lovasz79_IT}
since the Shrikhande graph is vertex-transitive, but it is {\em not} self-complementary.
\end{example}

\vspace*{0.1cm}
\begin{example}
\label{example: Hall-Janko}
Consider the Hall-Janko graph $\Gr{G}$, which is strongly regular with parameters
$\SRG(100, 36, 14, 12)$ \cite[Section~10.32]{BrouwerM22}. By Corollary~\ref{cor4: Lovasz number for srg},
$\theta(\Gr{G}) = 10$. The independence number of this graph is $\indnum{\Gr{G}} = 10$ (this
was obtained by the SageMath software \cite{SageMath}). Hence, the Shannon capacity
of the Hall-Janko graph is $\Theta(\Gr{G}) = 10$. This graph is vertex-transitive
but it is not self-complementary so, similarly to Example~\ref{example: Shrikhande},
this result does not follow from the capacity result in \cite[Theorem~12]{Lovasz79_IT}
for finite graphs that are vertex-transitive and self-complementary.
\end{example}

\vspace*{0.1cm}
\begin{example}
\label{example: Hoffman-Singleton}
The Hoffman-Singleton graph $\Gr{G}$ is a strongly regular graph
with parameters $\SRG(50, 7, 0, 1)$ \cite[Section~10.19]{BrouwerM22}.
By Corollary~\ref{cor4: Lovasz number for srg}, it follows that
$\theta(\Gr{G}) = 15$. The independence number of this graph is
$\indnum{\Gr{G}} = 15$, so the Shannon capacity of
this graph is $\Theta(\Gr{G}) = 15$.
\end{example}

\vspace*{0.1cm}
\begin{example}
\label{example: Schlaefli}
The Schl\"{a}fli graph $\Gr{G}$ is strongly regular with
parameters $\SRG(27,16,10,8)$ \cite[Section~10.10]{BrouwerM22}.
By Corollary~\ref{cor4: Lovasz number for srg}, it follows that
$\theta(\Gr{G}) = 3$. The independence number of this graph is
$\indnum{\Gr{G}} = 3$, so the Shannon capacity of
this graph is $\Theta(\Gr{G}) = 3$.
\end{example}

\vspace*{0.1cm}
\begin{example}
\label{example: Sims-Gewirtz}
The Sims-Gewirtz graph $\Gr{G}$ is strongly regular with parameters
$\SRG(56, 10, 0, 2)$ (a.k.a. the Gewirtz graph) \cite[Section~10.20]{BrouwerM22}.
By Corollary~\ref{cor4: Lovasz number for srg}, $\theta(\Gr{G}) = 16$, and
also $\indnum{\Gr{G}} = 16$ (it can be verified by the SageMath software \cite{SageMath}).
The Shannon capacity of the Sims-Gewirtz graph is therefore equal to
$\Theta(\Gr{G}) = 16$.
\end{example}

\vspace*{0.1cm}
\begin{example}
\label{example: M22}
The $M_{22}$ graph (a.k.a. Mesner graph) $\Gr{G}$ is strongly regular
with parameters $\SRG(77, 16, 0, 4)$ \cite[Section~10.27]{BrouwerM22}.
Its independence number is $\indnum{\Gr{G}} = 21$, and by
Corollary~\ref{cor4: Lovasz number for srg} also $\theta(\Gr{G}) = 21$.
Consequently, the Shannon capacity of the $M_{22}$ graph is $\Theta(\Gr{G}) = 21$.
\end{example}

\vspace*{0.1cm}
\begin{example}
\label{example: Cameron}
The Cameron graph $\Gr{G}$ is strongly regular with parameters $\SRG(231, 30, 9, 3)$
(see \cite[Section~10.54]{BrouwerM22}). By Corollary~\ref{cor4: Lovasz number for srg}, $\theta(\Gr{G}) = 21$,
which coincides with the independence number of $\Gr{G}$, i.e., $\indnum{\Gr{G}} = 21$.
The Shannon capacity of the Cameron graph is, hence, equal to $\Theta(\Gr{G}) = 21$.
\end{example}

\vspace*{0.1cm}
\begin{example}
\label{example: Chang - capacity}
Consider the three non-isomorphic Chang graphs. These are strongly regular graphs
with the same set of parameters $\SRG(28,12,6,4)$ \cite[Section~10.11]{BrouwerM22}.
By Corollary~\ref{cor4: Lovasz number for srg}, the Lov\'{a}sz $\theta$-function
of these graphs is equal to~4, which coincides with their independence number.
The Shannon capacity of each of these three graphs is therefore equal to~4.
\end{example}

\vspace*{0.1cm}
\begin{remark}
Examples~\ref{example: Petersen}--\ref{example: Chang - capacity} provide several
strongly regular graphs for which their Lov\'{a}sz $\theta$-function coincides
with their independence number, thus determining the Shannon capacity of these
graphs. This, however, is not the case in general for strongly regular
graphs. For example, in continuation to Example~\ref{example: Schlaefli}, the
Lov\'{a}sz $\theta$-function of the complement of the Schl\"{a}fli graph $\CGr{G}$
(it is a strongly regular graph $\SRG(27,10,1,5)$) is equal to $\theta(\CGr{G}) = 9$.
Indeed, $\theta(\Gr{G}) \, \theta(\CGr{G}) = 27$ is the order of $\Gr{G}$.
It was proved, however, by the rank-bound of Haemers \cite{Haemers79} that the
capacity of the complement of the Schl\"{a}fli graph is at most~7. It can be also
verified that the independence number of $\CGr{G}$ (i.e., the clique number of the
Schl\"{a}fli graph $\Gr{G}$) is equal to~6. This overall gives that
$6 = \clnum{\Gr{G}} \leq \Theta(\CGr{G}) \leq 7 < 9 = \theta(\CGr{G})$.
\end{remark}

We end this section by studying some capacity results of affine polar graphs.

\vspace*{0.1cm}
\begin{example}
\label{example: affine polar graphs}
{\em Affine polar graphs} (see, e.g., \cite[Section~3.3]{BrouwerM22})
are constructed by considering a vector space $\mathsf{V}$ of dimension $d$
over a finite field $\mathbb{F}_q$ with $q$ elements, equipped with a non-degenerate
quadratic form $Q$. The vertices in the graph are represented by the vectors in
$\mathsf{V}$, and any two vectors $u$ and $v$ (i.e., a pair of vertices in the graph)
are adjacent if $Q(u-v) = 0$. The resulting graph is denoted by
$\mathsf{VO}^{+}(d,q)$, $\mathsf{VO}^{-}(d,q)$, and $\mathsf{VO}(d,q)$ when the quadratic
form $Q$ is hyperbolic, elliptic or parabolic, respectively. In the first two cases,
$d$ is even, and in the third case, $d$ is odd; in all cases, $q$ is a positive integral
power of a prime number (as it is the cardinality of the finite field $\mathbb{F}_q$).
Under these conditions on $d$ and $q$, the graphs $\mathsf{VO}^{+}(d,q)$ and
$\mathsf{VO}^{-}(d,q)$ are strongly regular. Let $d = 2e$ with $e \in \naturals$.
The parameters of $\Gr{G}^{+} = \mathsf{VO}^{+}(2e,q)$, which is a strongly regular graph
$\SRG(n^{+}, d^{+}, \lambda^{+}, \mu^{+})$, are given by (see \cite[Section~3.3]{BrouwerM22})
\begin{IEEEeqnarray}{rCl}
\label{eq:06.11.22a1}
&& n^{+} = q^{2e}, \\
\label{eq:06.11.22a2}
&& d^{+} = (q^{e-1}+1) (q^e-1), \\
\label{eq:06.11.22a3}
&& \lambda^{+} = q(q^{e-2}+1)(q^{e-1}-1)+q-2, \\
\label{eq:06.11.22a4}
&& \mu^{+} = q^{e-1} (q^{e-1}+1), \\
\label{eq:06.11.22b1}
&& \Eigval{2}{\Gr{G}^{+}} = q^{e}-q^{e-1}-1,\\
\label{eq:06.11.22b2}
&& \Eigval{n^{+}}{\Gr{G}^{+}} = -q^{e-1}-1,
\end{IEEEeqnarray}
and $\Gr{G}^{-} = \mathsf{VO}^{-}(2e,q)$, a strongly regular graph
$\SRG(n^{-}, d^{-}, \lambda^{-}, \mu^{-})$, has the parameters
\begin{IEEEeqnarray}{rCl}
\label{eq:06.11.22a5}
&& n^{-} = q^{2e}, \\
\label{eq:06.11.22a6}
&& d^{-} = (q^{e-1}-1) (q^e+1), \\
\label{eq:06.11.22a7}
&& \lambda^{-} = q(q^{e-2}-1)(q^{e-1}+1)+q-2, \\
\label{eq:06.11.22a8}
&& \mu^{-} = q^{e-1} (q^{e-1}-1), \\
\label{eq:06.11.22b3}
&& \Eigval{2}{\Gr{G}^{-}} = q^{e-1}-1, \\
\label{eq:06.11.22b4}
&& \Eigval{n^{-}}{\Gr{G}^{-}} = -q^{e}+q^{e-1}-1.
\end{IEEEeqnarray}
By Corollary~\ref{cor4: Lovasz number for srg}, combined with the parameters (and eigenvalues) in
\eqref{eq:06.11.22a1}--\eqref{eq:06.11.22b4}, it follows that
\begin{IEEEeqnarray}{rCl}
\label{eq:06.11.22a9}
&& \theta(\Gr{G}^{+}) = q^{e} = \sqrt{n^{+}}, \\
\label{eq:06.11.22a10}
&& \theta(\overline{\Gr{G}^{+}}) = q^{e}, \\
\label{eq:06.11.22a11}
&& \theta(\Gr{G}^{-}) = q(q^{e}-q^{e-1}+1), \\
\label{eq:06.11.22a12}
&& \theta(\overline{\Gr{G}^{-}}) = \frac{q^{2e-1}}{q^{e}-q^{e-1}+1}.
\end{IEEEeqnarray}
Consequently, it can be verified that for some of these affine polar graphs (with the free parameters
$q$ and $e$ as above), their Lov\'{a}sz $\theta$-function
and independence number coincide. (The numerical computations of the independence numbers of these
graphs were performed by the SageMath software \cite{SageMath}.) This gives, for example, that
\begin{IEEEeqnarray}{rCl}
\label{16.10.22c1}
&& \Theta\bigl(\mathsf{VO}^{+}(4,2)\bigr) = 4, \qquad \Theta\bigl(\mathsf{VO}^{+}(6,2)\bigr) = 8,   \\
\label{16.10.22c2}
&& \Theta\bigl(\mathsf{VO}^{+}(4,3)\bigr) = 9, \qquad \Theta\bigl(\mathsf{VO}^{+}(6,3)\bigr) = 27,
\end{IEEEeqnarray}
as the exact values of the Shannon capacities of these strongly regular graphs, whose parameters
are, respectively, $\SRG(16,9,4,6)$, $\SRG(64,35,18,20)$, $\SRG(81,32,13,12)$, and
$\SRG(729,260,97,90)$.
\end{example}

\section{Proofs}
\label{section: proofs}

\subsection{Proofs for Section~\ref{subsection: bounds on theta}}
\label{subsection: proofs of bounds on theta}

\subsubsection{Proof of Proposition~\ref{prop1: bounds on theta}}
\label{subsubsection: proof of prop1}
Let $\Gr{G}$ be a $d$-regular graph on $n$ vertices. The rightmost inequality in \eqref{eq:21.10.22a1}
is provided in \cite[Theorem~9]{Lovasz79_IT}, together with a sufficient condition that it holds with equality
if $\Gr{G}$ is edge-transitive. Another sufficient condition for that inequality to hold with equality is
obtained later in this proof.

We first prove the leftmost inequality in \eqref{eq:21.10.22a1}, and also obtain sufficient conditions
that it holds with equality. By \eqref{eq: Lovasz79 - Corollary 2} (see \cite[Corollary~2]{Lovasz79_IT}),
\begin{IEEEeqnarray}{rCl}
\label{eq: 23.10.22a1}
\theta(\Gr{G}) \geq \frac{n}{\theta(\CGr{G})}
\end{IEEEeqnarray}
and, by \cite[Theorem~8]{Lovasz79_IT}, equality holds in \eqref{eq: 23.10.22a1} if $\Gr{G}$
is vertex-transitive (or, equivalently, if $\CGr{G}$ is vertex-transitive). Since $\CGr{G}$
is an $(n-d-1)$-regular graph of order $n$, an application of \eqref{eq: Lovasz79 - Theorem 9}
(see \cite[Theorem~9]{Lovasz79_IT}) to $\CGr{G}$ gives
\begin{IEEEeqnarray}{rCl}
\label{eq: 23.10.22a2}
\theta(\CGr{G}) \leq -\frac{n \Eigval{n}{\CGr{G}}}{(n-d-1)-\Eigval{n}{\CGr{G}}},
\end{IEEEeqnarray}
with equality in \eqref{eq: 23.10.22a2} if $\CGr{G}$ is edge-transitive. Combining \eqref{eq: 23.10.22a1}
and \eqref{eq: 23.10.22a2} gives
\begin{IEEEeqnarray}{rCl}
\label{eq: 23.10.22a3}
\theta(\Gr{G}) \geq 1 - \frac{n-d-1}{\Eigval{n}{\CGr{G}}},
\end{IEEEeqnarray}
with equality in \eqref{eq: 23.10.22a3} if $\CGr{G}$ is
vertex-transitive and edge-transitive (recall that a vertex-transitive
graph is regular). By the regularity of $\Gr{G}$, we have from \eqref{eq: 21.11.2022a3},
\begin{IEEEeqnarray}{rCl}
\label{eq: 23.10.22a4}
\Eigval{n}{\CGr{G}} = -1-\Eigval{2}{\Gr{G}}.
\end{IEEEeqnarray}
Combining \eqref{eq: 23.10.22a3} and \eqref{eq: 23.10.22a4} gives the leftmost
inequality in \eqref{eq:21.10.22a1}, together with the conclusion that it
holds with equality if $\CGr{G}$ is both vertex-transitive and edge-transitive.
Next, \eqref{eq:21.10.22a2} is obtained from \eqref{eq:21.10.22a1} by replacing
$\Gr{G}$ with $\CGr{G}$ (so $d$ is replaced by $(n-d-1)$ since $\Gr{G}$ is $d$-regular,
and $\CGr{G}$ is $(n-d-1)$-regular), and by relying on equality \eqref{eq: 23.10.22a4}.
A sufficient condition that the leftmost inequality in \eqref{eq:21.10.22a2}
holds with equality is therefore the same condition that the leftmost inequality
in \eqref{eq:21.10.22a1} holds with equality, while replacing $\Gr{G}$ with $\CGr{G}$;
this means that the leftmost inequality in \eqref{eq:21.10.22a2} holds with
equality if $\Gr{G}$ is vertex-transitive and edge-transitive. Likewise, the rightmost
inequality in \eqref{eq:21.10.22a2} holds with equality under the same condition
that the rightmost inequality in \eqref{eq:21.10.22a1} holds with equality, while
the graph $\Gr{G}$ in \eqref{eq:21.10.22a1} is replaced by its complement $\CGr{G}$ in
\eqref{eq:21.10.22a2}. This means that the rightmost inequality in \eqref{eq:21.10.22a1}
holds with equality if $\CGr{G}$ is edge-transitive (recall that, unlike vertex-transitivity
that is a property of $\Gr{G}$ if and only if it is a property of its complement $\CGr{G}$,
this is not the case for edge-transitivity).

We finally prove that if $\Gr{G}$ is a strongly regular graph, then the four inequalities
in \eqref{eq:21.10.22a1} and \eqref{eq:21.10.22a2} hold with equalities. Let $\Gr{G}$
be a strongly regular graph $\SRG(n,d,\lambda,\mu)$. Then, the largest,
second-largest, and smallest eigenvalues of $\Gr{G}$ are given by
\begin{IEEEeqnarray}{rCl}
\label{eq: lambda_1 srg}
&& \Eigval{1}{\Gr{G}} = d, \\[0.1cm]
\label{eq: lambda_2 srg}
&& \Eigval{2}{\Gr{G}} = r \triangleq \tfrac12 \, \Bigl[ \lambda - \mu + \sqrt{ (\lambda-\mu)^2 + 4(d-\mu) } \Bigr],  \\[0.1cm]
\label{eq: lambda_n srg}
&& \Eigval{n}{\Gr{G}} = s \triangleq \tfrac12 \, \Bigl[ \lambda - \mu - \sqrt{ (\lambda-\mu)^2 + 4(d-\mu) } \Bigr],
\end{IEEEeqnarray}
where $p_{1,2}$ in \eqref{eq: eigs s.r.g.} are replaced here with $r$ and $s$ in \eqref{eq: lambda_2 srg}
and \eqref{eq: lambda_n srg}, respectively. (Recall that if the strongly regular graph $\Gr{G}$ is disconnected,
then $\mu=0$ and $\lambda = d-1$ by \eqref{eq: relation pars. SRG}, which then indeed gives that
$\Eigval{2}{\Gr{G}} = d = \Eigval{1}{\Gr{G}}$).
The substitution of \eqref{eq: lambda_2 srg} and \eqref{eq: lambda_n srg} into \eqref{eq:21.10.22a1} and
\eqref{eq:21.10.22a2} gives
\begin{IEEEeqnarray}{rCl}
\label{eq:22.10.22a1}
\frac{n-d+r}{1+r} \leq \theta(\Gr{G}) \leq -\frac{ns}{d-s},
\end{IEEEeqnarray}
and
\begin{IEEEeqnarray}{rCl}
\label{eq:22.10.22a2}
1 - \frac{d}{s} \leq \theta(\CGr{G}) \leq \frac{n(1+r)}{n-d+r}.
\end{IEEEeqnarray}
We show as follows that the rightmost and leftmost terms in \eqref{eq:22.10.22a1} coincide, turning the
two inequalities in \eqref{eq:22.10.22a1} into equalities. The coincidence of the rightmost and leftmost
terms in \eqref{eq:22.10.22a1} is equivalent to proving the equality
\begin{IEEEeqnarray}{rCl}
\label{eq:30.10.22b1}
(n-d+r) \, (d-s) + ns \, (1+r) = 0.
\end{IEEEeqnarray}
Indeed, we get
\begin{IEEEeqnarray}{rCl}
&& (n-d+r) \, (d-s) + ns \, (1+r) \nonumber \\
\label{eq:30.10.22b2}
&& \hspace*{0.3cm} = nd - d^2 + d \, (r+s) + (n-1) rs \\
\label{eq:30.10.22b3}
&& \hspace*{0.3cm} = nd - d^2 + d \, (\lambda - \mu) + (n-1) \, (\mu-d) \\
\label{eq:30.10.22b4}
&& \hspace*{0.3cm} = -d \, (d-\lambda-1) + \mu \, (n-d-1) \\
\label{eq:30.10.22b5}
&& \hspace*{0.3cm} = 0,
\end{IEEEeqnarray}
where \eqref{eq:30.10.22b2} and \eqref{eq:30.10.22b4} hold by straightforward algebra;
\eqref{eq:30.10.22b3} holds by \eqref{eq: lambda_2 srg} and \eqref{eq: lambda_n srg};
\eqref{eq:30.10.22b5} holds by the identity in \eqref{eq: relation pars. SRG}.

The same conclusion also holds for the rightmost and leftmost terms in \eqref{eq:22.10.22a2},
which can be shown to coincide if $\Gr{G}$ is a strongly regular graph, thus turning the
two inequalities in \eqref{eq:22.10.22a2} into equalities. It holds in a similar way to the previous
part of the proof since the complement of a strongly regular graph is also strongly regular, and the second-largest
and smallest eigenvalues of $\Gr{G}$ in \eqref{eq:22.10.22a1} (i.e., respectively, $r$ and $s$) are replaced
in \eqref{eq:22.10.22a1} by the second-largest and smallest eigenvalues of $\CGr{G}$ (i.e., respectively,
$-1-s$ and $-1-r$). This proves that each of the four inequalities in \eqref{eq:21.10.22a1}
and \eqref{eq:21.10.22a2} holds with equality if $\Gr{G}$ is a strongly regular graph.

\vspace*{0.1cm}
\subsubsection{Proof of Corollary~\ref{cor4: Lovasz number for srg}}
\label{subsubsection: proof of cor4}

Let $\Gr{G}$ be a strongly regular graph $\SRG(n,d,\lambda,\mu)$.
By Proposition~\ref{prop1: bounds on theta},
\begin{IEEEeqnarray}{rCl}
\label{eq:30.10.22 a3.1}
\theta(\Gr{G}) &=& \frac{n-d+r}{1+r} \\[0.1cm]
\label{eq:30.10.22 a3.2}
&=& -\frac{ns}{d-s},
\end{IEEEeqnarray}
and
\begin{IEEEeqnarray}{rCl}
\label{eq: 30.10.22a4.1}
\theta(\CGr{G}) &=& 1 - \frac{d}{s} \\[0.1cm]
\label{eq: 30.10.22 a4.2}
&=& \frac{n(1+r)}{n-d+r},
\end{IEEEeqnarray}
where $r$ and $s$ are given in \eqref{eq: lambda_2 srg} and \eqref{eq: lambda_n srg}, respectively.
Their substitution into the right-hand sides of
\eqref{eq:30.10.22 a3.2} and \eqref{eq: 30.10.22a4.1} gives
\eqref{eq:30.10.22a1} and \eqref{eq:30.10.22a2}, respectively,
with the auxiliary parameter $t$ as defined in \eqref{eq: t - srg}.
Equality \eqref{eq: Lovasz equality for srg} is trivial by
\eqref{eq:30.10.22 a3.1}--\eqref{eq: 30.10.22 a4.2}.
This verifies \eqref{eq:30.10.22a1}--\eqref{eq:30.10.22a2}
for all strongly regular graphs.

Since the multiplicities of the distinct eigenvalues of a strongly regular
graph $\Gr{G}$ are integers, it follows from \eqref{eq: multiplicities s.r.g}
that if $2d + (n-1) (\lambda-\mu) \neq 0$, then $t$ in \eqref{eq: t - srg}
should be an integer. This implies that, under the latter condition, $\theta(\Gr{G})$
and $\theta(\CGr{G})$ are necessarily rational numbers.

\subsection{Proofs for Section~\ref{subsection: EIGs and Ramanujan graphs}}
\label{subsection: proofs of EIGs ineq., and Raman. Graphs}

\subsubsection{Proof of Corollary~\ref{cor0: ineq. for graph eigvals}}
\label{subsubsection: proof of cor0}

Let $\Gr{G}$ be a $d$-regular graph of order $n$, which is non-complete and also non-empty.
Straightforward algebra gives \eqref{eq1: cor0} by combining the rightmost
and leftmost inequalities in \eqref{eq:21.10.22a1}. Inequalities~\eqref{eq1: cor0}
and~\eqref{eq2: cor0} can be verified to be equivalent
by relying on the inequality $d+(n-1) \Eigval{n}{\Gr{G}} < 0$
(this holds since $\Eigval{n}{\Gr{G}} \leq -1$ and $d < n-1$ for a non-complete
and non-empty $d$-regular graph $\Gr{G}$ on $n$ vertices).
Equality holds in \eqref{eq1: cor0} if and only if the rightmost and leftmost
terms in \eqref{eq:21.10.22a1} are equal to $\theta(\Gr{G})$ (i.e., if the two
inequalities in \eqref{eq:21.10.22a1} are satisfied with equalities).
According to Item~(a) of Proposition~\ref{prop1: bounds on theta}, equality
holds in \eqref{eq1: cor0} and \eqref{eq2: cor0} if $\Gr{G}$ is a strongly regular graph.

We next show that the latter sufficient condition is also necessary.
To that end, we provide a second proof that relies on linear algebra.
Let $\A = \A(\Gr{G})$ and $\overline{\A} = \A(\CGr{G})$ be the adjacency matrices
of $\Gr{G}$ and $\CGr{G}$, respectively. Let $\J{n}$ and $\I{n}$ denote the $n$-times-$n$
all-ones and identity matrices, respectively. Then, by \eqref{eq:adjacency matrices},
$\overline{\A} = \J{n} - \I{n} - \A$.
Let ${\bf{j}}_n$ be the all-ones $n$-length
column vector, so ${\bf{j}}_n \, {\bf{j}}_n^{\mathrm{T}} = \J{n}$, and
\begin{IEEEeqnarray}{rCl}
\label{eq:051122a1}
\overline{\A} + c_1 \I{n} + c_2 \J{n}
= (c_1 - 1) \I{n} + (c_2 + 1) \J{n} - \A,
\end{IEEEeqnarray}
for all $c_1, c_2 \in \Reals$.
Since $\Gr{G}$ is a $d$-regular graph, the largest eigenvalue of $\A$ is equal to $d$
with ${\bf{j}}_n$ as an eigenvector. The vector ${\bf{j}}_n$ is also an eigenvector of
$\J{n}$ and $\I{n}$ with eigenvalues $n$ and $1$, respectively. By \eqref{eq:051122a1}, and
since $\Gr{G}$ is $d$-regular,
\begin{IEEEeqnarray}{rCl}
\label{eq:051122a2}
(\overline{\A} + c_1 \I{n} + c_2 \J{n}) \, {\bf{j}}_n = \bigl[ c_1-1 + n(c_2+1) - d \bigr] \, {\bf{j}}_n,
\end{IEEEeqnarray}
which means that ${\bf{j}}_n$ is an eigenvector of the matrix $\overline{\A} + c_1 \I{n} + c_2 \J{n}$. The
other $n-1$ eigenvectors of the symmetric matrix $\overline{\A} + c_1 \I{n} + c_2 \J{n}$ can be made orthogonal
to the eigenvector ${\bf{j}}_n$ (since eigenvectors of a symmetric matrix, which correspond to distinct eigenvalues,
are orthogonal with respect to the standard inner product in $\Reals^n$; furthermore, eigenvectors
that correspond to the same eigenvalue can be made orthogonal by the Gram-Schmidt procedure).
Let ${\bf{v}}$ be such an eigenvector, different from ${\bf{j}}_n$.
Then $\J{n} {\bf{v}} = {\bf{j}}_n ( {\bf{j}}_n^{\mathrm{T}} {\bf{v}}) = {\bf{0}}$, so by \eqref{eq:051122a1}
\begin{IEEEeqnarray}{rCl}
\label{eq:051122a3}
(\overline{\A} + c_1 \I{n} + c_2 \J{n}) \, {\bf{v}} = (c_1-1) {\bf{v}} - \A {\bf{v}},
\end{IEEEeqnarray}
which means that ${\bf{v}}$ is also an eigenvector of the adjacency matrix $\A$ (since, by assumption,
$\bf{v}$ is an eigenvector of the matrix $\overline{\A} + c_1 \I{n} + c_2 \J{n}$ in the left-hand side of
\eqref{eq:051122a3}). The equality $\A {\bf{v}} = \lambda {\bf{v}}$ holds, and
$\lambda \in \{\Eigval{2}{\Gr{G}}, \ldots, \Eigval{n}{\Gr{G}}\}$ with
$\Eigval{2}{\Gr{G}}$ being the largest among them. This gives that $\lambda \leq \Eigval{2}{\Gr{G}}$.
The symmetric matrix $\overline{\A} + c_1 \I{n} + c_2 \J{n}$ is positive semi-definite if
and only if all its eigenvalues are nonnegative, i.e. (see \eqref{eq:051122a2} and \eqref{eq:051122a3}),
\begin{IEEEeqnarray}{rCl}
\label{eq:051122a4}
&& c_1-1 + n(c_2+1) - d \geq 0, \\
\label{eq:051122a5}
&& c_1-1 - \Eigval{2}{\Gr{G}} \geq 0.
\end{IEEEeqnarray}
Setting the two conditions in \eqref{eq:051122a4} and \eqref{eq:051122a5} to be satisfied with equalities gives
\begin{IEEEeqnarray}{rCl}
\label{eq:051122a6}
c_1 = 1 + \Eigval{2}{\Gr{G}}, \qquad c_2 = -\frac{n-d+\Eigval{2}{\Gr{G}}}{n},
\end{IEEEeqnarray}
which implies that (see the left-hand side of \eqref{eq:051122a1})
\begin{IEEEeqnarray}{rCl}
\label{eq:31.10.2022a1}
\overline{\A} + \bigl(1 + \Eigval{2}{\Gr{G}} \bigr) \, \I{n}
- \biggl( \frac{n-d+\Eigval{2}{\Gr{G}}}{n} \biggr) \, \J{n} \succeq 0
\end{IEEEeqnarray}
is a positive semi-definite matrix. Clearly, also
\begin{IEEEeqnarray}{rCl}
\label{eq:31.10.2022a2}
\A - \Eigval{n}{\Gr{G}} \I{n} \succeq 0
\end{IEEEeqnarray}
is positive semi-definite.
The trace of a product of two $n$-times-$n$ positive semi-definite matrices is nonnegative,
so it follows from \eqref{eq:31.10.2022a1} and \eqref{eq:31.10.2022a2} that
\begin{IEEEeqnarray}{rCl}
\label{eq:31.10.2022a3}
\tr \Biggl( \biggl( \overline{\A} + \bigl(1 + \Eigval{2}{\Gr{G}} \bigr) \, \I{n}
- \biggl( \frac{n-d+\Eigval{2}{\Gr{G}}}{n} \biggr) \, \J{n} \biggr) \,
\bigl( \A - \Eigval{n}{\Gr{G}} \I{n} \bigr) \Biggr) \geq 0.
\end{IEEEeqnarray}
For an adjacency matrix $\A$ of a $d$-regular graph $\Gr{G}$,
\begin{IEEEeqnarray}{rCl}
\label{eq:31.10.2022a4}
&& \tr( \overline{\A} \A ) = \tr(\A) = \tr(\overline{\A}) = 0,  \\
\label{eq:31.10.2022a5}
&& \tr(\J{n} \A) = nd,
\end{IEEEeqnarray}
and
\begin{IEEEeqnarray}{rCl}
\label{eq:31.10.2022a6}
\hspace*{-1.6cm} \tr(\I{n}) = \tr(\J{n}) = n.
\end{IEEEeqnarray}
By \eqref{eq:31.10.2022a4}--\eqref{eq:31.10.2022a6}, expanding the left-hand side in \eqref{eq:31.10.2022a3}
gives the inequality
\begin{IEEEeqnarray}{rCl}
\label{eq:31.10.2022a7}
-d \, \bigl(n-d+\Eigval{2}{\Gr{G}}\bigr) - \Eigval{n}{\Gr{G}} \, \bigl(d + (n-1) \Eigval{2}{\Gr{G}} \bigr) \geq 0,
\end{IEEEeqnarray}
which proves \eqref{eq1: cor0} (we have $d + (n-1) \Eigval{2}{\Gr{G}} > 0$, by the assumption
that the $d$-regular graph $\Gr{G}$ is non-complete and non-empty, so $\Eigval{2}{\Gr{G}} \geq 0$).
This alternative proof is next used to find the necessary and sufficient condition for
the satisfiability of \eqref{eq1: cor0} with equality. An equality in \eqref{eq1: cor0}
holds if and only if \eqref{eq:31.10.2022a3} holds with equality (i.e., the trace of the
product of the two positive semi-definite matrices in the left-hand sides of
\eqref{eq:31.10.2022a1} and \eqref{eq:31.10.2022a2} is equal to zero).
This holds if and only if the column spaces of the two matrices in the left-hand sides of
\eqref{eq:31.10.2022a1} and \eqref{eq:31.10.2022a2} are orthogonal. [Clarification: if $\bf{B}, \bf{C} \succeq 0$ are
positive semi-definite matrices, then $\bf{B} = \bf{S} \bf{S}^{\mathrm{T}}$ and
$\bf{C} = \bf{Q} \bf{Q}^{\mathrm{T}}$ for some matrices $\bf{S}$ and $\bf{Q}$. Hence, under the assumption that
$\tr(\bf{B} \bf{C})=0$,
\begin{IEEEeqnarray}{rCl}
\hspace*{-1cm}
0 = \tr({\bf{B}}{\bf{C}}) = \tr({\bf{S}} {\bf{S}}^{\mathrm{T}} {\bf{Q}} {\bf{Q}}^{\mathrm{T}})
= \tr({\bf{S}}^{\mathrm{T}} {\bf{Q}} {\bf{Q}}^{\mathrm{T}} {\bf{S}})
= \tr({\bf{S}}^{\mathrm{T}} {\bf{Q}} \, ({\bf{S}}^{\mathrm{T}} {\bf{Q}})^{\mathrm{T}}) = \sum_{i=1}^n \sum_{j=1}^n
\bigl( {\bf{S}}^{\mathrm{T}} {\bf{Q}} \bigr)_{i,j}^2
\end{IEEEeqnarray}
implies that ${\bf{S}}^{\mathrm{T}} {\bf{Q}} = \bf{0}$, so also ${\bf{B}}{\bf{C}}
= {\bf{S}} ({\bf{S}}^{\mathrm{T}} {\bf{Q}}) {\bf{Q}}^{\mathrm{T}} = \bf{0}$.
Since ${\bf{B}}$ and ${\bf{C}}$ are symmetric matrices with ${\bf{B}}{\bf{C}} = \bf{0}$,
it means that the column spaces of $\bf{B}$ and $\bf{C}$ are orthogonal].
These column spaces, however, can be orthogonal only if $\A$ has no eigenvalues other than
$d$, $\Eigval{2}{\Gr{G}}$ and $\Eigval{n}{\Gr{G}}$; otherwise, there is a joint eigenvector of
$\A$ in the two column spaces of the matrices in the left-hand sides of \eqref{eq:31.10.2022a1}
and \eqref{eq:31.10.2022a2}, so these two column spaces cannot be orthogonal in the latter case.
[Clarification: suppose that $\A$ has an eigenvector ${\bf{v}}$ that corresponds to an eigenvalue
$\lambda^\ast \not\in \{d, \Eigval{2}{\Gr{G}}, \Eigval{n}{\Gr{G}}\}$. Then, by \eqref{eq:051122a3}
and the value of $c_1$ in \eqref{eq:051122a6}, the right-hand side of \eqref{eq:051122a3} is
equal to $(\Eigval{2}{\Gr{G}} - \lambda^\ast) {\bf{v}}$, and
$(\A - \Eigval{n}{\Gr{G}} \I{n}) {\bf{v}} = (\lambda^\ast - \Eigval{n}{\Gr{G}}) {\bf{v}}$.
Both coefficients of ${\bf{v}}$ in these two expressions are nonzero, so $\bf{v}$ belongs
to the two column spaces of the two matrices in the left-hand sides of \eqref{eq:31.10.2022a1}
and \eqref{eq:31.10.2022a2}. These column spaces are therefore not orthogonal under the assumption
of the existence of a distinct eigenvalue $\lambda^\ast$ of $\A$, as above].
We next distinguish between two cases in regard to the connectivity of the graph $\Gr{G}$.
\begin{enumerate}[(a)]
\item If the $d$-regular graph $\Gr{G}$ is connected, then $\Eigval{1}{\Gr{G}} = d > \Eigval{2}{\Gr{G}}$.
By assumption, $\Gr{G}$ is also non-complete and non-empty graph, so $\Eigval{2}{\Gr{G}} > \Eigval{n}{\Gr{G}}$.
The connected regular graph $\Gr{G}$ thus has exactly three distinct eigenvalues, so it is strongly regular.
\item If the $d$-regular graph $\Gr{G}$ is disconnected, then $\Eigval{1}{\Gr{G}} = d = \Eigval{2}{\Gr{G}}$.
If, by assumption, inequality \eqref{eq1: cor0} holds with equality, then $\Eigval{n}{\Gr{G}} = -1$. This means
that $\Gr{G}$ is a disjoint union of equal-sized complete graphs $\CoG{d+1}$, so it is an imprimitive strongly
regular graph (i.e., there are no common neighbors of any pair of non-adjacent vertices in $\Gr{G}$).
\end{enumerate}
We therefore conclude that, for a non-complete and non-empty $d$-regular graph on $n$ vertices, the
condition that $\Gr{G}$ is strongly regular is also necessary (and not only sufficient, as it
is shown in the first proof) for the inequality in \eqref{eq1: cor0} to hold with equality. Equivalently,
$\Gr{G}$ being strongly regular is a necessary condition for the inequality in \eqref{eq2: cor0}
to hold with equality. This completes the proof of the necessity and sufficiency of the condition
on the strong regularity of $\Gr{G}$.

\vspace*{0.1cm}
\subsubsection{Proof of Corollary~\ref{cor1: ineq. for graph eigvals}}
\label{subsubsection: proof of cor1 - eigenvals ineq.}
Let $\Gr{G}$ be a $d$-regular graph of order~$n$, which is non-complete and non-empty.

{\em Proof of Item~(a)}:
By the definition of $g_2(\Gr{G})$ and $g_n(\Gr{G})$ in \eqref{eq: 17.11.22a1},
and in light of the equalities $\Eigval{n}{\CGr{G}} = -1 - \Eigval{2}{\Gr{G}}$
and $\Eigval{2}{\CGr{G}} = -1 - \Eigval{n}{\Gr{G}}$ (see \eqref{eq: 21.11.2022a2}),
it can be readily verified that the rightmost and leftmost inequalities in
\eqref{eq: 17.11.22a2} are equivalent to \eqref{eq1: cor0} and \eqref{eq2: cor0},
respectively. Furthermore, in light of Corollary~\ref{cor0: ineq. for graph eigvals},
each inequality in \eqref{eq: 17.11.22a2} holds with equality if and only if $\Gr{G}$
is strongly regular.

{\em Proof of Item~(b)}:
Let $\Gr{G}$ be a strongly regular graph. We next prove the
claim about the dichotomy in the number of distinct values in the
sequence $\{g_\ell(\Gr{G})\}_{\ell=1}^n$.
\begin{enumerate}[(1)]
\item For a $d$-regular graph on $n$ vertices, $\Eigval{1}{\Gr{G}} = d$, and $\Eigval{1}{\CGr{G}} = n-d-1$.
Substituting these eigenvalues into \eqref{eq: 17.11.22a1} gives that $g_1(\Gr{G}) = n-d-2 \geq 0$
($\Gr{G}$ is non-complete, so $d \leq n-2$).
\item The graph $\Gr{G}$ and its complement $\CGr{G}$ are both strongly regular, so each one of them
has at most three distinct eigenvalues.
\item For a strongly regular graph, by Item~(a), $g_2(\Gr{G}) = -1 = g_n(\Gr{G})$, and $g_\ell(\Gr{G})=-1$ if either
(i)~$\Eigval{\ell}{\Gr{G}} = \Eigval{2}{\Gr{G}}$ and $\Eigval{\ell}{\CGr{G}} = \Eigval{2}{\CGr{G}}$,
or (ii)~$\Eigval{\ell}{\Gr{G}} = \Eigval{n}{\Gr{G}}$ and $\Eigval{\ell}{\CGr{G}} = \Eigval{n}{\CGr{G}}$.
\item By assumption, $\Gr{G}$ is a strongly regular graph, which implies that so is $\CGr{G}$.
Due to their regularity,
$\Eigval{1}{\Gr{G}} \geq  \Eigval{2}{\Gr{G}}$, and $\Eigval{1}{\CGr{G}} \geq  \Eigval{2}{\CGr{G}}$
with equalities, respectively, if and only if $\Gr{G}$ or $\CGr{G}$ are disconnected graphs.
The sequence $\{g_\ell(\Gr{G})\}_{\ell=1}^n$ gets an additional (third) distinct value if and only
if the multiplicities of the smallest and the second-largest eigenvalues of $\Gr{G}$ in the {\em subsequence}
$(\Eigval{2}{\Gr{G}}, \ldots, \Eigval{n}{\Gr{G}})$ are distinct.
Indeed, in the latter case, only one of the following two options is possible:
(iii) $\Eigval{\ell}{\Gr{G}} = \Eigval{2}{\Gr{G}}$ and $\Eigval{\ell}{\CGr{G}} = \Eigval{n}{\CGr{G}}$,
or (iv)~$\Eigval{\ell}{\Gr{G}} = \Eigval{n}{\Gr{G}}$ and $\Eigval{\ell}{\CGr{G}} = \Eigval{2}{\CGr{G}}$.
This holds since, by \eqref{eq: 21.11.2022a2}, the multiplicity of the second-largest eigenvalue of
$\Gr{G}$ is equal to the multiplicity of the smallest eigenvalue of $\CGr{G}$, and similarly,
the multiplicity of the smallest eigenvalue of $\Gr{G}$ is equal to the multiplicity of the
second-largest eigenvalue of $\CGr{G}$. It therefore follows that the third
distinct value (as above) is attained by the sequence $\{g_\ell\}_{\ell=1}^n$ a number of times that
is equal to the absolute value of the difference between the multiplicities of the second-largest
and the smallest eigenvalues of $\Gr{G}$ in the subsequence $(\Eigval{2}{\Gr{G}}, \ldots, \Eigval{n}{\Gr{G}})$
(provided that the latter two multiplicities are distinct).
\end{enumerate}

{\em Proof of Item~(c)}: Let $\Gr{G}$ be self-complementary and $d$-regular on $n$ vertices. Then,
\begin{IEEEeqnarray}{rCl}
\label{eq: 17.11.22a5}
d = \tfrac12 (n-1), \quad \Eigval{2}{\CGr{G}} = \Eigval{2}{\Gr{G}},
\quad \Eigval{n}{\CGr{G}} = \Eigval{n}{\Gr{G}}.
\end{IEEEeqnarray}
Combining the rightmost inequality in \eqref{eq: 17.11.22a2} and the equalities
in \eqref{eq: 17.11.22a5} readily gives
\begin{IEEEeqnarray}{rCl}
\label{eq: 17.11.22a6}
\frac{2 \Eigval{2}{\Gr{G}}^{\, 2} - \tfrac12 (n+1)}{1 + 2 \Eigval{2}{\Gr{G}}} \geq -1.
\end{IEEEeqnarray}
Since $\Gr{G}$ is non-complete and non-empty, we get $\Eigval{2}{\Gr{G}} > 0$, which then gives
from \eqref{eq: 17.11.22a6} the quadratic inequality
\begin{IEEEeqnarray}{rCl}
\label{eq: 17.11.22a7}
2 \Eigval{2}{\Gr{G}}^{\, 2} + 2 \Eigval{2}{\Gr{G}} - \tfrac12 (n-1) \geq 0.
\end{IEEEeqnarray}
Its solution gives \eqref{eq: 17.11.22a3} ($n>1$ as otherwise, $\Gr{G}=\CoG{1}$,
but $\Gr{G}$ is by assumption non-complete). Hence, it also follows that
\begin{IEEEeqnarray}{rCl}
\label{eq: 17.11.22a8}
\Eigval{n}{\Gr{G}} &=& \Eigval{n}{\CGr{G}} \\
\label{eq: 17.11.22a9}
&=& -1 - \Eigval{2}{\Gr{G}} \\
\label{eq: 17.11.22a10}
&\leq& -1 - \tfrac12 (\sqrt{n}-1) \\[0.1cm]
\label{eq: 17.11.22a11}
&=& -\tfrac12 (\sqrt{n}+1),
\end{IEEEeqnarray}
where \eqref{eq: 17.11.22a8} holds since $\Gr{G}$ is (by assumption) self-complementary;
\eqref{eq: 17.11.22a9} holds by \eqref{eq: 21.11.2022a2} since $\Gr{G}$ is regular, and
\eqref{eq: 17.11.22a10} holds by \eqref{eq: 17.11.22a3}.

{\em Proof of Item~(d)}:
If $\Gr{G}$ is self-complementary and strongly regular, then
in light of Items~(b) and~(c) here, both inequalities in \eqref{eq: 17.11.22a3}
and \eqref{eq: 17.11.22a4} hold with equality.

\vspace*{0.1cm}
\subsubsection{Proof of Corollary~\ref{cor1: clique and chromatic nums}}
\label{subsubsection: proof of cor1 - clique and chromatic nums}
Let $\{\Gr{G}_\ell\}_{\ell \in \naturals}$ be a sequence of
regular graphs where $\Gr{G}_{\ell}$ is $d_\ell$-regular of order
$n_\ell$, such that $n_\ell \to \infty$ and $\frac{d_\ell}{n_\ell} \to 0$
as we let $\ell$ tend to infinity. Then,
\begin{IEEEeqnarray}{rCl}
\label{eq: 24.10.22a1}
\limsup_{\ell \to \infty} \, \clnum{\Gr{G}_{\ell}} &\leq&
\limsup_{\ell \to \infty} \, \theta(\CGr{G}_{\ell}) \\
\label{eq: 24.10.22a2}
&\leq& \limsup_{\ell \to \infty} \,
\frac{n_\ell \bigl(1+\Eigval{2}{\Gr{G}_\ell}\bigr)}{n_\ell-d_\ell+\Eigval{2}{\Gr{G}_\ell}} \\[0.1cm]
\label{eq: 24.10.22a3}
&=& 1 + \limsup_{\ell \to \infty} \Eigval{2}{\Gr{G}_\ell},
\end{IEEEeqnarray}
where \eqref{eq: 24.10.22a1} holds by the leftmost inequality in \eqref{eq1b: sandwich};
\eqref{eq: 24.10.22a2} holds by the rightmost inequality in \eqref{eq:21.10.22a2};
\eqref{eq: 24.10.22a3} holds by the assumption that $n_\ell \to \infty$, and since
the eigenvalues of $\Gr{G}_\ell$ are bounded (in absolute value) by $d_\ell$ with
$\underset{\ell \to \infty}{\lim} \frac{d_\ell}{n_\ell}=0$.
This leads to inequality \eqref{eq:20.10.22b1}, by a floor operation in the right-hand side
of \eqref{eq: 24.10.22a3}, since clique numbers are integers.

We next prove inequality \eqref{eq:20.10.22b2}. For any graph $\Gr{G}$ on $n$ vertices,
\begin{IEEEeqnarray}{rCl}
\label{eq: well-known1}
\indnum{\Gr{G}} \chrnum{\Gr{G}} \geq n.
\end{IEEEeqnarray}
(This well-known inequality holds since
the independence number $\indnum{\Gr{G}}$ denotes the size of a largest independent set
in $\Gr{G}$, and in coloring the vertices in $\Gr{G}$ with $\chrnum{\Gr{G}}$ colors,
all color classes are independent). Additionally, $\clnum{\Gr{G}} = \indnum{\CGr{G}}$, so
\begin{IEEEeqnarray}{rCl}
\label{eq: prod}
\clnum{\Gr{G}} \chrnum{\CGr{G}} \geq n.
\end{IEEEeqnarray}
This gives
\begin{IEEEeqnarray}{rCl}
\label{eq: 24.10.22a4}
\liminf_{\ell \to \infty} \,  \frac{\chrnum{\CGr{G}_\ell}}{n_\ell} &\geq&
\liminf_{\ell \to \infty} \,  \frac1{\clnum{\Gr{G}_\ell}} \\[0.1cm]
\label{eq: 24.10.22a5}
&=& \frac1{ \underset{\ell \to \infty}{\limsup} \, \clnum{\Gr{G}_\ell}} \\
\label{eq: 24.10.22a6}
&\geq& \frac1{1 + \underset{\ell \to \infty}{\limsup} \lfloor \Eigval{2}{\Gr{G}_\ell} \rfloor} \, ,
\end{IEEEeqnarray}
where \eqref{eq: 24.10.22a4} holds by \eqref{eq: prod};
\eqref{eq: 24.10.22a5} is trivial, and \eqref{eq: 24.10.22a6} holds by \eqref{eq:20.10.22b1}.

\vspace*{0.1cm}
\subsubsection{Proof of Corollary~\ref{cor2: clique and chromatic nums}}
\label{subsubsection: proof of cor2}

Inequalities~\eqref{eq:21.10.22c1} and~\eqref{eq:21.10.22c2} readily
follow from Corollary~\ref{cor1: clique and chromatic nums} since
if $\{\Gr{G}_{\ell}\}_{\ell \in \naturals}$ is a sequence of
Ramanujan $d$-regular graphs ($d$ is a fixed degree of the vertices),
then (by definition)
\begin{IEEEeqnarray}{rCl}
\label{eq: 24.10.22a7}
\Eigval{2}{\Gr{G}_\ell} \leq 2 \sqrt{d-1},
\end{IEEEeqnarray}
for all $\ell \in \naturals$. By the assumption that the graph $\Gr{G}_\ell$
has order $n_\ell$ with $\underset{\ell \to \infty}{\lim} n_\ell = \infty$,
inequalities \eqref{eq:21.10.22c1} and \eqref{eq:21.10.22c2} are obtained
by combining, respectively, \eqref{eq:20.10.22b1} and \eqref{eq:20.10.22b2}
with \eqref{eq: 24.10.22a7}.

Inequality~\eqref{eq:21.10.22d5} is obtained by combining the leftmost
inequality in \eqref{eq:21.10.22a1} with \eqref{eq: 24.10.22a7}. Indeed,
since $| \Eigval{2}{\Gr{G}_\ell} | \leq d$ (where the degree of the vertices
of $\Gr{G}_\ell$ is, by assumption, equal to a fixed value $d$), it follows
that
\begin{IEEEeqnarray}{rCl}
\label{eq: 27.11.2022a1}
\theta(\Gr{G}_\ell) \geq \frac{n_\ell-2d}{1+\Eigval{2}{\Gr{G}_\ell}}
\geq \frac{n_\ell-2d}{1+2\sqrt{d-1}},
\end{IEEEeqnarray}
which then yields \eqref{eq:21.10.22d5}.

\vspace*{0.1cm}
\subsubsection{Proof of Corollary~\ref{cor3: clique number - Ramanujan}}
\label{subsubsection: proof of cor3}

Let $\Gr{G}$ be a Ramanujan $d$-regular graph on $n$ vertices. If $\Gr{G} = \CoG{n}$
is the complete graph, which is a Ramanujan $(n-1)$-regular graph if $n \geq 3$, then
inequality \eqref{eq:21.10.22b1} clearly holds with equality (since $\clnum{\CoG{n}} = n$).
Otherwise, if $\Gr{G}$ is non-complete, then combining the leftmost inequality in
\eqref{eq1b: sandwich} and the rightmost inequality in \eqref{eq:21.10.22a2} gives that
\begin{eqnarray}
\label{eq: 24.10.22a8}
\clnum{\Gr{G}} \leq \bigg\lfloor
\frac{n \bigl(1+\Eigval{2}{\Gr{G}}\bigr)}{n-d+\Eigval{2}{\Gr{G}}} \bigg\rfloor,
\end{eqnarray}
where the floor operation in the right-hand side of \eqref{eq: 24.10.22a8} is justified
since $\clnum{\Gr{G}}$ is an integer. Since $0 \leq \Eigval{2}{\Gr{G}} \leq 2 \sqrt{d-1}$
is satisfied for every Ramanujan $d$-regular graph $\Gr{G}$ that is non-complete,
and since the function $f_1 \colon (-1, \infty) \to (0, \infty)$ that is given by
\begin{eqnarray}
\label{eq: f1 func.}
f_1(x) \triangleq \frac{n(1+x)}{n-d+x}, \quad x > -1,
\end{eqnarray}
is monotonically increasing, inequality \eqref{eq:21.10.22b1} then follows from
\eqref{eq: 24.10.22a8} and the monotonicity of the function $f_1$.
Eq.~\eqref{eq:21.10.22b3} follows from \eqref{eq:21.10.22b1},
\eqref{eq: prod}, and since the chromatic number is an integer.

Inequality~\eqref{eq:21.10.22d3} holds with equality if $\Gr{G} = \CoG{n}$
(both sides are then equal to~1). Otherwise, \eqref{eq:21.10.22d3}
holds by the leftmost inequality in \eqref{eq:21.10.22a1}, since (as mentioned above)
$0 \leq \Eigval{2}{\Gr{G}} \leq 2 \sqrt{d-1}$ is satisfied for every Ramanujan
$d$-regular graph $\Gr{G}$ that is non-complete, and since
the function $f_2 \colon (-1, \infty) \to (0, \infty)$ that is given by
\begin{eqnarray}
\label{eq: f2 func.}
f_2(x) \triangleq \frac{n-d+x}{1+x}, \quad x > -1,
\end{eqnarray}
is monotonically decreasing.

\subsection{Proofs for Section~\ref{subsection: eigenvalues of strong products}}
\label{subsection: Proofs of eig bounds for strong products}

\subsubsection{Proof of Proposition~\ref{prop: bounds on eigvals - strong products}}
\label{subsubsection: proof of prop. 2}

Let $\Gr{G}_1, \ldots, \Gr{G}_k$ be regular graphs such that, for all $\ell \in \OneTo{k}$,
the graph $\Gr{G}_\ell$ is $d_\ell$-regular on $n_\ell$ vertices.
Let $\Gr{G} = \Gr{G}_1 \boxtimes \ldots \boxtimes \Gr{G}_k$ be the strong product of
these $k$ regular graphs. We next prove the two items of
Proposition~\ref{prop: bounds on eigvals - strong products}.
\begin{enumerate}[(a)]
\item
By the leftmost inequality in \eqref{eq:21.10.22a1},
unless $\Gr{G} = \CoG{n}$,
\begin{IEEEeqnarray}{rCl}
\label{eq:07.11.22a1}
\theta(\Gr{G}) \geq \frac{n(\Gr{G})-d(\Gr{G})+\Eigval{2}{\Gr{G}}}{1+\Eigval{2}{\Gr{G}}},
\end{IEEEeqnarray}
where $n(\Gr{G})$ and $d(\Gr{G})$ denote, respectively, the order and valency of the
strong product, which is a regular graph (since, by assumption, each factor is regular).
The following equalities hold as a result of the strong product operation:
\begin{IEEEeqnarray}{rCl}
\label{eq:07.11.22a2}
&& n(\Gr{G}) = \prod_{\ell=1}^k n_\ell, \\
\label{eq:07.11.22a3}
&& d(\Gr{G}) = \prod_{\ell=1}^k (1+d_\ell) - 1, \\
\label{eq:07.11.22a4}
&& \theta(\Gr{G}) = \prod_{\ell=1}^k \theta(\Gr{G}_\ell).
\end{IEEEeqnarray}
Indeed, equality~\eqref{eq:07.11.22a2} holds since the cardinality of a Cartesian
product of finite sets is equal to the product of the cardinalities of each set;
equality~\eqref{eq:07.11.22a3} can be justified by first verifying the special case
of a strong product of two regular graphs, and then proceeding by a mathematical
induction on $k$. Finally, equality~\eqref{eq:07.11.22a4} holds by
\eqref{eq: Lovasz79 - Theorem 7} (see \cite[Theorem~7]{Lovasz79_IT}).
Combining the bound in \eqref{eq:07.11.22a1} with equalities
\eqref{eq:07.11.22a2}--\eqref{eq:07.11.22a4} gives
\begin{IEEEeqnarray}{rCl}
\label{eq:07.11.22a5}
\prod_{\ell=1}^k \theta(\Gr{G}_\ell) \geq
1 + \frac{\overset{k}{\underset{\ell=1}{\prod}} n_\ell
-\overset{k}{\underset{\ell=1}{\prod}} (1+d_\ell)}{1+\Eigval{2}{\Gr{G}}}.
\end{IEEEeqnarray}
Unless all $\Gr{G}_\ell$ (with $\ell \in \OneTo{k}$) are complete graphs,
the left-hand side of \eqref{eq:07.11.22a5} is strictly larger than~1, and
then rearrangement of the terms in \eqref{eq:07.11.22a5} gives the lower
bound on $\Eigval{2}{\Gr{G}}$ in \eqref{eq: eig2 - LB1}. Next, the possible
loosening of the lower bound in the right-hand side of \eqref{eq: eig2 - LB1} to
the lower bound in the right-hand side of \eqref{eq: eig2 - LB2} holds by
\eqref{eq: Lovasz79 - Theorem 9} (see \cite[Theorem~9]{Lovasz79_IT}).
Inequality~\eqref{eq: eig2 - LB2} holds with equality if each regular factor
$\Gr{G}_\ell$ is either edge-transitive (by \cite[Theorem~9]{Lovasz79_IT})
or strongly regular (by Item~(a) of Proposition~\ref{prop1: bounds on theta}).

\item
Combining \eqref{eq: Lovasz79 - Theorem 9} with equalities
\eqref{eq:07.11.22a2}--\eqref{eq:07.11.22a4} gives, with
$n = n(\Gr{G})$ and $d = d(\Gr{G})$,
\begin{IEEEeqnarray}{rCl}
\label{eq:07.11.22a6}
\prod_{\ell=1}^k \theta(\Gr{G}_\ell) & = & \theta(\Gr{G}) \\
\label{eq:07.11.22a7}
& \leq & -\frac{n \Eigval{n}{\Gr{G}}}{d - \Eigval{n}{\Gr{G}}} \\
\label{eq:07.11.22a8}
& = & - \frac{\overset{k}{\underset{\ell=1}{\prod}} n_\ell \cdot
\Eigval{n}{\Gr{G}}}{ \overset{k}{\underset{\ell=1}{\prod}} (1+d_\ell) - 1 - \Eigval{n}{\Gr{G}}}.
\end{IEEEeqnarray}
Unless all $\Gr{G}_\ell$ (with $\ell \in \OneTo{k}$) are empty graphs,
the denominator in the right-hand side of \eqref{eq:07.11.22a8} is strictly positive.
This gives \eqref{eq: eig_min - UB1} after rearrangement of terms. Finally, the
transition from \eqref{eq: eig_min - UB1} to \eqref{eq: eig_min - UB2} is justified if
\begin{IEEEeqnarray}{rCl}
\label{eq:07.11.22a9}
\theta(\Gr{G}_\ell) = -\frac{n_\ell \Eigval{\min}{\Gr{G}_\ell}}{d_\ell - \Eigval{\min}{\Gr{G}_\ell}},
\qquad \forall \, \ell \in \OneTo{k}.
\end{IEEEeqnarray}
As above (the end of the proof of Item~(a)), the condition in \eqref{eq:07.11.22a9} holds
if the regular graph $\Gr{G}_\ell$ is either edge-transitive or strongly regular.
\end{enumerate}

\vspace*{0.1cm}
\subsubsection{Proof of Corollary~\ref{corollary: bounds on eigvals - strong powers}}
\label{subsubsection: proof of cor. - bounds on eigs., strong product}

Let $\Gr{G}$ a $d$-regular graph of order $n$. The lower bound on the second-largest eigenvalue
$\Eigval{2}{\Gr{G}^{\boxtimes \, k}}$ in the right-hand side of \eqref{eq: eig2 LB - graph powers}
follows from \eqref{eq: eig2 - LB1} by setting there $\Gr{G}_1, \ldots, \Gr{G}_k$ to be all
identical to $\Gr{G}$. The upper bound on the smallest eigenvalue $\Eigval{\min}{\Gr{G}^{\boxtimes \, k}}$
in the right-hand side of \eqref{eq: eig_min UB - graph powers} follows in a similar way from
\eqref{eq: eig_min - UB1}.

\vspace*{0.1cm}
\subsubsection{Proof of Proposition~\ref{prop.: not Ramanujan}}
\label{subsection: proof of proposition on non-Ramanujan graphs}

Let $\Gr{G}$ be a connected $d$-regular graph on $n$ vertices, which is non-empty and non-complete,
and let $k \in \naturals$. Then, $\Gr{G}^{\boxtimes \, k}$ is a connected regular graph,
which is non-complete and non-empty (so, its largest eigenvalue is of multiplicity~1).
By \eqref{eq: eig2 LB - graph powers},
\begin{eqnarray}
\label{eq: 09.11.22a1}
\Eigval{2}{\Gr{G}^{\boxtimes \, k}} \geq \frac{n^k - (1+d)^k}{\theta(\Gr{G})^k-1} - 1.
\end{eqnarray}
In order to prove that the $k$-fold strong power $\Gr{G}^{\boxtimes \, k}$ is non-Ramanujan,
it is sufficient to show that the lower bound on its second-largest eigenvalue in
the right-hand side of \eqref{eq: 09.11.22a1} is larger than $2 \sqrt{d_k-1}$ (see the
right-hand side of \eqref{eq:Ramanujan}); here, $d_k = (1+d)^k - 1$ is the valency of
the strong power $\Gr{G}^{\boxtimes \, k}$. Since $\Gr{G}$ is $d$-regular and non-complete
(i.e., $d < n-1$ and $\theta(\Gr{G})>1$), the right-hand side of \eqref{eq: 09.11.22a1}
scales asymptotically like $\Bigl( \frac{n}{\theta(\Gr{G})} \Bigr)^k$ (for a sufficiently
large $k$), whereas the expression $2 \sqrt{d_k-1}$ scales asymptotically like
$2 (1+d)^{\frac{k}{2}}$. Comparing these two exponents gives that if
\begin{eqnarray}
\label{eq: 09.11.22a2}
\theta(\Gr{G}) < \frac{n}{\sqrt{1+d}},
\end{eqnarray}
then the exponential growth rate (in $k$) of the right-hand side of \eqref{eq: 09.11.22a1}
is larger than that one of $2 \sqrt{d_k-1}$. Hence, for sufficiently large $k$,
the strong power $\Gr{G}^{\boxtimes \, k}$ is a (highly) non-Ramanujan graph under
the condition in \eqref{eq: 09.11.22a2}. This means that there exists
$k_0 \in \naturals$ such that the strong power $\Gr{G}^{\boxtimes \, k}$
is non-Ramanujan for all $k \geq k_0$. We next obtain an explicit value of
such $k_0$, which is not necessarily the smallest one, proving that such
a valid value for $k_0$ is given by~\eqref{eq: k_0}. To that end, based on the
above explanation, one needs to deal with the inequality
\begin{eqnarray}
\label{eq: 09.11.22a3}
\frac{n^k - (1+d)^k}{\theta(\Gr{G})^k-1} - 1 > 2 \sqrt{(1+d)^k - 2}.
\end{eqnarray}
In order to obtain a closed-form solution, we strengthen the condition in
\eqref{eq: 09.11.22a3} to
\begin{eqnarray}
\label{eq: 09.11.22a4}
\frac{n^k - (1+d)^k}{\theta(\Gr{G})^k} - 1 \geq 2 (1+d)^{\frac{k}{2}}.
\end{eqnarray}
Dividing both sides of \eqref{eq: 09.11.22a4} by $(1+d)^{\frac{k}{2}}$ gives
\begin{eqnarray}
\label{eq: 10.11.22c1}
\biggl( \frac{n}{\sqrt{1+d} \; \theta(\Gr{G})} \biggr)^k
- \biggl( \frac{\sqrt{1+d}}{\theta(\Gr{G})} \biggr)^k
- (1+d)^{-\frac{k}{2}} \geq 2.
\end{eqnarray}
Let $k \geq 3$. The condition imposed in \eqref{eq: 10.11.22c1} can be further
strengthened to
\begin{eqnarray}
\label{eq: 10.11.22c2}
\biggl( \frac{n}{\sqrt{1+d} \; \theta(\Gr{G})} \biggr)^k
- \biggl( \frac{\sqrt{1+d}}{\theta(\Gr{G})} \biggr)^k
\geq 2 + (1+d)^{-\frac{3}{2}}.
\end{eqnarray}
Since $d<n-1$ for a non-complete $d$-regular graph of order $n$, for all $k \geq 3$,
\begin{align}
\biggl( \frac{n}{\sqrt{1+d} \; \theta(\Gr{G})} \biggr)^k
- \biggl( \frac{\sqrt{1+d}}{\theta(\Gr{G})} \biggr)^k
& = \biggl( \frac{n}{\sqrt{1+d} \; \theta(\Gr{G})} \biggr)^k
\, \Biggl[1 - \biggl(\frac{1+d}{n} \biggr)^k \Biggr] \nonumber \\[0.1cm]
\label{eq: 10.11.22c3}
& \geq \biggl( \frac{n}{\sqrt{1+d} \; \theta(\Gr{G})} \biggr)^k \,
\Biggl[1 - \biggl( \frac{1+d}{n} \biggr)^3 \Biggr] > 0,
\end{align}
which, by combining \eqref{eq: 10.11.22c2} and \eqref{eq: 10.11.22c3}, gives the stronger condition
\begin{eqnarray}
\label{eq: 10.11.22c4}
\frac{n^3-(1+d)^3}{n^3} \, \biggl( \frac{n}{\sqrt{1+d} \; \theta(\Gr{G})} \biggr)^k
\geq 2 + (1+d)^{-\frac{3}{2}},
\end{eqnarray}
with $k \geq 3$. Solving inequality \eqref{eq: 10.11.22c4} implies that
inequality \eqref{eq: 09.11.22a3} is satisfied for all $k \geq k_0$,
with the closed-form expression of $k_0$ in \eqref{eq: k_0}.
It therefore gives that if $\Gr{G}$ is a $d$-regular graph on $n$ vertices, which satisfies
the condition in \eqref{eq: 09.11.22a2}, then $\Gr{G}^{\boxtimes \, k}$ is non-Ramanujan for
all $k \geq k_0$ (it becomes, in fact, a highly non-Ramanujan graph since both sides of
inequality \eqref{eq: 09.11.22a3} have different exponential growth rates in $k$, so the condition
for a Ramanujan graph is sharply violated for the strong power $\Gr{G}^{\boxtimes \, k}$
when the value of $k$ is increased).

We next specialize this result for graphs that are self-complementary and vertex-transitive.
For $n=1$, the complete graph $G=\CoG{1}$ is a self-complementary and vertex-transitive
graph, whose all strong powers are also isomorphic to $\CoG{1}$, so they are therefore
non-Ramanujan graphs.

Let $\Gr{G}$ be a graph of order $n>1$ that is self-complementary and vertex-transitive, so it is
$d$-regular with $d = \tfrac12 (n-1)$.
Additionally, for such a graph $\Gr{G}$, the Lov\'{a}sz $\theta$-function is equal to $\theta(\Gr{G}) = \sqrt{n}$,
and it coincides with the Shannon capacity of $\Gr{G}$ (see \cite[Theorems~8 and~12]{Lovasz79_IT}). Then,
\begin{eqnarray}
\label{eq: 10.11.22a1}
\frac{n}{\sqrt{d+1}} = \sqrt{\frac{2n^2}{n+1}} > \sqrt{n} = \theta(\Gr{G}),
\end{eqnarray}
so the condition in \eqref{eq: 09.11.22b1} is fulfilled by graphs
of order $n$ that are self-complementary and vertex-transitive.
The value of $k_0$ in \eqref{eq: k_0} is specialized for such graphs to
\begin{align}
\label{eq: 10.11.22a2}
k_0 &= \max\left\{ 3, \left\lceil \frac{2 \log \Bigl( \frac{8n^3}{8n^3-(n+1)^3} \Bigr)
+ 2 \log \Bigl( 2 + \sqrt{\frac{8}{(n+1)^3}} \, \Bigr)}{\log \bigl( \frac{2n}{n+1} \bigr)} \right\rceil \right\} \\[0.2cm]
&=
\begin{dcases}
\label{eq: 11.11.22a1}
5  & \mbox{if} \; n = 5, \\
4  & \mbox{if} \; n = 9, \\
3  & \mbox{if} \; n \geq 13 \; \mbox{with} \; n \equiv 1 \, ( \hspace*{-0.35cm} \mod 4).
\end{dcases}
\end{align}
The constraint on $n$ in \eqref{eq: 11.11.22a1} is the
necessary condition on $n$ in Remark~\ref{remark: n's for sc-ve graphs}.

\subsection{Proofs for Section~\ref{subsection: LB on the chromatic number of strong products}}
\label{subsection: proofs - chromatic numbers of strong products}

\subsubsection{Proof of Proposition~\ref{prop. chromatic numbers}}
\label{subsubsection: proof of prop. - chromatic numbers of strong products}

\begin{enumerate}[(a)]
\item
Let $\Gr{G}_1, \ldots, \Gr{G}_k$ be $k$ simple, finite and undirected graphs,
$\bigcard{\V{\Gr{G}_\ell}} = n_\ell$ for $\ell \in \OneTo{k}$,
and let $\Gr{G} = \Gr{G}_1 \boxtimes \ldots \boxtimes \Gr{G}_k$.
We provide two alternative simple proofs of \eqref{eq:26.10.22a1}.

{\em First proof:}
\begin{IEEEeqnarray}{rCl}
\label{eq: 13.11.2022a1}
\chrnum{\Gr{G}} &\geq& \theta(\CGr{G}) \\[0.1cm]
\label{eq: 13.11.2022a2}
&\geq& \frac{\card{\V{\Gr{G}}}}{\theta(\Gr{G})} \\
\label{eq: 13.11.2022a3}
&=& \prod_{\ell=1}^k  \frac{\card{\V{\Gr{G}_\ell}}}{\theta(\Gr{G}_\ell)},
\end{IEEEeqnarray}
where \eqref{eq: 13.11.2022a1} holds by \eqref{eq1b: sandwich},
\eqref{eq: 13.11.2022a2} holds by \eqref{eq: Lovasz79 - Corollary 2},
and equality \eqref{eq: 13.11.2022a3} holds by \eqref{eq: Lovasz79 - Theorem 7}
and since $\V{\Gr{G}} = \V{\Gr{G}_1} \times \ldots \times \V{\Gr{G}_k}$.
The ceiling operation can be add to the right-hand side of
\eqref{eq: 13.11.2022a3} since a chromatic number is an integer.

{\em Second proof:}
\begin{IEEEeqnarray}{rCl}
\label{eq: 13.11.2022a4}
\chrnum{\Gr{G}} &\geq& \frac{\card{\V{\Gr{G}}}}{\indnum{\Gr{G}}} \\[0.1cm]
\label{eq: 13.11.2022a5}
&\geq& \frac{\card{\V{\Gr{G}}}}{\theta(\Gr{G})} \\
\label{eq: 13.11.2022a6}
&=& \prod_{\ell=1}^k  \frac{\card{\V{\Gr{G}_\ell}}}{\theta(\Gr{G}_\ell)},
\end{IEEEeqnarray}
where \eqref{eq: 13.11.2022a4} holds by \eqref{eq: well-known1},
\eqref{eq: 13.11.2022a5} holds by \eqref{eq1: capacity bounds},
and \eqref{eq: 13.11.2022a6} is \eqref{eq: 13.11.2022a3}.

We next prove \eqref{eq:26.10.22a2}.
\begin{IEEEeqnarray}{rCl}
\label{eq: 13.11.2022a7}
\chrnum{\CGr{G}} &\geq& \theta(\Gr{G}) \\
\label{eq: 13.11.2022a8}
&=& \prod_{\ell=1}^k \theta(\Gr{G}_\ell),
\end{IEEEeqnarray}
where \eqref{eq: 13.11.2022a7} holds by \eqref{eq1a: sandwich},
and \eqref{eq: 13.11.2022a8} holds by \eqref{eq: Lovasz79 - Theorem 7}.

\item Let $\Gr{G}_1, \ldots, \Gr{G}_k$ be regular graphs, where $\Gr{G}_\ell$ is $d_\ell$-regular
of order $n_\ell$ for all $\ell \in \OneTo{k}$. Inequality \eqref{eq:26.10.22a3} is \eqref{eq:26.10.22a1}.
Inequality \eqref{eq:26.10.22a4} follows from  \eqref{eq: Lovasz79 - Theorem 9} and \eqref{eq:26.10.22a3}.
Furthermore, by Item~(a) in Proposition~\ref{prop1: bounds on theta}, inequality \eqref{eq:26.10.22a4} holds
with equality if each regular graph $\Gr{G}_\ell$, for $\ell \in \OneTo{k}$, is either edge-transitive
or strongly regular.

\item By \eqref{eq: 13.11.2022a3}, with $\card{\V{\Gr{G}_\ell}} = n_\ell$,
\begin{IEEEeqnarray}{rCl}
\label{eq: 13.11.2022a9}
\frac{\card{\V{\Gr{G}}}}{\theta(\Gr{G})} &=& \prod_{\ell=1}^k  \frac{n_\ell}{\theta(\Gr{G}_\ell)}.
\end{IEEEeqnarray}
Suppose that, for all $\ell \in \OneTo{k}$, $\Gr{G}_\ell$ is $d_\ell$-regular, and it is also
either edge-transitive or strongly regular. By Item~(a) in
Proposition~\ref{prop1: bounds on theta}, for all $\ell \in \OneTo{k}$,
\begin{IEEEeqnarray}{rCl}
\label{eq: 13.11.2022a10}
\theta(\Gr{G}_\ell) = -\frac{n_\ell \, \Eigval{\min}{\Gr{G}_\ell}}{d_\ell -  \Eigval{\min}{\Gr{G}_\ell}}.
\end{IEEEeqnarray}
Combining \eqref{eq: 13.11.2022a9} and \eqref{eq: 13.11.2022a10} gives
\begin{IEEEeqnarray}{rCl}
\label{eq: 13.11.2022a11}
\frac{\card{\V{\Gr{G}}}}{\theta(\Gr{G})} &=& \prod_{\ell=1}^k
\biggl( 1 - \frac{d_\ell}{\Eigval{\min}{\Gr{G}_\ell}} \biggr).
\end{IEEEeqnarray}
On the other hand, since $\Gr{G} = \Gr{G}_1 \boxtimes \ldots \boxtimes \Gr{G}_\ell$
is $d$-regular, with $d \triangleq d(\Gr{G})$ as given in \eqref{eq:26.10.22a6}, it
follows from \eqref{eq: Lovasz79 - Theorem 9} that
\begin{IEEEeqnarray}{rCl}
\label{eq: 13.11.2022a12}
\theta(\Gr{G}) \leq -\frac{\card{\V{\Gr{G}}} \; \Eigval{\min}{\Gr{G}}}{d(\Gr{G}) - \Eigval{\min}{\Gr{G}}}.
\end{IEEEeqnarray}
It should be noted, in regard to \eqref{eq: 13.11.2022a12}, that even if all $\Gr{G}_\ell$'s
are regular and edge-transitive graphs, their strong product $\Gr{G}$ is not necessarily
edge-transitive. In fact, $\Gr{G}$ is not edge-transitive, unless all the $k$ factors
$\{\Gr{G}_\ell\}_{\ell=1}^k$ are complete graphs (see \cite[Theorem~3.1]{HammackIK16}).
For this reason, \eqref{eq: 13.11.2022a12} does not hold in general with equality (see
\cite[Theorem~9]{Lovasz79_IT}).
Finally, combing \eqref{eq: 13.11.2022a11} and \eqref{eq: 13.11.2022a12} gives inequality
\eqref{eq:26.10.22a5}.

\item
Let $\Gr{G}_1, \ldots, \Gr{G}_k$ be regular graphs, where $\Gr{G}_\ell$ is $d_\ell$-regular
on $n_\ell$ vertices for all $\ell \in \OneTo{k}$. Then, under the assumptions of Item~(d),
\begin{enumerate}[(1)]
\item
\begin{IEEEeqnarray}{rCl}
\label{eq: 14.11.22a1}
\prod_{\ell=1}^k  \frac{\card{\V{\Gr{G}_\ell}}}{\theta(\Gr{G}_\ell)}
& = & \prod_{\ell=1}^k  \theta(\CGr{G_\ell}) \\
\label{eq: 14.11.22a2}
& \geq & \prod_{\ell=1}^k \clnum{\Gr{G_\ell}}
\end{IEEEeqnarray}
where \eqref{eq: 14.11.22a1} holds since, by assumption,
each of the graphs $\Gr{G}_1, \ldots, \Gr{G}_k$ is vertex-transitive
or a strongly regular graph
(this is because \cite[Theorem~8]{Lovasz79_IT} and \eqref{eq: Lovasz equality for srg}
provide different sufficient conditions for inequality
\eqref{eq: Lovasz79 - Corollary 2} to hold with equality).
Inequality~\eqref{eq: 14.11.22a2} holds by the leftmost inequality
in \eqref{eq1b: sandwich}.

\item
\begin{IEEEeqnarray}{rCl}
\label{eq: 14.11.22a3}
\prod_{\ell=1}^k \biggl( 1 - \frac{d_\ell}{\Eigval{\min}{\Gr{G}_\ell}} \biggr)
&=& \prod_{\ell=1}^k  \frac{\card{\V{\Gr{G}_\ell}}}{\theta(\Gr{G}_\ell)} \\
\label{eq: 14.11.22a4}
& \geq & \prod_{\ell=1}^k \clnum{\Gr{G_\ell}}
\end{IEEEeqnarray}
where, by Item~(a) of Proposition~\ref{prop1: bounds on theta}, equality \eqref{eq: 14.11.22a3}
holds since (by assumption), for all $\ell \in \OneTo{k}$, the graph $\Gr{G}_\ell$ is either
regular and edge-transitive, or a strongly regular graph.
Inequality \eqref{eq: 14.11.22a4} holds under the same reasoning as of \eqref{eq: 14.11.22a1}
and \eqref{eq: 14.11.22a2}.
\end{enumerate}
To summarize, it shows that under proper assumptions, the lower bound on the chromatic
number of $\Gr{G}$ in the right-hand side of \eqref{eq:26.10.22a1}, or even its loosened
bound in the right-hand side of \eqref{eq:26.10.22a4}, are larger than or equal to the
lower bound $\overset{k}{\underset{\ell=1}{\prod}} \clnum{\Gr{G}_\ell}$.

\item
By \eqref{eq:26.10.22a2},
\begin{IEEEeqnarray}{rCl}
\label{eq: 14.11.22a5}
\chrnum{\CGr{G}} \geq \prod_{\ell=1}^k \theta(\Gr{G}_\ell).
\end{IEEEeqnarray}
Let, for all $\ell \in \OneTo{k}$, the graph $\Gr{G}_\ell$ be
$d_\ell$-regular on $n_\ell$ vertices, and suppose that it is
either edge-transitive or strongly regular.
Then, by Item~(a) of Proposition~\ref{prop1: bounds on theta},
\begin{IEEEeqnarray}{rCl}
\label{eq: 14.11.22a6}
\theta(\Gr{G}_\ell) = -\frac{n_\ell \, \Eigval{\min}{\Gr{G}_\ell}}{d_\ell
- \Eigval{\min}{\Gr{G}_\ell}}, \quad \forall \, \ell \in \OneTo{k}.
\end{IEEEeqnarray}
Combining \eqref{eq: 14.11.22a5} and \eqref{eq: 14.11.22a6},
followed by taking a ceiling operation on the lower bound on the chromatic
number $\chrnum{\CGr{G}}$, gives \eqref{eq:26.10.22a7}.

\item
By the assumption that $\Gr{G}_1, \ldots, \Gr{G}_k$ are self-complementary,
\begin{IEEEeqnarray}{rCl}
\label{eq: 14.11.22a7}
\theta(\Gr{G}_\ell) = \theta(\CGr{G}_\ell), \quad \forall \, \ell \in \OneTo{k}.
\end{IEEEeqnarray}
Furthermore, by the assumption that for all $\ell \in \OneTo{k}$,
$\Gr{G}_\ell$ is a graph on $n_\ell$ vertices that is either
vertex-transitive or strongly regular,
\begin{IEEEeqnarray}{rCl}
\label{eq: 14.11.22a8}
\theta(\Gr{G}_\ell) \, \theta(\CGr{G}_\ell) = n_\ell, \quad \forall \, \ell \in \OneTo{k}.
\end{IEEEeqnarray}
Combining \eqref{eq: 14.11.22a7} and \eqref{eq: 14.11.22a8} gives
\begin{IEEEeqnarray}{rCl}
\label{eq: 14.11.22a9}
\theta(\Gr{G}_\ell) = \sqrt{n_\ell}, \quad \forall \, \ell \in \OneTo{k}.
\end{IEEEeqnarray}
Consequently, by \eqref{eq:26.10.22a1} and \eqref{eq: 14.11.22a9},
\begin{IEEEeqnarray}{rCl}
\label{eq: 14.11.22a10}
\chrnum{\Gr{G}} &\geq& \Bigg\lceil \prod_{\ell=1}^k \frac{n_\ell}{\theta(\Gr{G}_\ell)} \Bigg\rceil \\
\label{eq: 14.11.22a11}
&=& \Big\lceil \prod_{\ell=1}^k \sqrt{n_\ell} \, \Big\rceil \\
\label{eq: 14.11.22a12}
&=& \big\lceil \sqrt{n} \, \big\rceil,
\end{IEEEeqnarray}
and, from \eqref{eq:26.10.22a2} and \eqref{eq: 14.11.22a9},
\begin{IEEEeqnarray}{rCl}
\label{eq: 14.11.22a13}
\chrnum{\CGr{G}} &\geq& \Bigg\lceil \prod_{\ell=1}^k \theta(\Gr{G}_\ell) \Biggr\rceil \\
\label{eq: 14.11.22a14}
&=& \Big\lceil \prod_{\ell=1}^k \sqrt{n_\ell} \, \Big\rceil \\
\label{eq: 14.11.22a15}
&=& \big\lceil \sqrt{n} \, \big\rceil.
\end{IEEEeqnarray}
This proves \eqref{eq:26.10.22a8.1} and \eqref{eq:26.10.22a8.2},
and it completes the proof of Proposition~\ref{prop. chromatic numbers}.
\end{enumerate}

\subsubsection{Proof of Corollary~\ref{cor8: chrnum strong product srg}}
\label{subsubsection: proof of cor. chr. strong prod. srg}

The rightmost inequality in \eqref{eq:02.11.2022b1} is a well-known upper bound
on the chromatic number of strong products (see \cite{Berge73}, \cite[Theorem~3]{Borowiecki72}).
The leftmost inequality in \eqref{eq:02.11.2022b1} gives a lower bound on the chromatic
number of a strong product of (not necessarily distinct) non-complete, and
strongly regular graphs. It readily follows by combining equality~\eqref{eq:30.10.22a1}
in Corollary~\ref{cor4: Lovasz number for srg}, together with
inequality~\eqref{eq:26.10.22a1} in Proposition~\ref{prop. chromatic numbers}.
Finally, by Part~1 of Item~(d) in Proposition~\ref{prop. chromatic numbers},
the leftmost term in \eqref{eq:02.11.2022b1} is larger than or equal to the product
of the clique numbers of $\{\Gr{G}_\ell\}_{\ell=1}^k$.

\section*{Acknowledgments}
The author wishes to acknowledge the two anonymous reviewers for
helpful and timely reports.

\end{document}